\documentclass[10pt]{amsart}
\usepackage{a4wide}
\usepackage{amsmath,amssymb,url,stmaryrd,enumerate}
\usepackage[all]{xy}
\usepackage{srcltx}
\usepackage{color}

\def\La{\Lambda}

\def\O{\mathcal{O}} 
\def\F{\mathcal{F}}

\def\*{^\times }
\def\G{\mathcal{G}}

\def\g{\gamma}
\def\w{\omega}
\def\a{\alpha}
\def\b{\beta}

\def\ph{\varphi}
\def\Q{\mathbb{Q}}
\def\Qp{\mathbb{Q}_p}

\def\Z{\mathbb{Z}}
\def\N{\mathbb{N}}
\def\Hom{\text{\rm Hom}}
\def\hom{\mathcal{H}om}
\def\spec{\text{Spec}}

\def\limp{\underset{\longleftarrow}{\text{ lim }}\!}
\def\limi{\underset{\longrightarrow}{\text{ lim }}\!}

\def\mod{{\rm mod}\ }

\def\im{\text{im}}

\def\u{\underline}

\def\coker{\mathrm{coker}\ \!}
\def\im{\mathrm{im}}

\def\fp{\mathfrak{p}}
\def\cA{\mathcal{A}}
\def\fa{\mathfrak{a}}
\def\fb{\mathfrak{b}}

\def\la {\leftarrow}

\def\cC{\mathcal{C}}	
\def\cB{\mathcal{B}}
\def\cD{\mathcal{D}}
\def\cF{\mathcal{F}}
\def\cG{\mathcal{G}}

\def\cK{\mathcal{K}}
\def\cI{\mathcal{I}}

\def\cL{\mathcal{L}}
\def\cM{\mathcal{M}}
\def\cN{\mathcal{N}}

\def\wt{\widetilde}
\def\wh{\widehat}

\def\ov{\overline}

\def\lra{\longrightarrow}

\def\fil{{\rm Fil}}

\def\pr{{\noindent{\it Proof} :} }

\def\da{\downarrow}

  \def\rp{{\rm Rep}}

 \def\rptor{{\rp^{\rm tor}_{B^+_{dR}}(G_K)}}
  \def\rpdr{{\rp_{B_{dR}}(G_K)}}
  \def\rpdrp{{\rp_{B^+_{dR}}(G_K)}}

\def\co{\cC^0(G_K)}
\def\ci{\cC^{\infty}(G_K)}

\def\oxm{{Coh(\O_X)}}

\def\hra{\hookrightarrow}

 \def\Ob{{\rm Ob}}

\newtheorem{theo}{Theorem}[section]

\newtheorem{thm}[theo]{Theorem}

\newtheorem{prop}[theo]{Proposition}

\newtheorem{lemma}[theo]{Lemma}
\newtheorem{coro}[theo]{Corollaire}
\newtheorem{cor}[theo]{Corollary}

\newtheorem*{thA}{Theorem A}

\newtheorem{rem}[theo]{Remarque}

\newtheorem{remk}[theo]{Remark}
\newtheorem{remks}[theo]{Remarks}

\title{Almost $C_p$ Galois representations and vector bundles}
\author{Jean-Marc Fontaine}

\date{\today}

\begin{document}

\maketitle

\noindent Abstract: {\it Let $K$ be a finite extension of $\Q_p$ and $G_K$ the absolute Galois group. Then $G_K$ acts on the fundamental curve $X$ of $p$-adic Hodge theory (cf. \cite{FF}) and we may consider the abelian category $\cM(G_K)$ of coherent $\O_X$-modules equipped with a continuous and semi-linear action of $G_K$. 

An {\rm almost $C_p$-representation of $G_K$}  is a $p$-adic Banach space $V$ equipped with a linear and continuous action of $G_K$ such that there exists $d\in\N$, two $G_K$-stable finite dimensional sub-$\Q_p$-vector spaces $U_+$ of $V$, $U_-$ of $C_p^d$, and a $G_K$-equivariant isomorphism
$$V/U_+\lra C_p^d/U_-\ .$$
These representations form an abelian category $\cC(G_K)$ (cf. \cite{Fo03}). The main purpose of this paper is to prove that $\cC(G_K)$ can be recovered from $\cM(G_K)$ by a simple construction (and conversely) inducing, in particular,  an equivalence of triangulated categories
$$D^b(\cM(G_K))\lra D^b(\cC(G_K))$$}

\section{Introduction}

\subsection{} We fix a prime number $p$, an algebraic closure $\ov\Q_p$ of $\Q_p$ and a finite extension $K$ of $\Q_p$ contained in $\ov\Q_p$. We set $G_K={\rm Gal}(\ov\Q_p/K)$ and $C_p$ the $p$-adic completion of $\ov\Q_p$ on which $G_K$ acts by continuity.

In \cite{FF}, Laurent Fargues and I introduced the fundamental curve $X_{\Q_p,C_p^\flat}$ of $p$-adic Hodge theory, denoted $X$ below, which is a separated noetherian regular scheme of dimension $1$ defined over $\Q_p$ (i.e. we have $H^0(X,\O_X)=\Q_p$). The structural sheaf is naturally equipped with a topology (cf. \S \ref{top}): if $U$ is any open subset of $X$, then $\O_X(U)$ is a locally convex $\Q_p$-algebra.  There is a natural action of $G_K$ on $X$ which is continuous. We may consider the abelian category $\cM(G_K)$ of {\it $G_K$-equivariant coherent $\O_X$-modules}, that is of coherent $\O_X$-modules equipped with a semi-linear and continuous action of $G_K$. 

Any non zero $\cF\in\Ob(\cM(G_K))$ as a degree $\deg(\cF)\in\Z$ and a rank ${\rm rk}(\cF)\in\N$, hence also a slope $s(\cF)=\deg(\cF)/{\rm rk}(\cF)\in \Q\cup\{+\infty\}$ (with the convention that $s(\cF)=+\infty$ if $\cF$ is a torsion $\O_X$-module). As in the classical case, one says that a coherent $\O_X[G_K]$-module $\cF$ is {\it semi-stable} if $\cF\not=0$ and if $s(\cF')\leq s(\cF)$ for any non zero subobject $\cF'$ of $\cF$.

 We may consider the full sub-category $\cM^0(G_K)$ of $\cM(G_K)$ whose objects are semi-stable of slope $0$. One of the main results of \cite{FF} is that, if $\cF$ is any object of $\cM^0(G_K)$, then $\cF(X)=H^0(X,\cF)$ is a finite dimensional $\Q_p$-vector space, hence is an object of the abelian category $\rp_{\Q_p}(G_K)$ of $p$-adic representations of $G_K$ (that is of finite dimensional $\Q_p$-vector spaces equipped with a linear and continuous action of $G_K$) and that the functor
$$\cM^0(G_K)\lra \rp_{\Q_p}(G_K)\ \ ,\ \ \cF\mapsto \cF(X)$$
is an equivalence of categories (with $V\mapsto \O_X\otimes_{\Q_p}V$ as a quasi-inverse).

The main purpose of this paper is to discuss the following question: Is there an extension of this result enabling us to give an analogous Galois description of all objects of $\cM(G_K)$ ?

\subsection{}\label{introcgk} In \cite{Fo03}, I introduced the category of {\it almost $C_p$-representations of $G_K$}: A {\it Banach representation of $G_K$} is a $p$-adic Banach space (i.e. a topological $\Q_p$-vector space whose topology can be defined by a norm and which is complete) equipped with a linear and continuous action of $G_K$. With an obvious definition of morphisms, Banach representations of $G_K$ form an additive category $\cB(G_K)$ containing the category $\rp_{\Q_p}(G_K)$ as a full subcategory. By continuity, $G_K$ acts on the $p$-adic completion $C_p$ of $\ov\Q_p$ and $C_p$ has a natural structure of a Banach representation. The category $\cC(G_K)$ of {\it almost $C_p$-representations of $G_K$} is the full subcategory of $\cB(G_K)$ whose objects are those $V$'s for which one can find $d\in\N$, two $G_K$-stable finite dimensional sub-$\Q_p$-vector spaces $U_+$ of $V$ and $U_-$ of $C_p^d$ and an isomorphism $V/U_+\to C_p^d/U_-$ in $\cB(G_K)$. This category turns out to be abelian ({\it loc.cit.}).\\

The curve $X$ has only one closed point $\infty$ which is $G_K$-stable and the orbit under $G_K$ of any other closed point is infinite. This implies that a torsion object of $\cM(G_K)$ is supported at $\infty$. As the completion of $\O_{X,\infty}$ is the ring $B^+_{dR}$ of $p$-adic periods, the category $\cM^\infty(G_K)$ of torsion objects of $\cM(G_K)$ ($\iff$ semi-stable objects of slope $\infty$) can be identified to the category $\rp_{B^+_{dR}}^{\rm tor}(G_K)$ of $B^+_{dR}$-modules of finite length equipped with a semi-linear and continuous action of $G_K$. The topology of any $B^+_{dR}$-module of finite length is the topology of a $p$-adic Banach space and we may consider the forget functor
$$\rp_{B^+_{dR}}^{\rm tor}(G_K)\lra \cB(G_K)\ .$$
We proved in {\it loc.cit.} that this functor is fully faithful and that the essential image $\cC^\infty(G_K)$ is contained in $\cC(G_K)$. Hence, setting $\cC^0(G_K)=\rp_{\Q_p}(G_K)$, we see that for $s\in\{0,\infty\}$, the functor
$\cF\mapsto\cF(X)$ induces an equivalence of categories
$$\cM^s(G_K)\lra \cC^s(G_K)\ .$$
\
Similarly, as for a smooth projective curve over a field, we defined in \cite{FF} the Harder-Narasimhan filtration of any $\cF \in\cM(G_K)$: this  is the unique filtration
$$0=\cF^0\subset\cF^1\subset\ldots\subset\cF^{r-1}\subset\cF^{r}=\cF$$
such that all the $\cF^i/\cF^{i-1}$ are semi-stable and that $s(\cF^i/\cF^{i-1})>s(\cF^{i+1}/\cF^i)$ for $0<i<r$. We call the  $s(\cF^i/\cF^{i-1})$, for $1\leq i\leq r$, the {\it HN-slopes of $\cF$}.

Let $\cM^{\geq 0}(G_K)$ the full subcategory of $\cM(G_K)$ whose objects are {\it effective}, i.e. such that all their HN-slopes are $\geq 0$.

Similarly let $\cC^{\geq 0}(G_K)$ the full subcategory of $\cC(G_K)$ whose objects are {\it effective}, i.e. those $V$'s which are isomorphic to a subobject (in $\cC(G_K)$) of an object of $\cC^\infty(G_K)$.

If $\cF$ is any coherent $\O_X[G_K]$-module, then $\cF(X)$ is a topological $\Q_p$-vector space equipped with a linear and continuous action of $G_K$. The main result of this paper is:

\begin{thA} \label{mainthrm} If $\cF$ is any coherent $\O_X[G_K]$-module, then $\cF(X)$ is an effective almost $C_p$-representation of $G_K$. By restriction to $\cM^{\geq 0}(G_K)$ the functor $\cF\mapsto\cF(X)$ induces an equivalence of categories
$$\cM^{\geq 0}(G_K)\lra \cC^{\geq 0}(G_K)\ .$$
\end{thA}

This equivalence doesn't extend to an equivalence between $\cM(G_K)$ and $\cC(G_K)$. Nevertheless each of these two categories can be reconstructed from the other: The above functor induces an equivalence of triangulated categories
$$D^b(\cM(G_K))\lra D^b(\cC(G_K))$$
and each of them can be reconstructed as the heart of a $t$-structure. More precisely:

-- Denote $\cM^{<0}(G_K)$ the full sub-category of $\cM(G_K)$ whose objects are those for which all HN-slopes are $<0$. Then $t=(\cM^{\geq 0}(G_K),\cM^{<0}(G_K))$ is what is called a torsion pair on $\cM(G_K)$. From this torsion pair, we can construct an other abelian category f$(\cM(G_K))^t$ which is the full subcategory of $D^b(\cM(G_K))$ whose objects are those $\cF^\bullet$'s such that
$$\cF^i=0 \text{ for } i\not\in\{0,1\}\ \ ,\ \ H^0(\cF^\bullet) \text{ is an object of } \cM^{<0}(G_K)\ \ ,\ \  H^1(\cF^\bullet) \text{ is an object of } \cM^{\geq 0}(G_K)\ .$$
There is a natural equivalence $(\cM(G_K))^t\to\cC(G_K)$.

-- Similarly, denote $\cC^{<0}(G_K)$ the full sub-category of $\cC(G_K)$ whose objects are those $V$'s for which $\Hom(V,W)=0$ for all $W$ in $cC^\infty(G_K)$. Then $t'=(\cC^{<0}(G_K),\cC^{\geq 0}(G_K))$ is a torsion pair on $\cC(G_K)$ which can be used to define the abelian  subcategory $(\cC(G_K))^{t'}$ which is the full subcategory of $D^b(\cC(G_K)$ whose objects are those $V^\bullet$ such that
$$V^i=0 \text{ for } i\not\in\{0,1\}\ \ ,\ \ H^0(V^\bullet) \text{ is an object of } \cC^{\geq 0}(G_K)\ \ ,\ \  H^1(V^\bullet) \text{ is an object of } \cC^{< 0}(G_K) \ .$$  
There is a natural equivalence $(\cC(G_K))^{t'}\to\cM(G_K)$.\\

\subsection{Contents}\

In \S 2, we recall and slightly extend the results of \cite{Fo03} on almost $C_p$-representations. We first recall (\S 2.1) some basic facts about locally convex spaces over a non archimedean field. We introduce (\S 2.2) the category of ($p$-adic) ind-Fr\'echet representations (of $G_K$). Then (\S 2.3), we recall some basic facts about the ring of periods $B^+_{dR}$ and $B_{dR}$ that we equip with a locally convex topology. In \S 2.4, we discuss some properties of $B^+_{dR}$-representations and  $B_{dR}$-representations (of $G_K)$.

 We describe (\S 2.5)  the main properties of the category $\cC(G_K)$ of almost $C_p$-representations and of its full subcategories $\cC^0(G_K)$ of finite dimensional $p$-adic representations and $\cC^\infty(G_K)$ of $B^+_{dR}$-representations of finite length. In \S 2.5, we also introduce  the category $\wh\cC(G_K)$ of representations of $G_K$ which are suitable limits (in the category of locally convex $p$-adic representations of $G_K$) of almost $C_p$-representations.
In \S 2.6, we recall the notion of almost split exact sequence of $\cB(G_K)$ and the fact that an extension in $\cB(G_K)$ of two almost $C_p$-representations is an almost $C_p$-representation if and only if the associated short exact sequence almost splits.\\

 In \S 3, we study the category $\rp_{B_e}(G_K)$ of {\it $B_e$-representations of $G_K$}. We also recall and precise some of the results of \cite{FF} on {\it coherent $\O_X[G_K]$-modules}.  We first recall (\S 3.1) some basic facts about the sub-$\Q_p$-algebras $B^+_{cris}$ and $B_e$ of $B_{dR}$ which are stable under the action of $G_K$ and equipped with a natural topology of locally convex algebras. Then we introduce (\S 3.2) $\rp_{B_e}(G_K)$  and show that this is a $\Q_p$-linear abelian category.
 
 We recall (\S 3.3) the definition of the fundamental curve $X=X_{\Q_p,C_p^\flat}$ of $p$-adic Hodge theory introduced in \cite{FF} on which $G_K$ acts and give a description of the category  $\oxm$ of coherent $\O_X$-modules. We discuss (\S 3.4) the topology on the structural sheaf $\O_X$ and give a description of the category $\cM(G_K)$ of coherent $\O_X[G_K]$-modules (\S 3.5). We describe (\S 3.6) the Harder-Narasimhan filtration on any $\cF\in\cM(G_K)$.
 
 We consider two full subcategories  of $\cM(G_K)$:
 
 -- the category $\cM^0(G_K)$, of the semi-stable objects of slope $0$,
 
 -- the category $\cM^\infty(G_K)$ of objects whose underlying $\O_X$-module is torsion.
 
 We show (\S 3.7) that the global sections functor induces equivalence of categories 
 $$\cM^0(G_K)\lra\cC^0(G_K)\ \ \text{ and }\ \ \cM^\infty(G_K)\lra \cC^\infty(G_K)\ .$$
 
In \S 3.8, we introduce two kinds of twists of the objects of $\cM(G_K)$, the {\it Tate twists} and the {\it Harder-Narasimhan twists}.

Say that a $B_e$-representation $\La$ is {\it  trivialisable} if there exists $U\in\cC^0(G_K)$ and an isomorphism $B_e\otimes_{\Q_p}U\to\La$.  In \S 3.9, we show that $\rp_{B_e}(G_K)$ is the smallest subcategory of itself containing trivialisable $B_e$-representations and stable under taking extensions and direct summands.

In \S 3.10, we show that, if $\La$ is a $B_e$-representation of $G_K$, then the underlying topological $\Q_p$-vector space equipped with its action of $G_K$ is an object of $\wh\cC(G_K)$ and that the forget functor
$$\rp_{B_e}(G_K)\lra \wh\cC(G_K)$$
 is exact and fully faithful.
 
 We conclude this section by discussing the cohomology of coherent $\O_X$-modules (\S 3.11) and of coherent $\O_X[G_K]$-modules (\S 3.12). We show that, taking the global sections, we get a functor
 $$\cM(G_K)\lra\cC(G_K)\ \ ,\ \ \cF\mapsto \cF(X)=H^0(X,\O_X)$$ whose essential image is contained in  $\cC^{\geq 0}(G_K)$.\\
 
 The aim of \S 4 is to construct a left adjoint
 $$\cC(G_K)\lra \cM(G_K)\ \ ,\ \ V\mapsto \cF_V$$
 of the functor $\cF\mapsto\cF(X)$.
 
 We show (\S 4.3) that any almost $C_p$-representation $V$ has a {\it $B_e$-hull}, i.e. there is a pair $V_e=(V_e,\iota_e^V)$ with $V_e$ a $B_e$-representation (of $G_K$) and $\iota_e^V:V\to V_e$ a morphism in $\wh\cC(G_K)$ such that, for all $\La\in\rp_{B_e}(G_K)$, the map
 $$\Hom_{\rp_{B_e}(G_K)}(V_e,\La)\lra \Hom_{\wh\cC(G_K)}(V,\La)$$
 induced by $\iota_e^V$ is bijective. 
 
 Similarly with obvious definitions, we show that $V$ has a $B^+_{dR}$-hull $V^+_{dR}$ and a $B_{dR}$-hull $V_{dR}$.
 
 Using the existence of these hulls and the relations between them and knowing the description of $\cM(G_K)$ given in \S 3.5, the construction of the functor $V\mapsto \cF_V$ is quite simple.
 
 The proof of the existence of these hulls relies heavily on the description of all extensions in $\cC(G_K)$ of an object of $\cC^\infty(G_K)$ by an object of $\cC^0(G_K)$ which is given in \S 4.2.\\
 
 The aim of \S 5 is to prove theorem A. 
 
 We first show (\S 5.1 (resp. \S 5.2)) that $\cM^{\geq 0}(G_K)$ (resp. $\cC^{\geq 0}(G_K)$) is the smallest full subcategory of $\cM(G_K)$ (resp. $\cC(G_K)$) containing $\cM^0(G_K)$ and $\cM^\infty(G_K)$ (resp. $\cC^0(G_K)$ and $\cC^\infty(G_K)$) and stable under extensions and direct summands.
 
 Then (\S 5.3), we prove by d\'evissage that the functor
 $$\cM^{\geq 0}(G_K)\lra \cC^{\geq 0}(G_K)\ \ ,\ \ \cF\mapsto \cF(X)$$
 is an equivalence of exact categories (see \S \ref{conv}), the functor $V\mapsto \cF_V$ being a quasi-inverse.\\
 
 The purpose of \S 6 is to extend the main result to the categories $\cM(G_K)$ and $\cC(G_K)$. 
 
After some general nonsense on derived categories of exact subcategories of abelian categories (\S 6.1), we first extend the main result to an equivalence of of triangulated categories (\S 6.2),
$$D^b(\cM(G_K))\lra D^b(\cC(G_K))\ .$$

To go further, we need to introduce the full subcategories $\cM^{<0}(G_K)$ of $\cM(G_K)$ and $\cC^{<0}(G_K)$ of $\cC(G_K)$ of coeffective objects. The main theorem said that, if $\cF\in\cM^{\geq 0}(G_K)$, then $H^0(X,\cF)$ has a natural structure of an object of $\cC^{\geq 0}(G_K)$ and this structure determines $\cF$. We prove (\S 6.3) that, if $\cF\in\cM^{<0}(G_K)$, then $H^1(X,\cF)$ has a natural structure of an object of $\cC^{< 0}(G_K)$ and this structure determines $\cF$.

Using this result, we can build $\cC(G_K)$ from $\cM(G_K)$ and conversely. We give two different recipes (with independent proofs) for that. We first (\S 6.4) describe explicitely the heart of the $t$-structure on $D^b(\cM(G_K))$ corresponding to $\cC(G_K)$ and of the $t$-structure on $D^b(\cC(G_K))$ corresponding to $\cM(G_K))$. In \S 6.5, we explain that $(\cM^{\geq 0}(G_K),\cM^{<0}(G_K))$ is a torsion pair on $\cM(G_K)$. One can use it to construct a new abelian category equipped with a torsion pair. Up to equivalence, it is $\cC(G_K)$ equipped with the torsion pair $(\cC^{<0}(G_K),\cC^{\geq 0}(G_K))$.\\

\subsection{A remark on possible generalisations}\

The results of this paper are obviously a special case of a much more general result where $K$ is replaced by any reasonable rigid analytic, Berkovich or adic space. Let's sketch a description of the case where $K$ is now any field complete with respect to a non archimedean non trivial absolute value with perfect residue field of characteristic $p$. 

-- We can define the abelian category $Coh(\O_{X_K})$ of {\it coherent modules on the curve $X_K$}. When $K$ is a perfectoid field, $X_K$ is the curve $X_{\Q_p,K^\flat}$ constructed in \cite{FF}. If $K$ is not perfectoid, then $X_K$ doesn't exist but one can define the category of coherent modules over this virtual curve. When $K$ is a finite extension of $\Q_p$, there is a natural equivalence of categories 
$$Coh(\O_{X_K})\lra \cM(G_K)\ .$$

-- We still have the Harder-Narasimhan filtration on $Coh(\O_{X_K})$ and may consider its exact subcategories $Coh^{\geq 0}(\O_{X_K})$ and $Coh^{<0}(\O_{X_K})$ which form a torsion pair $t$ on $Coh(\O_{X_K})$.

-- The construction of the curve $X_K$ is functorial in $K$. If $C$ is the completion of a separable closure $K^s$ of $K$, for any coherent $\O_{X_K}$-module $\cF$, we may consider the pull-back $f^\ast\cF$
of $\cF$ via $f:X_C\to X_K$. 

If $G_K={\rm Gal}(K^s/K)$, we may consider the exact category $\cB(G_K)$ of $p$-adic Banach representations of $G_K$ and we have exact and faithful functors
$$\begin{matrix} Coh^{\geq 0}(\O_{X_K})\lra \cB(G_K)&,&\cF\mapsto H^0(X_C,f^\ast\cF)\cr Coh^{<0}(\O_{X_K})\lra \cB(G_K)&,&\cF\mapsto H^1(X_C,f^\ast\cF)\end{matrix}$$

	But, in general, these functors are not fully faithful. Working with $\cB(G_K)$ amounts to work over the small pro-\'etale site of $K$ and we need to work with a bigger site. A possibility is to use the big pro-\'etale site $K_{pro\acute{e}t}$ of $K$ as defined by Scholze (\cite{Scho2}), \S 8 \footnote{More precisely, we fix an uncountable cardinal $\kappa$ satisfying the properties of \cite{Scho2}, lemma 4.1. The underlying category is the category of perfectoid spaces over $K$ which are $\kappa$-small ({\it loc. cit.}, def.4.3) and covering are as defined in {\it loc. cit.}, def. 8.1 (the only difference with the big pro-\'etale site of Scholze is that we restrict ourself to perfectoid spaces lying over the given non archimedean field $K$).} and to replace
	
	-- $\cB(G_K)$ with the category $Vect_{\Qp}(K)$ of $\Q_p$-sheaves over  $K_{pro\acute{e}t}$ ;
	
	-- and $\cC(G_K)$ with  the category of {\it pseudo-geometric $\Q_p$-sheaves}, an abelian full subcategory of $Vect_{\Q_p}(K)$ defined by imitating the definition of $\cC(G_K)$ as a full subcategory of $\cB(G_K)$.
	
	The correspondence $K\mapsto X_K$ can be extended to a functor 
$$U\mapsto X_U$$
 from the category of perfectoid spaces to the category of $\Q_p$-schemes. We also have exact and faithful functors
$$\begin{matrix}
Coh^{\geq 0}(\O_{X_K})\lra Vect_{\Q_p}(K)&,&\cF\mapsto \big{(} U\mapsto H^0(X_U,f_U^\ast\cF)\big{)}\cr
Coh^{<0}(\O_{X_K})\lra Vect_{\Q_p}(K)&,&\cF\mapsto \big{(} U\mapsto H^1(X_U,f_U^\ast\cF)\big{)}\end{matrix}$$
(where $f_U:X_U\to X_K$ is the structural morphism).

It seems likely (and not so hard to prove) that these functors are fully faithful and that one can describe their essential images $Vect_{\Q_p}^{pg,\geq 0}(K)$ and $Vect_{\Q_p}^{pg,<0}(K)$. These two functors seem to induce an equivalence of categories
$$(Coh(\O_{X_K}))^t\lra Vect_{\Q_p}^{pg}(K)$$
the induced torsion pair on $Vect_{\Q_p}^{pg}(K)$ being $t'=(Vect_{\Q_p}^{<0}(K),Vect_{\Q_p}^{\geq 0}(K))$.

In the case where $K$ is the $p$-adic completion of an algebraic closure of $\Q_p$, this result has been proved by Le Bras \cite{LB}. We hope to come back soon on this generalisation.

\subsection{Conventions and notations}\label{conv}\

 If $\cC$ is a category, we often write $C\in\cC$ for $C\in{\rm Ob}(\cC)$.\\
 
 An {\it exact subcategory of an abelian category $\cA$} is a strictly full subcategory of $\cA$ containing $0$ and stable under extensions.

If $\cB$ is an exact subcategory of $\cA$, we say that a sequence of morphisms of $\cA$ is {\it exact} if it is exact as a sequence of morphisms in $\cA$. In particular, we have the obvious notion of a {\it short exact sequence}. It is easy to see that, equipped with this class of short exact sequences, $\cB$ is an exact category in the  sense of Quillen (cf. \cite{Qu}, see also \cite{Lau}). Actually, any exact category $\cB$ in the sense of Quillen can be viewed as an exact subcategory of an abelian category (cf \cite{Qu}, \S2).\\

 As usual $\Z_p(1)$ is the Tate module of the multiplicative group, and, for all $n\in\N$,
$$\Z_p(n)={\rm Sym}_{\Z_p}^n\Z_p(1)\ \,\ \ \Z_p(-n)=\cL_{\Z_p}(\Z_p(n),\Z_p)\ .$$
If $M$ is any $\Z_p$-module equipped with a linear action of $G_K$, for all $n\in\Z$,
$$M(n)=M\otimes_{\Z_p}\Z_p(n)\ .$$

\subsection{} I would like to thank Laurent Fargues for helpful discussions. 

\section{Representations of $G_K$}

In this paper, each time we say "representation", we mean "representation of $G_K$". In this section, we introduce a few categories of such representations and describe some of their properties. Most of them are already known (see in particular \cite{Fo03}) or easy consequences of known properties.

\subsection{Banach, Fr\'echet, ind-Banach and ind-Fr\'echet}\

We refer to \cite{Emerton} and \cite{Schneider} for basic facts about $p$-adic functional analysis. All results of this paragraph are either contained or easy consequences of results contained
 in at least one of these two memoirs.
 
 We fix a non archimedean field $E$, i.e. a field complete with respect to a non trivial non archimedean absolute value,  and denote $\O_E$ its valuation ring. In the applications in this paper,   $E$ will be $\Q_p$.\\
 
 \begin{itemize}
\item A {\it locally convex $E$-vector space} is a topological $E$ vector space $V$ such that the open sub-$\O_E$-modules of $V$ form a fundamental system of neighbourhood of $0$,
\item A {\it Fr\'echet $E$-vector space} or an {\it $E$-Fr\'echet} is a locally convex $E$-vector space which is metrisable and complete,
\item A {\it Banach $E$-vector space} or an {\it $E$-Banach} is a Fr\'echet vector space whose topology can be defined by a norm,
\item An {\it ind-Fr\'echet} (resp. {\it ind-Banach}) $E$-vector space or an {\it ind-$E$-Fr\'echet} (resp. {\it ind-$E$-Banach}) is a locally convex $E$-vector space $V$, such that one can find an increasing sequence $(V_n)_{n\in\N}$ of closed sub-$E$-vector spaces such that 

i)  $V=\cup_{n\in\N}V_n$,

ii) each $V_n$, with the induced topology, is an $E$-Fr\'echet (resp. an $E$-Banach), 

iii) the topology of $V$ is the coarsest locally convex topology with these properties.
\end{itemize}

Condition (iii) is equivalent to the fact that a sub-$\O_E$-module $L$ of $V$ is open if and only if $L\cap V_n$ is open in $V_n$ for all $n\in\N$.\\

If $V$ is a topological $E$-vector space, it is an $E$-Fr\'echet if and only if

i) its topology can be defined by a countable family $(q_n)_{n\in\N}$ of semi-norms,

ii) $V$ is complete.

In this situation, replacing each $q_n$ by $q'_n=\sup_{0\leq i\leq n} q_i$, we may assume that $q_n\leq q_{n+1}$ for all $n$. Then, if $\ov V_n$ is the Hausdorff completion of $V$, with respect to $q_n$, this is an $E$-Banach and we have an homeomorphism
$$V\mapsto \limp_{n\in\N}\ov V_n$$
(with the inverse limit topology on the RHS). Conversely, any inverse limit, indexed by $\N$, of $E$-Banach is an $E$-Fr\'echet.\\

 Let $V$ be a topological  $E$-vector space. We say that a decreasing filtration $(F^nV)_{n\in\Z}$ by closed sub-$E$-vector spaces of $V$ is {\it admissible} if
 
 i) $\cup_{n\in\Z}F^nV=V$ and $\cap_{n\in\Z}F^nV=0$,
 
 ii) if $m\in\Z$ and $r\in\N$, then $F^mV/F^{m+r}V$, equipped with the induced topology, is an $E$-Banach,
 
 iii) if $m\in\Z$, then the natural map
 $$F^mV\lra \limp_{r\in\N} F^mV/F^{m+r}V$$
 is an homeomorphism (with the inverse limit topology on the RHS),
 
iv)  a sub-$\O_E$-module $L$ of $V$ is open if and only if $L\cap F^nV$ is open in $F^nV$ for all $n$.\\

The following result is obvious:

\begin{prop} Let $V$ be a topological $E$-vector space. Then

i)  $V$ is an ind-$E$-Fr\'echet if and only if it has an admissible filtration,

ii) $V$ is an $E$-Banach (resp. an $E$-Fr\'echet, resp. an ind-$E$-Banach) if and only if has an admissible filtration $(F^nV)_{n\in\Z}$ such that $F^0V=V$ and $F^1V=0$ (resp. $F^0V=V$, resp. $F^1V=0$).
\end{prop}

\begin{prop}\label{filadm} Let $V_1$ and $V_2$ two ind-$E$-Fr\'echet, $(F^nV_1)_{n\in\Z}$ an admissible filtration of $V_1$ and $(F^nV_2)_{n\in\Z}$ an admissible filtration of $V_2$. Let $u:V_1\to V_2$ an $E$-linear map. The following are equivalent:

i) the map $u$ is continuous,

ii) for all $m\in\Z$, there exists $n\in\Z$ such that $u(F^mV_1)\subset F^nV_2$ and the induced map 
$$F^mV_1\to F^nV_2$$ is continuous.
\end{prop}

\pr ii)$\implies$i): It's enough to show that, if $L$ is an open lattice in $V_2$, then $f^{-1}(L)$ is open in $V_1$ which means that if $m\in\Z$, then $f^{-1}(L)\cap F^mV_1$ is open in $F^mV_1$ which is indeed true as, if $n$ is such that $f(F^mV_1)\subset F^nV_2$, this is the inverse image of the continuous map $F^mV_1\to F^nV_2$ which is induced by $f$.

i)$\implies$ii): All the $F^nV_2$ are $E$-Fr\'echet. For each fixed $m$, so is $F^mV_1$ and the existence of such an $n$ is explained by Schneider in \cite{Schneider} (cor. 8.9). \qed

\begin{coro} Let $V$ be an ind-$E$-Fr\'echet  and $(F^nV)_{n\in\Z}$ an admissible filtration. Then $V$ is an $E$- Banach (resp. an $E$-Fr\'echet, resp. an ind-$E$-Banach) if and only if  there exists $m\leq n$ such that $F^mV=V$ and $F^nV=0$ (resp. $m$ such that $F^mV=V$, resp. $n$ such that $F^nV=0$).
\end{coro}

\begin{coro} Let $V$ be an ind-$E$-Fr\'echet and $(F_1^nV)_{n\in\Z}$ and $(F_2^nV)_{n\in\Z}$ two admissible filtrations.
For all $m\in\Z$, there exists $n\in\Z$ such that $F_1^mV\subset F_2^nV$.
\end{coro}
\

An {\it ind Fr\'echet $E$-algebra} is a topological $E$-algebra $B$ which has a {\it multiplicative admissible filtration}, i.e. an admissible filtration $(F^nB)_{n\in\Z}$ of the underlying topological $E$-vector space such that, if $m,n\in\Z$, and, if $b\in F^mB$, $b'\in F^nB$, then $bb'\in F^{m+n}B$.

A {\it Banach} (resp. {\it Fr\'echet}, resp. {\it ind-Banach}) $E$-algebra is an ind Fr\'echet $E$-algebra $B$ which has a multiplicative admissible filtration  $(F^nB)_{n\in\Z}$ such that $F^0B=B$ and $F^1B=0$ (resp. $F^0B=B$, resp. $F^1B=0$).

\subsection{Ind-Fr\'echet representations}\

From now on $E$ will be $\Q_p$. We'll say {\it Banach}, {\it Fr\'echet}, {\it ind-Banach}, {\it ind-Fr\'echet} instead of $\Q_p$-Banach, $\Q_p$-Fr\'echet, ind-$\Q_p$-Banach, ind-$\Q_p$-Fr\'echet. We'll say {\it Banach algebra}, {\it Fr\'echet algebra},... instead of $\Q_p$-Banach algebra, $\Q_p$-Fr\'echet algebra,...\\

The category $\cI\cF(G_K)$ of {\it ind-Fr\'echet representations (of $G_K$)} is the category whose

--  objects are ind-Fr\'echet equipped with a $\Q_p$-linear and continuous action of $G_K$,

-- morphisms are $G_K$-equivariant continuous $\Q_p$-linear map.

The category $\cI\cF(G_K)$ is an additive $\Q_p$-linear category and any morphism
$$f:V_1\lra V_2$$
 has a kernel and a cokernel: the kernel is the $G_K$-stable closed sub-$\Q_p$-vector space which is the kernel of the underlying $\Q_p$-linear map. The cokernel is the quotient of $V_2$ by the $G_K$-stable closed sub-$\Q_p$-vector space which is the closure of $f(V_1)$.
 
 We say that a morphism $f$ is {\it strict} if the map $${\rm Coim}(f)\to {\rm Im}(f)$$ is an homeomorphism.\\

 Similarly one can define in an obvious way the categories $\cB(G_K)$, $\cI\cB(G_K)$ and $\cF(G_K)$ of {\it Banach, ind-Banach, Fr\'echet representations (of $G_K)$}. This is coherent with the definition of $\cB(G_K)$ already given in the introduction.
   
\subsection{The rings $B^+_{dR}$ and $B_{dR}$ and their topologies}\

We denote $B_{dR}$ the usual field of $p$-adic periods. Recall (cf. eg. \cite{Fo88}, \S 1.5) that it is the fraction field of a discrete valuation ring $B^+_{dR}$, that $G_K$ acts naturally on these two $\Q_p$-algebras and  that $\Z_p(1)$  is naturally a $G_K$-stable sub-$\Z_p$-module of $B^+_{dR}$. We chose a generator $t$  of $\Z_p(1)$. This is also a generator of the maximal ideal of $B^+_{dR}$. Therefore, for all $d\in\Z$, the $d^{\rm th}$-power of this ideal is
$$\fil^d B_{dR}= B^+_{dR}.t^d=B^+_{dR}(d)$$
and is stable under $G_K$. For each $d\geq 0$, we set 
$$B_d=B^+_{dR}/\fil^dB_{dR}\ .$$
Recall  (\cite{Fo88}, \S 1.5.3)  that $B_d$ has a natural structure of a Banach algebra on which the action of $G_K$ is continuous, that, in particular, $B_1=C_p$, and that, for each $d\in\N$, the projection $B_{d+1}\to B_d$ is also continuous. Equipped with the topology of the inverse limit, $B^+_{dR}$ becomes a Fr\'echet algebra on which $G_K$ acts continuously. 

For all $n\in\Z$, multiplication by $t^n$ defines a bijection $B^+_{dR}\to\fil^nB_{dR}$ and we equip $\fil^nB_{dR}$ with the induced topology (for which the action of $G_K$ is continuous (but not that multiplication by $t^n$ does not commute with the action of $G_K$). 

If $n\in\Z$, then $\fil^{n+1}B_{dR}$ is closed in $\fil^nB_{dR}$ and we equip $B_{dR}$ with its natural locally convex topology (a sub-$\Z_p$-module $L$ of $B_{dR}$ is open if and only if, for all $n\in\Z$, the $\Z_p$-module $L\cap \fil^nB_{dR}$ is open in $\fil^nB_{dR}$).

We see that $B_{dR}$ is an ind-Fr\'echet $K$-algebra, with $(\fil^nB_{dR})_{n\in\Z}$ as a $G_K$-equivariant multiplicative admissible filtration. In particular $B_{dR}$ has a natural structure of an ind-Fr\'echet $K$-representation of $G_K$.

\subsection{$B^+_{dR}$ and $B_{dR}$-representations}\label{bbdr}

 Any $B^+_{dR}$-module of finite type has a natural structure of a $K$-Fr\'echet and any finite dimensional $B_{dR}$-vector space has a natural structure of an ind-Fr\'echet $K$-vector space. 
 
 A $B^+_{dR}$-repr\'esentation (resp. a {\it $B_{dR}$-representation) (of $G_K$)} is a $B^+_{dR}$-module of finite type (resp. a finite dimensional $B_{dR}$-vector space) equipped with a semi-linear and continuous action of $G_K$.
With the $G_K$-equivariant $B^+_{dR}$-linear maps as morphisms, these representations form a category that we denote $\rpdrp$ (resp. $\rpdr$).

The category $\rptor=\cC^\infty(G_K)$ of {\it torsion $B^+_{dR}$-representations (of $G_K$)} defined in the introduction (\S \ref{introcgk}) is the full subcategory of $\rpdrp$ whose objects are such that the underlying $B^+_{dR}$-module is torsion ($\iff$ of finite length).\\

Recall (e.g. \cite{stack},  Tag 02MN) that a {\it Serre subcategory} $\cC$ of an abelian category $\cA$ is a strictly full subcategory of $\cA$ containing $0$ which is stable under subobjects, quotients and extensions. In particular, this is an abelian category. Given $\cA$ and $\cC$, one can define the quotient category $\cA/\cC$ which is an abelian category, solution of the obvious universal problem.

\begin{prop} \label{propbdr}The category $\ci$ is a Serre subcategory of $\rpdrp$. The functor
$$\rpdrp\to\rpdr\ \ ,\ \ W\mapsto B_{dR}\otimes_{B^+_{dR}}W$$
is essentially surjective and induces an equivalence 
 $$\rpdrp/\ci\overset{\simeq}\lra \rpdr\ .$$
\end{prop}	

\pr The essential surjectivity comes from the fact that, for any $B_{dR}$-representation $W$, there is a $G_K$-stable lattice $B^+_{dR}$-lattice $W^+$. This result itself comes from the fact that if $W_0^+$ is a $B^+_{dR}$-lattice of $W$, then $W_0$ is an ind-Fr\'echet $K$-vector space with  $(t^nW_0^+)_{n\in\Z}$ forming an admissible filtration. For each $w\in W$, the $g(w)$'s for $g\in G_K$ form a compact subset of $W$, hence it is bounded which implies (\cite{Schneider}, prop.5.6) that it is contained in $t^{-n}W_0^+$ for $n\gg 0$. Hence, if $e_1,e_2,\ldots,e_d$ is a basis of $W$ over $B_{dR}$, there exists $n\in\N$ such that $g(e_i)\in t^{-n}W_0^+$ for $1\leq i\leq d$ and $g\in G_K$. Therefore the sub-$B^+_{dR}$-module $W^+$ of $W$ generated by all these $g(e_i)$'s is also contained in $t^{-n}W_0^+$ and is a $G_K$-stable $B^+_{dR}$-lattice of $W$. The continuity of the action of $G_K$ on $W$ implies the continuity of the action on $W^+$ which is an object of $\rpdrp$. We have an obvious identification of $B_{dR}\otimes_{B^+_{dR}}W^+$ to $W$ and the functor is essentially surjective.

The rest of the proof is straightforward. \qed\\

If $W$ is any object of $\ci$, there is an integer $d$ such that the underlying $B^+_{dR}$-module is a $B_d$-module of finite type. As $B_d$ is a Banach $\Q_p$-algebra,  the underlying topological $\Q_p$-vector space is a Banach and $W$ has a natural structure of a $p$-adic Banach representation.\\

In \cite{Fo03}, we proved (th.3.1):
\begin{prop}\label{propforget} The forget functor
 $$\ci\lra\cB(G_K)$$
 is fully faithful.
 \end{prop}

In other words, given a $p$-adic Banach representation $W$ of $G_K$, there is at most one structure of $B^+_{dR}$-module of finite length on $W$ extending the action of $\Q_p$ such that $W$ becomes a torsion $B^+_{dR}$-representation. 

We use this result to identify $\ci$ to a full subcategory of $\cB(G_K)$. \\

We denote $$\wh\cC^\infty(G_K)$$ the full subcategory of $\cI\cF(G_K)$ whose objects are those $W$'s which admit a $G_K$-equivariant admissible filtration $(F^nW)_{n\in\Z}$ such that $F^mW/F^nW\in\cC^\infty(G_K)$ for all $m\leq n$ in $\Z$. By passing to the limit, the previous proposition implies that, on such a $W$, there is a unique structure of $B^+_{dR}$-module such that the action of $G_K$ is semi-linear and each $F^mW$ is a sub-$B^+_{dR}$-module (and this structure is independent of the choice of $(F^nW)_{n\in\Z}$). We also see that $\wh\cC^\infty(G_K)$ is an abelian category and that any morphism of $\wh\cC^\infty(G_K)$ is $B^+_{dR}$-linear.

Moreover $\rp_{B^+_{dR}}(G_K)$ can be identified to a full subcategory of $\wh\cC^\infty(G_K)$. The proposition \ref{propbdr} implies that this is also true for $\rp_{B_{dR}}(G_K)$.\\

\begin{prop}\label{drdeviss} Let $d\in\N$.

 i) Let $W_1$ be an object of $\ci$ such that ${\rm length}_{B^+_{dR}}W_1\geq d$. There exists a finite extension $K'$ of $K$ contained in $\ov\Q_p$ and a $G_{K'}$-stable sub-$B^+_{dR}$-module $W_1'$ of $W_1$ of length $d$.

ii) Let $W_2$ be an object of $\rp_{B^+_{dR}}(G_K)$ with ${\rm length}_{B^+_{dR}}W_2\geq d$. There exists a finite extension $K'$ of $K$ contained in $\ov\Q_p$ and a $G_{K'}$-stable sub-$B^+_{dR}$-module $W'_2$ of $W_2$ such that ${\rm length}_{B^+_{dR}}W_2/W_2'=d$.
\end{prop}

\pr i) Via an obvious induction, we see that it is enough to check it for $d=1$. Replacing $W_1$ by the kernel of the multiplication by $t$ in $W_1$, we may assume that $W_1$ is a $C_p$-representation.\\

Recall some basic facts of Sen's theory \cite{Sen}:

Let $\chi:G_K\to \Z_p^\ast$ be the cyclotomic character, $H_K$ the kernel of $\chi$ and $L=(C_p)^{H_K}$ which is also the completion of  $K_\infty=\ov\Q_p^{H_K}$. We set $\Gamma_K=G_K/H_K={\rm Gal}(K_\infty/K)$. The character $\chi$ factors through a character $\Gamma_K\to\Z_p^\ast$ that we still denote $\chi$.

For any $C_p$-representation $W$ (of $G_K$), denote $W_K^f$ the union of the finite dimensional sub-$K$-vector spaces of $W^{H_K}$ stable under the action of $G_K$ (acting through $\Gamma_K$). This is a finite dimensional $K_\infty$-vector space equipped with a semi-linear action of $\Gamma_K$. With obvious notations,

-- the functor 
$$\rp_{C_p}(G_K)\lra \rp_{K_\infty}(\Gamma_K)\ \ ,\ \ W\mapsto W_K^f$$
 is exact and fully faithful,
 
-- for any $W\in\rp_{C_p}(G_K)$, the obvious map 
$$C_p\otimes_{K_\infty}W_K^f\lra W$$
is an isomorphism,

-- for all $W\in\rp_{C_p}(G_K)$, there exists a unique endomorphism $\a_{W,K}$ of the $K_\infty$-vector space $W_K^f$ such that:

for all $w\in W_K^f$, there is an open subgroup $\Gamma_w$ of $\Gamma_K$ such that, if $\g\in\Gamma_w$, then
$$ \g(w)=\exp(\log(\chi(\g)).\a_{W,K})(w)$$
(the series $\exp(\lambda\a_{W,K})$ converges to an endomorphism of $W_K^f$ for all small enough $\lambda\in\Z_p$).

It is easy to see that, if $K_1$ is a finite extension of $K$ contained in $\ov\Q_p$, then $W^f_{K_1}$ can be identified to $(K_1)_\infty\otimes_{K_\infty}W^f_K$ and that $\a_{W,K_1}$ is the $(K_1)_\infty$-endomorphism of $W^f_{K_1}$ deduced from $\a_{W,K}$ by scalar extension.\\

Chose such a $K_1$ containing an eigenvalue $\lambda$ of $\a_{W,K}$, hence also of $\a_{W,K_1}$ and chose a non zero eigenvector $\w_0\in\ W_{K_1}^f$ for $\a_{W,K_1}^f$. 
There is a finite extension $K'$ of $K_1$ contained in $\ov\Q_p$ such that, for all $\g\in\Gamma_{K'}$, we have
$$\g(w_0)=\exp(\log(\chi(\g)).\lambda).w$$

We can view $w_0$ as a non zero element of $W_{K'}^f$ and we see that for all $b\in K'$ and all $\g\in\Gamma_{K'}$, we  have
$$\g(bw_0)=\g(b).\exp(\log(\chi(\g)).\lambda).w\ ,$$
hence the $K'$-line of $W_{K'}^f$ generated by $w_0$ is stable under the action of $\Gamma_{K'}$. Therefore the $C_p$-line $W'_1$ of $W_1$ generated by $w_0$ is stable under the action of $G_{K'}$.\\

ii) Replacing $W_2$ by $W_2/t^rW_2$ with $r$ big enough, we may assume that $W_2$ is an object of $\ci$. The result follows by duality from the assertion (i) applied to the Pontryagin dual $W=\cL_{B^+_{dR}}(W_2,B_{dR}/B^+_{dR})$ of $W_2$.\qed
\\

\subsection{Almost $C_p$-representations}\label{secalmostcp}\

If $V_1$ and $V_2$ are two objects of $\cI\cF(G_K)$, an {\it almost isomorphism}
$$f:V_1\rightsquigarrow V_2\ \ ,\ \text{also denoted }\wt f: V_1/U_1\lra V_2/U_2,$$
is a triple $f=(U_1,U_2,\wt f)$ where $U_1$ is a  finite dimensional $G_K$-stable sub $\Q_p$-vector space  of $V_1$, $U_2$ is a  finite dimensional $G_K$-stable sub $\Q_p$-vector space  of $V_2$ and $$\wt f:V_1/U_1\to V_2/U_2$$ an isomorphism of ind-Fr\'echet representations.

We say that two objects $V_1$ and $V_2$  of $\cI\cF(G_K)$ are {\it almost isomorphic} if there exists  an almost isomorphism
$$f:V_1\rightsquigarrow V_2\ .$$

We have (cf. \cite{Fo03}, th.5.3):

\begin{prop} \label{almost} Let $V$ be an object of $\cB(G_K)$. The following are equivalent:

i) $V$ is almost isomorphic to a torsion $B^+_{dR}$-representation,

ii) $V$ is almost isomorphic to a $C_p$-representation,

iii) there is $d\in\N$ such that $V$ is almost isomorphic to $C_p^d$ (equipped with the natural action of $G_K$).
\end{prop}

We denote $\cC(G_K)$ the category of {\it almost $C_p$-representations (of $G_K$)}, that is the full subcategory of $\cB(G_K)$ whose objects satisfy the equivalent conditions of the previous proposition. This is coherent with the definition given in the introduction (\S \ref{introcgk}).

The category $\cC(G_K)$ contains $\ci=\rp_{B^+_{dR}}^{\rm tor}(G_K)$ and $\cC^0(G_K)=\rp_{\Q_p}(G_K)$ as full subcategories.\\

A  {\it weak Serre subcategory} $\cB$ of an abelian category $\cA$ is a strictly full subcategory which is abelian, such that the inclusion functor is exact and which is closed under taking extensions.\\

The following results are essentially contained in \cite{Fo03}:

\begin{thm}\label{thalmost} The category $\cC(G_K)$ is abelian and any morphism of $\cC(G_K)$ is strict as a morphism of $\cB(G_K)$.  A sequence of morphisms of $\cC(G_K)$ is exact if and only if the underlying sequence of $\Q_p$-vector spaces is exact. The category $\co$ is a Serre subcategory of $\cC(G_K)$and $\ci$ is a weak Serre subcategory of $\cC(G_K)$.

Moreover:

i) If $U\in\co$ and $W\in\ci$, then $\Hom_{\cC(G_K)}(W,U)=0$.

ii There exists two additive functions
$$d: {\rm Ob}\ \cC(G_K)\lra \N\ \text{ and }\ h: {\rm Ob}\ \cC(G_K)\lra \Z$$
uniquely determined respectively by $d(U)=0$ if $U\in\cC^0(G_K)$ and $d(C_p)=1$ (resp. $h(U)=\dim_{\Q_p}(U)$ if $U\in\cC^0(G_K)$ and $h(C_p)=0$).

Furthermore, if $W\in\ci$, then $d(W)={\rm length}_{B^+_{dR}}(W)$ and $h(W)=0$.

\end{thm}

\pr This is essentially the theorem 5.1 of \cite{Fo03} except for:

-- The fact that $\cC^0(G_K)$ is a Serre's subcategory of $\cC(G_K)$ which is a triviality.

-- The fact that $\cC^\infty(G_K)$ is a weak Serre's subcategory of $\cC(G_K)$. The only thing which is not obvious is the stability under extensions of $\cC^\infty(G_K)$ inside of $\cC(G_K)$, which is contained in {\it loc.cit.}, prop.6.3.

-- The fact that if $U\in\co$ and $W\in\ci$, then $\Hom_{\cC(G_K)}(W,U)=0$, which is the corollary of the theorem 5.1 of {\it loc.cit.}. \qed\\

For instance, we see that, if $U$ is a $G_K$-stable finite dimensional sub $\Q_p$-vector space of $C_p$, then $d(C_p/U)=1$ and $h(C_p/U)=-\dim_{\Q_p}U$.

If $V\in\cC(G_K)$, $W\in\cC^\infty(G_K)$ and $\wt f:V/U_+\to W/U_-$ is an almost isomorphism, from the diagramme
$$\begin{matrix} 0&\lra& U_+&\lra& V& \lra &V/U_+&\lra& 0\cr
&&&&&&\downarrow\simeq\cr
0&\lra&U_-&\lra& W&\lra &W/U_-&\lra&0\end{matrix}$$
whose lines are exact, we deduce that
$$d(V)=d(W)\ \ ,\ \ h(V)=h(U_+)-h(U_-)=\dim_{\Q_p}(U_+)-\dim_{\Q_p}(U_-)\ .$$

\begin{coro}\label{vw} i) For any $V\in\cC(G_K)$, we have $V\in\cC^0(G_K)\iff d(V)=0$ (in which case $h(V)=\dim_{\Q_p}V\geq 0$).

ii) If $g:V\to W$ is a monomorphism of $\cC(G_K)$ with $W\in\cC^\infty(G_K)$ such that $d(V)=d(W)$, then $g$ is an isomorphism.
\end{coro}

\pr : Looking at an almost isomorphism as above, the first assertion is immediate.

For the second, let $U$ be the cokernel of $g$. We have $d(U)=0$, hence $U\in\cC^0(G_K)$, hence $U=0$, as there is no non trivial morphism from $W$ to an object of $\cC^0(G_K)$. \qed\\

\begin{remk} {\rm As $\co$ is a Serre subcategory of $\cC(G_K)$, we may consider the quotient  $$\wt\cC(G_K)=\cC(G_K)/\co$$
 It is known (\cite{Fo03}, prop.7.1) that this abelian category  is semi-simple with exactly one isomorphism class of simple objects which is the class of $C_p$ viewed as an object of this category. Hence $\wt\cC(G_K)$  is completely determined, up to equivalence, by the somewhat mysterious huge skew field $\cD_K$ of the endomorphisms of $C_p$ in this category ({\it loc.cit.}, prop.7.2).}
\end{remk}\

We denote $$\wh\cC(G_K)$$ the full subcategory of $\cI\cF(G_K)$ whose objects are those $V$'s which admit a $G_K$-equivariant admissible filtration $(F^nV)_{n\in\Z}$ such that $F^mV/F^nV\in\cC(G_K)$ for all $m\leq n$ in $\Z$. 

By passing to the limit, we see that the previous theorem implies:

 \begin{prop} \label{cifab} Any morphism of $\wh\cC(G_K)$ is strict (as a morphism of $\cI\cF(G_K)$) and this category is abelian. A sequence of morphisms of $\wh\cC(G_K)$ is exact if and only if the underlying sequence of $\Q_p$-vector spaces is exact.The category $\cC(G_K)$ is a Serre's subcategory of $\wh\cC(G_K)$ though $\wh\cC^\infty(G_K)$ is a weak Serre's subcategory.
 \end{prop}
 
 \begin{rem} \label{serrebdr} {\rm As $\rp_{B^+_{dR}}(G_K)$ and $\rp_{B_{dR}}(G_K)$ are Serre's subcategories of $\wh\cC^\infty(G_K)$, these two categories are also weak Serre's subcategory of $\wh\cC(G_K)$.}
 \end{rem}

\subsection{Almost split exact sequences}\

We say that a sequence of morphisms of $\cI\cF(G_K)$ is {\it exact} if the underlying sequence of $\Q_p$-vector spaces is exact.

An {\it almost splitting} of a short exact sequence $$0\lra V'\lra V\lra V''\lra 0$$
in $\cI\cF(G_K)$  is a $G_K$-stable closed sub-$\Q_p$-vector space $S$ of $V$ such that \begin{enumerate}
\item the compositum $S\subset V\to V''$ is onto,
\item the $\Q_p$-vector space $S\cap V'$ is finite dimensional.
\end{enumerate}

We say that such an exact sequence {\it almost splits} if there exists such an almost splitting. This is equivalent to saying that there exists a $G_K$-stable finite dimensional sub-$\Q_p$-vector space $U$ of $V'$ such that the sequence
$$0\lra V'/U\lra V/U\lra V''\lra 0$$
splits. 

We observe that any almost splitting $S$ of a short exact sequence
$$0\lra V'\lra V\lra V''\lra 0$$
defines, in an obvious way, almost isomorphisms
$$V\rightsquigarrow V'\oplus V''\rightsquigarrow S\oplus V''\ .$$

Let's recall two results of \cite{Fo03} (th.5.2 and prop.5.2):

 \begin{prop} Let 
$$0\lra W'\lra W\lra W''\lra 0$$
a short exact sequence in $\cI\cF(G_K)$ with $W'$ and $W''$ in $\ci$. Then $W\in \ci$ if and only if the sequence almost splits.
\end{prop}

\begin{prop} \label{almostsplit} Let 
$$0\lra V'\lra V\lra V''\lra 0$$
a short exact sequence in $\cI\cF(G_K)$ with $V'$ and $V''$ in $\cC(G_K)$. Then $V$ is in $\cC(G_K)$ if and only if the sequence almost splits.
\end{prop}

\begin{cor} Among the strictly full subcategories of $\cB(G_K)$ which are abelian, containing $C_p$ and $\cC^0(G_K)$ and stable under almost split extensions, there is a smallest one. This is $\cC(G_K)$.
\end{cor}
 
 \pr Clear ! \qed

\section{$B_e$-representations and coherent $\O_X[G_K]$-modules}

\subsection{The topological $\Q_p$-algebras $B^+_{cris}$ and  $B_e$}\label{bcris}\

Recall (cf. e.g. \cite{Fo88}, \S 2.3 and 4.1) that $B^+_{cris}$ is a Banach algebra equipped with a continuous endomorphism $\ph$ and 
a continuous action of $G_K$ commuting with $\ph$. 
There is a natural $G_K$-equivariant continuous injective homomorphism of topological $\Q_p$-algebras
$$B^+_{cris}\lra B^+_{dR}$$
that we use to identify $B^+_{cris}$ to a subring of $B^+_{dR}$ containing $t$.

For each $d\in\N$, we set 
$$P^d=\{b\in B^+_{cris}\mid \ph(b)=p^db\}\ .$$
This is a $G_K$-stable closed sub-$\Q_p$-vector space of $B^+_{cris}$ as well as of $B^+_{dR}$. Moreover $B^+_{cris}$ and $B^+_{dR}$ induce the same topology on $P^d$ which can be viewed as a Banach representation of $G_K$.  We have a a canonical short exact sequence
$$0\lra \Q_p(d)\lra P^d\lra B_d\lra 0$$
where $Q_p(d)=\Q_pt^d$  and $P^d\to B_d$ is the compositum $P^d\subset B^+_{cris}\subset B^+_{dR}\overset{{\rm proj}}\lra B_d$. In particular we see that $P^d$ is an almost $C_p$-representation with $d(P^d)=d$ and $h(P^d)=1$.\\

As usual, we set $B_{cris}=B^+_{cris}[1/t]$ that we can view as a $G_K$-stable subring of $B_{dR}$.

We have $\ph(t)=pt$ and $\ph$ extends uniquely to $B_{cris}$. Moreover the natural map $B_{cris}\to B_{dR}=B^+_{dR}[1/t]$ is still injective and we use it to identify $B_{cris}$ to a $G_K$-stable sub-$\Q_p$-algebra of $B_{dR}$.

Recall that 
$$B_e=\{b\in B_{cris}\mid\ph(b)=b\}$$
is also a $G_K$-stable sub-$\Q_p$-algebra of $B_{dR}$. We endow it with the topology induced by  the (locally convex) topology of $B_{dR}$.

Then, we have
$$B_e=\limi_{d\in\N}\fil^{-d}B_e=\cup_{d\in\N}\fil^{-d}B_e$$
where, for all $d\in\N$,
$$ \fil^{-d}B_e=B_e\cap B^+_{dR}t^{-d}=P^d.t^{-d}=P^d(-d)$$ is an  almost $C_p$-representation (with $d(P^d(-d))=d$ and $h(P^d(-d))=1$) homeomorphic to $P^d$ as a Banach. Setting $P^d=P^d(-d)=0$ for $d>0$, we see that $B_e$ is an ind-Banach algebra with $(P^{-n}(n))_{n\in\Z}$ a $G_K$-stable multiplicative admissible  filtration.

\subsection{$B_e$-representations}\

The topology of $B_e$ induces on each $B_e$-module of finite type a natural topology for which it is an ind-Fr\'echet (actually an ind-Banach). A {\it $B_e$-representation (of $G_K$)} is a $B_e$-module of finite type equipped with a semi-linear and continuous action of $G_K$. With the $G_K$-equivariant $B_e$-linear maps as morphisms, $B_e$-representations form a category that we denote $\rp_{B_e}(G_K)$.

\begin{prop} The $B_e$-module underlying any $B_e$-representation is free of finite rank. The category $\rp_{B_e}(G_K)$ is a  $Q_p$-linear abelian category.
\end{prop}

\pr 

Recall (\cite{FF}, th.6.5.2) that $B_e$ is a principal ideal domain. In particular it is a noetherian ring and the fact that $\rp_{B_e}(G_K)$ is a $\Q_p$-linear abelian category is obvious.

 Moreover ({\it loc.cit.}, prop.10.1.1), for any maximal ideal $\fp$ of $B_e$, the orbit of $\fp$ under the action of $G_K$ is infinite. This implies that there is no non trivial $G_K$-equivariant ideal of $B_e$.
If $\La$ is any non zero $B_e$-representation of $G_K$, the annihilator of its torsion sub-module is a proper $G_K$-equivariant ideal and must be $0$. Therefore the $B_e$-module underlying $\La$ is torsion free, hence free of finite rank.
 \qed\\

\begin{remk}
{\rm Let $C_e$ the fraction field of $B_e$. This is the union of the fractional ideals of $B_e$. For each such ideal $\fa$,  the choice of a generator $a$ defines a bijection
$$B_e\lra \fa\ \ ,\ \ b\mapsto ba$$
and we put on $\fa$ the topology defined by transport de structure, which is independent of the choice of the generator. Hence each $\fa$ is naturally an ind-Banach ($\Q_p$-vector space). If $\fa\subset\fb$ are two fractional ideals, this inclusion is continuous and $\fa$ is a closed sub-$\Q_p$-vector space of $\fb$. Hence we may endow $C_e$ with the coarsest locally convex topology such that, for all fractional ideal $\fa$, the map $\fa\to C_e$ is continuous (a lattice $\cL$ in $C_e$ is open if and only if $\cL\cap\fa$ is open in $\fa$ for all $\fa$).

The action of $G_K$ on $C_e$ is continuous for this topology (but $C_e$ doesn't seem to be an object of $\cI\cF(G_K)$) and we may consider the category $\rp_{C_e}(G_K)$ of $C_e$-representations (of $G_K)$, that is of finite dimensional $C_e$-vector spaces equipped with  a semi-linear an continuous action of $G_K$. This is obviously a $\Q_p$-linear abelian category.

We have an obvious exact $\Q_p$-linear functor
$$\rp_{B_e}(G_K)\lra \rp_{C_e}(G_K)\ \ ,\ \ \La\mapsto C_e\otimes_{B_e}\La\ .$$
This functor is fully faithful: if $M\in \rp_{C_e}(G_K)$ is a $C_e$-representation of dimension $d$, there is at most one $G_K$-equivariant sub-$B_e$-module of rank $d$ because if $\La_1$ and $\La_2$ are two of them, so are $\La_1+\La_2$ and $(\La_1+\La_2)/\La_1$ is torsion, 
hence $0$.}
\end{remk}

\begin{remk}
{\rm If $\La$ is any $B_e$-representation of $G_K$, the underlying $\Q_p$-vector space is locally convex and $\La$ inherits a natural structure of an object of $\cI\cF(G_K)$. We'll see later  that the forget functor $$\rp_{B_e}(G_K)\lra \cI\cF(G_K)$$ is fully faithful (prop.\ref{forget}) and that its essential image is contained in $\wh\cC(G_K)$ (prop.\ref{forgete}).}
\end{remk}

\begin{prop}\label{wlambda} Let $W\in\cC^\infty(G_K)$ and $\Lambda\in\rp_{B_e}(G_K)$. Then
$$\Hom_{\cI\cF(G_K)}(W,\Lambda)=0\ .$$
\end{prop}

\pr Let $f :W\to\Lambda$ such a morphism. We see that $B_{dR}\otimes_{B_e}\Lambda$ is a $B_{dR}$-representation of $G_K$ and that
$$g:\Lambda\lra B_{dR}\otimes_{B_e}\Lambda\ \ ,\ \ \lambda\mapsto 1\otimes\lambda$$
is a morphism of $\cI\cF(G_K)$. But $gf:W\to B_{dR}\otimes_{B_e}\Lambda$ must be $B^+_{dR}$-linear (\S \ref{bbdr}). As the $B^+_{dR}$-module $W$ is torsion though $B_{dR}\otimes\Lambda$ is torsion free, we have $gf=0$, hence also $f=0$ as $g$ is injective. \qed

\subsection{Coherent $\O_X$-modules}\label{coh}\

We know that $B_e$ is a PID and we may consider the "open curve"
$$X_e=\spec\ B_e,$$
a noetherian regular affine scheme of dimension $1$ whose function field is the fraction field $C_e$ of $B_e$ that we can see as a subfield of $B_{dR}$. For each closed point $x$ of $X$, the local ring $\O_{X,x}$ is a DVR and we denote $v_x$ the corresponding valuation on $C_e$ normalised by $v_x(C_e^\ast)=\Z$.

Recall (cf. \cite{FF}, \S 6.5.1) that the  curve $X=X_{\Q_p,C_p^\flat}$  can be defined as the compactification at infinity of $X_e$. More precisely, as $B_{dR}$ is the fraction field of the discrete valuation ring $B^+_{dR}$, it is naturally equipped with a valuation $v_{dR}$: if $b\in B_{dR}$ is $\not=0$, then $v_{dR}(b)$ is the biggest $n\in\Z$ such that $b\in\fil^nB_{dR}$. We denote $v_{\infty}$ the restriction of $v_{dR}$ to $C_e$. The topological space underlying $X$ is obtained from the topological space underlying $X_e$ by adding the closed point $\infty$ defined by $v_\infty$. Hence, the function field of $X$ is $C_e$ and, if $U$ is any non empty open subspace of $X$, we have
$$\O_X(U)= \{b\in C_e\mid v_x(b)\geq 0,\ \forall x\in U\}\ .$$

 We have $X\setminus\{\infty\}=X_e$, the ring $B^+_{dR}$ is the completion of $\O_{X,\infty}$ and $B_{dR}$ is the completion of $C_e$ for the topology defined by $v_\infty$.\\

Consider the following category $\oxm$:

\begin{itemize}
\item an object of $\oxm$ is a triple $(\cF_e,\cF^+_{dR},\iota_\F)$ with $\cF_e$ a $B_e$-module of finite type, $\cF^+_{dR}$ a $B^+_{dR}$-module of finite type and
$$\iota_\cF :\cF^+_{dR}\lra B_{dR}\otimes_{B_e}\cF_e$$
a $B^+_{dR}$-linear map inducing an isomorphism of $B_{dR}$-vector spaces
$$ B_{dR}\otimes_{B^+_{dR}}\cF^+_{dR}\lra B_{dR}\otimes_{B_e}\cF_e\ .$$
\item a morphism  $(\cF_e,\cF^+_{dR},\iota_\F)\lra  (\cG_e,\cG^+_{dR},\iota_\cG)$ is a pair $(f_e,f^+_{dR})$ with $f_e:\cF_e\to\cG_e$ a $B_e$-linear map and $f^+_{dR}:\cF^+_{dR}\to\cG^+_{dR}$ a $B^+_{dR}$-linear map such that the obvious diagram commutes.
\end{itemize}\

To any coherent $\O_X$-module $\cF$, we can associate an object $(\cF_e,\cF^+_{dR},\iota_\F)$ of this category:
\begin{itemize}
\item $\cF_e=\cF(X_e)$,
\item $\cF^+_{dR}=B_{ dR}\otimes_{\O_{X,\infty}}\cF_\infty$, the completion of the fiber of $\cF$ at $\infty$,
\item the completion at $\infty$ of the general fiber is $B_{dR}\otimes_{B^+_{dR}}\cF^+_{dR}$ as well as $B_{dR}\otimes_{B_e}\cF_e$ and $\iota_\cF: \cF^+_{dR}\lra B_{dR}\otimes_{B_e}\cF_e$ is the natural map.
\end{itemize}\

This correspondence is obviously functorial and it's immediate to see that it gives an equivalence of categories. We use it to identify the category of coherent $\O_X$-modules to $\oxm$. In this equivalence we see that the category $Bund(X)$ of vector bundles over $X$, i.e. of torsion free coherent $\O_X$-modules, can be identified to the full subcategory of $\oxm$ whose objects are triples $(\cF_e,\cF^+_{dR},\iota_\F)$ such that the $B_e$-module $\cF_e$ and the $B^+_{dR}$-module $\cF^+_{dR}$ are torsion free ($\iff$ free).\\

\subsection{The topology on $\O_X$}\label{top}

The curve $X$ can be also described ([FF], \S 6.5.1) as
$$X={\rm Proj}\bigoplus_{d\in\N}P^d$$
and there is ({\it loc.cit.}, th.6.5.2) a one to one correspondence between the closed points of $X$ and the $\Q_p$-lines in $P^1$ (the map associating to such a line the prime ideal of $P=\oplus_{d\in\N}P^d$ that it generates is a bijection between the set of these lines and the set of non zero homogeneous prime ideals of $P$ different from $\oplus_{d>0}P^d$). 
In this correspondence $\infty$ corresponds to the line generated by $t$. 

Moreover, if $x_1,x_2,\ldots,x_r$ are closed points of $X$ and if, for $1\leq i\leq r$, we chose a generator $t_i$ of the $\Q_p$-line associated to $x_i$, we see that the $\Q_p$-algebra $\O_X(X\setminus\{x_1,x_2,\ldots,x_r\})$ has a natural topology: If we set $u=t_1t_2\ldots t_r$, we have$$O_X(X\setminus\{x_1,x_2,\ldots,x_r\})=\cup_{n\in\N}P^{nr}u^{-n}$$
and we see that it is an ind-Banach algebra with $(P^{nr}u^{-n)})_{n\in\N}$ a multiplicative admissible Banach filtration.
Thus we may consider $\O_X$ as a sheaf of ind-Banach algebras (the restriction maps are obviously continuous).

\subsection{The category $\cM(G_K)$}\label{oxgk}\

The group $G_K$ acts continuously on $X$ and it makes sense to speak of the category $\cM(G_K)$ of coherent $\O_X[G_K]$-modules, that is of coherent $\O_X$-modules equipped with a semi-linear and continuous action of $G_K$. 

 We see that:

-- the open subset $X_e=\spec\ B_e$ is stable under $G_K$ and 
$G_K$ acts continuously on the ind-Banach algebra $B_e$, 

-- the point $\infty$ is fixed by $G_K$ and the action of $G_K$ on  the Fr\'echet algebra $B^+_{dR}$ (resp. on the ind-Fr\'echet algebra $B_{dR}$), completion at $\infty$ of $\O_{X,\infty}$ (resp. of the function field $C_e$ of $X$) is continuous.

From the description of coherent $\O_X$-modules of the previous paragraph, we see that we can identify $\cM(G_K)$ to the category whose
\begin{itemize} 
\item objects are triple
$\cF=(\cF_e,\cF^+_{dR},\iota_\cF)$,
where $\cF_e$ is a $B_e$-representation, $\cF^+_{dR}$ is a $B^+_{dR}$-representation and  $$\iota_\cF : \cF^+_{dR}\lra B_{dR}\otimes_{B_e}\cF_e$$
is a $G_K$-equivariant homorphism of $B^+_{dR}$-modules such that the induced $B_{dR}$-linear map
$$B_{dR}\otimes_{B^+_{dR}}\cF^+_{dR}\lra B_{dR}\otimes_{B_e}\cF_e$$
is bijective,

\item a morphism $$f:(\cF_e,\cF^+_{dR},\iota_\F)\lra (\cG_e,\cG^+_{dR},\iota_\G)$$
is a pair $(f_e,f^+_{dR})$ with $f_e:\cF_e\to\cG_e$ (resp. $f^+_{dR}:\cF^+_{dR}\to\cG^+_{dR}$) a morphism of $B_e$-representations (resp. $B^+_{dR}$-representations) such that the obvious diagramme commutes.
\end{itemize}
\

When there is no ambiguity about the map $\iota_\cF$, we write abusively
$$\cF=(\cF_e,\cF^+_{dR})$$

We also denote 
$$\cF_{dR}=B_{dR}\otimes_{B_e}\cF_e=B_{dR}\otimes_{C_e}(C_e\otimes_{B_e}\cF_e)$$
the completion at $\infty$ of the generic fiber $\cF_\eta=C_e\otimes_{B_e}\cF_e$ of $\cF$.\\

The category ${\rm Bund}_X(G_K)$ of {\it $G_K$-equivariant vector bundles over $X$} is the full subcategory of $\cM(G_K)$ whose objects are those for which the underlying $\O_X$-module is torsion free. From the fact that any $B_e$-representation is torsion free, we see that, if $\cF$ is any coherent $\O_X[G_K]$-module,  there is no torsion away from $\infty$. Therefore ${\rm Bund}_X(G_K)$ is the full subcategory of $\cM(G_K)$ whose objects are those $\cF$ such that the $B^+_{dR}$-module $\cF^+_{dR}$ is free ($\iff$ torsion free).

\subsection{The Harder-Narasimhan filtration}\label{hnf}\

The abelian category $\oxm$ is equipped with two additive functions, the {\it rank} and the {\it degree} (cf \cite{FF}, chap.5)
$${\rm rk} :  \oxm\to \N\ \ ,\ \ \deg:\oxm\to \Z$$
The rank of $\cF=(\cF_e,\cF^+_{dR},\iota_\cF)$ is the rank of the $B_e$-module $\cF_e$. It is $0$ if and only if $\cF$ is torsion. It's more difficult to compute the degree. But this additive function is characterised by the following facts:
\begin{itemize}
\item  if $D$ is a divisor, then $$\deg(\cL(D))=\deg(D)=\sum_{{\rm closed\ points\ of }X}n_x\ \ ,\ \text{ if } D=\sum n_x[x]\ ,$$
\item if $\cF$ is a vector bundle of rank $r$, then
$$\deg(\cF)=\deg(\wedge^r\cF)\ ,$$
\item if $\cF$ is a torsion $\O_X$-module, then
$$\deg(\cF)=\sum_{{\rm closed\ points\ of }X} {\rm length}_{\O_{X,x}}\cF_x\ .$$
\end{itemize}

The slope of a non zero coherent $O_X$-module $\cF$ is
$${\rm slope}(\cF)=\deg(\cF)/{\rm rank}(\cF)\in\Q\cup\{+\infty\}$$
(with the convention that the slope of a non zero torsion coherent $\O_X$-module is $+\infty$). 

As in the classical case, \begin{itemize}
\item a coherent $\O_X$-module $\cF$ is {\it semi-stable} if it is non zero and if ${\rm slope}(\cF')\leq {\rm slope}(\cF)$ for any non zero coherent sub-$\O_X$-module of $\cF$,
\item the Harder-Narasimhan filtration of a coherent $\O_X$-module $\cF$ is the unique increasing filtration
$$0=\cF_0\subset\cF_1\subset\ldots\subset\cF_{m-1}\subset\cF_m=\cF$$
by coherent sub-$\O_X$-modules such that each $\cF_i/\cF_{i-1}$ is semi-stable with
$${\rm slope}(\cF_1/\cF_0)>{\rm slope}(\cF_2/\cF_1)>\ldots>{\rm slope}(\cF_{m-1}/\cF_{m-2})>{\rm slope}(\cF_m/\cF_{m-1})$$
\end{itemize}

The slopes of the $\cF_i/\cF_{i-1}$ for $1\leq i\leq m$ are called {\it the HN-slopes of $\cF$}.\\

The Harder-Narasimhan filtration splits continuously but not canonically.\\

If $\cF$ is an object of $\cM(G_K)$, the unicity of the Harder-Narasimhan filtration implies that this filtration is by subobjects in $\cM(G_K)$. In general, there is no $G_K$-equivariant splitting of this filtration.

\subsection{The equivalences $\cM^0(G_K)\to \cC^0(G_K)$ and $\cM^\infty(G_K)\to \ci$}\

For all $s\in\Q\cup\{+\infty\}$, we denote $\cM^s(G_K)$ the full subcategory of $\cM(G_K)$
 whose objects are semi-stable of slope $s$. We also write $\cM^\infty(G_K)=\cM^{+\infty}(G_K)$.\\
 
We have $H^0(X,\O_X)=\Q_p$. A central result of \cite{FF} (th.8.2.10)  is that a coherent $\O_X$-module $\cF$ is semi-stable of slope $0$ if and only if it is isomorphic to $\O_X^r$ for some positive integer $r$. From that we deduce:

\begin{prop}\label{moco} If $\cF\in\cM^0(G_K)$, then $\cF(X)\in\cC^0(G_K)$ and $${\rm rank}(\cF)=\dim_{\Q_p}\cF(X)\ .$$
The functor
$$\cM^0(G_K)\lra \cC^0(G_K)\ \ ,\ \ \cF\mapsto \cF(X)$$
is an equivalence of categories. The functor
$$\cC^0(G_K)\lra\cM^0(G_K)\ \ ,\ \ U\mapsto \O_X\otimes U= (B_e\otimes_{\Q_p}U,B^+_{dR}\otimes_{\Q_p}U)$$
is a quasi-inverse.
\end{prop}
\

 If $\cF\in\cM(G_K)$, as there is no torsion away from infinity, $\cF\in\cM^\infty(G_K)$ if and only if $\cF_e=0$. From that, we deduce:
 
 \begin{prop}\label{mici} If $\cF\in\cM^\infty(G_K)$, then $\cF(X)=\cF^+_{dR}$ and belongs to $\ci$. Moreover $$\deg(\cF)={\rm length}_{B^+_{dR}}\cF(X)\ .$$
 The functor
$$\cM^\infty(G_K)\lra \ci\ \ ,\ \ \cF\mapsto \cF(X)$$
is an equivalence of categories. The functor
$$\ci\lra\cM^\infty(G_K)\ \ ,\ \ W\mapsto \u W= (0,W)$$
is a quasi-inverse.
\end{prop}\

For any $s\in\Q$, we denote $\cM^{\geq s}(G_K)$ (resp. $\cM^{<s}(G_K)$) the full subcategory of $\cM(G_K)$ whose objects are those which have all their HN-slopes $\geq s$ (resp. $<s$).

For any $\cF\in\cM(G_K)$, we denote $\cF^{\geq 0}$ the biggest term of the Harder-Narasimhan filtration which belongs to $\cM^{\geq 0}(G_K)$ and $\cF^{<0}=\cF/\cF^{\geq 0}$.
We have a short exact sequence 
$$0\lra \cF^{\geq 0}\lra \cF\lra \cF^{<0}\lra 0$$
with $\cF^{\geq 0}\in\cM^{\geq 0}(G_K)$ and $\cF^{<0}\in\cM^{<0}(G_K)$.\\

The category $\cM(G_K)$ is equipped with a tensor product. From the classification of vector bundles over $X$ (cf. \cite{FF}, th.8.2.10), we get the fact that if $s,t\in\Q\cup\{+\infty\}$, if $\cF\in\cM^s(G_K)$ and if $\cG\in\cM^t(G_K)$, then $\cF\otimes\cG\in\cM^{s+t}(G_K)$ (with the convention that $s+t=+\infty$ if $s$ or $t$ is $+\infty$).

The additive category ${\rm Bund}_X(G_K)$ has an internal hom
$$(\cF,\cG)\mapsto Hom_{\O_X}(\cF,\cG)$$
We see that $(Hom_{\O_X}(\cF,\cG))_e=\cL_{B_e}(\cF_e,\cG_e)$ is the $B_e$-module of the $B_e$-linear maps  $\cF_e\to\cG_e$ though  $(Hom_{\O_X}(\cF,\cG))^+_{dR}=\cL_{B^+_{dR}}(\cF^+_{dR},\cG^+_{dR})$ is the $B^+_{dR}$-module of the $B^+_{dR}$-linear maps  $\cF^+_{dR}\to\cG^+_{dR}$.

In ${\rm Bund}_X(G_K)$, there is also a duality: The dual of $\cF$ is $\cF^\vee=Hom_{\O_X}(\cF,\O_X)$. If $\cF,\cG\in{\rm Bund}_X(G_K)$, then
$Hom_{\O_X}(\cF,\cG)=\cF^\vee\otimes\cG$. If $\cF$ is semi-stable of slope $s$, then $\cF^\vee$ is semi-stable of slope $-s$.

\subsection{Tate and Harder-Narasimhan twists}\label{twists}\

Recall that, for any $p$-adic vector space $V$ equipped with a linear action of $G_K$ and $n\in\Z$, the $n^{\rm th}$ Tate's twist of $V$ is $$V(n)=V\otimes_{\Q_p}\Q_p(n)$$
where $\Q_p(n)=\Q_pt^n\subset B_{dR}$. This construction is functorial.\\

For any $n\in\Z$, we denote $$\O_X(n)_{T}=\O_X\otimes \Q_p(n)=(B_e(n),B^+_{dR}(n))=(B_e.t^n,B^+_{dR}.t^n)$$ 
(where $B_e.t^n$ (resp. $B^+_{dR}.t^n$) is the sub $B_e$-module (resp. $B^+_{dR}$-module) of $B_{dR}$ generated by $t^n$)
the $G_K$-equivariant line bundle of slope $0$ associated to $\Q_p(n)$.\\

For $\cF\in\cM(G_K)$ and $n\in\Z$, the {\it $n^{\rm th}$ Tate twist} of $\cF$ is
$$\cF(n)_{T}=\cF\otimes\O_X(n)_{T}=(\cF_e(n),\cF^+_{dR}(n),\iota_\cF(n))\ .$$
It has the same degree, the same rank and the same slope as $\cF$. \\

For any $n\in\Z$, we consider the $G_K$-equivariant line bundle
$$\O_X(n)_{HN}= (B_e,B^+_{dR}(-n))=(B_e,B^+_{dR}.t^{-n})\ $$
There is an obvious short exact sequence
$$\begin{matrix} 0\lra \O_X\lra \O_X(n)_{HN}\lra (0,B_n(-n))\lra 0&\text{ if } n\geq 0,\cr
 0\lra \O_X(n)_{HN}\lra \O_X\lra (0,B_{-n})\lra 0&\text{ if } n< 0,\end{matrix}$$
In particular, $\O_X(n)_{HN}$ is a modification of $\O_X$ and is of degree $n$. It is semi-stable of slope $n$.

For $\cF\in\cM(G_K)$ and $n\in\Z$, we define the {\it $n^{\rm th}$ Harder-Narasimhan twist} of $\cF$ as
$$\cF(n)_{HN}=\cF\otimes\O_X(n)_{HN}=(\cF_e,\cF^+_{dR}(-n),\iota_\cF(-n))=(\cF_e,t^{-n}.\cF^+_{dR},\iota_\cF(-n))\ .$$

It has the same  rank as $\cF$. If $\cF$ is semi-stable of slope $s$, then $\cF(n)_{HN}$ is semi-stable of slope $s+n$.\\

These two construction are obviously functorial and commute with Harder-Narasimhan filtration. In particular,

-- if $\cF$ is semis-table of slope $s$, then $\cF(n)_T$ is semi-stable of slope $s$, though $\cF(n)_{HN}$ is semi-stable of slope $s+n$,

--  the HN-slopes of $\cF(n)_T$ are the same as the HN-slopes of $\cF$, though the HN-slopes of $\cF(n)_{HN}$ are the $s+n$ for $s$ describing the HN-slopes of $\cF$.\\

  These constructions commute: for $m,n\in\Z$, we have
$$\cF(m)_T(n)_{HN}=\cF(n)_{HN}(m)_T\ .$$\

\begin{remk} {\rm In \cite{FF} (def.8.2.1) the $G_K$-equivariant line bundle $\O_X(n)_{HN}(n)_T$ is denoted $\O_X(n)$. We have to avoid confusion between the three $G_K$-equivariant line bundles $\O_X(n)_T$, $\O_X(n)_{HN}$ and 
$$\O_X(n)=(B_e(n),B^+_{dR})=(B_e.t^n,B^+_{dR})\ .$$}
\end{remk}\

\subsection{Potentially trivialisable $B_e$-representations}\

Let $\La$ be a $B_e$-representation of $G_K$ and $K'$ a finite extension of $K$ contained in $\ov\Q_p$. We say that $\La$ is {\it $G_{K'}$-trivialisable} if there is $U\in\cC^0(G_{K'})$ and a $G_{K'}$-equivariant isomorphism of $B_e$-modules
$$B_e\otimes_{\Q_p}U\simeq \La\ .$$

We say that $\La$ is {\it trivialisable} if it is $G_K$-trivialisable and {\it potentially trivialisable} if there is a finite extension $K'$ of $K$ contained in $\ov\Q_p$ such that $\La$ is $G_{K'}$-trivialisable.

\begin{prop} Any  absolutely irreducible $B_e$-representation of $G_K$ is potentially trivialisable.\end{prop}

\pr Let $\La$ such a $B_e$-representation. Then $\La_{dR}=B_{dR}\otimes_{B_e}\La$ is a $B_{dR}$-representation. Let $\cL$ be the set of $G_K$-stable $B^+_{dR}$-lattices of $\La_{dR}$. We know (prop.\ref{propbdr}) that $\cL$ is not empty. For each $L\in\cL$, we may consider the $G_K$-equivariant vector bundle over $X$
$$\cF_L=(\La,L)\ .$$
Such an $\cF_L$ is semi-stable (otherwise the Harder-Narasimhan filtration would be non trivial and would induce a non trivial filtration of the $B_e$-representation $(\cF_L)_e=\La$ which is not possible as $\La$ is irreducible). 

Chose such an $\cF_L$. Replacing $\cF_L$ with $\cF_L(n)_{HN}$ with $n\in\N$ big enough, we may assume that the degree $d$ of $\cF_L$ is $\geq 0$. By proposition \ref{drdeviss}, we can find a finite extension $K'$ of $K$ contained in $\ov\Q_p$ and a $G_K'$-stable sub-$B^+_{dR}$-lattice $L_0\subset L$ such that ${\rm length}_{B^+_{dR}}(L/L_0)=d$. Then $\cF_{L_0}=(\La,L_0)$ is a $G_{K'}$-equivariant vector bundle over $X$ of degree $d-d=0$. As the $B_e$-representation $\La$ is absolutely irreducible, it is irreducible as a $B_e$-representation of $G_{K'}$. Hence, $\cF_{L_0}$ is semi-stable of slope $0$. By proposition \ref{moco}, there is a $\Q_p$-representation $U$ of $G_{K'}$ such that 
$$\cF_{L_0}\simeq \O_X\otimes U\ .$$
Therefore $\La$, as a $B_e$-representation of $G_{K'}$, is isomorphic to $B_e\otimes_{\Q_p}U$.
 \qed\\

\begin{coro}\label{dirtriv} The category $\rp_{B_e}(G_K)$ is the smallest full subcategory of itself containing potentially trivialisable $B_e$-representations and stable under taking extensions. This is also the smallest full subcategory of itself containing trivialisable $B_e$-representations and stable under taking extensions and direct summands.
\end{coro}

\pr For any $B_e$-representation $\La$ of $G_K$, one can find a finite extension $K_1$ of $K$ contained in $\ov\Q_p$ such that $\La$, viewed as a $B_e$-representation of $G_{K_1}$, can be viewed as a successive extension of absolutely irreducible $B_e$-representations of $G_{K_1}$ and the first assumption results from the previous proposition.
Hence we may find a finite extension $K'$ of $K$ contained in $\ov\Q_p$ such that $\La$, as a $B_e$-representation of $G_{K'}$, is a successive extension of $G_{K'}$-trivialisable $B_e$-representations. Therefore the induced $B_e$-representation of $G_K$
$$\La' = B_e[G_K]\otimes_{B_e[G_{K'}]}\La=\Q[G_K]\otimes_{\Q[G_{K'}]}\La$$
is a successive extension of trivialisable $B_e$-representations of $G_K$. But the obvious $G_K$-equivariant projection $\La'\to \La$ splits (as, if $\La^\vee$ denotes the $B_e$-dual of $\La$ and if $H={\rm Gal}(K'/K)$, we have a short exact sequence
$$0\lra \Hom_{\rp_{B_e}(G_K)}(\La,\La')\lra \Hom_{\rp_{B_e}(G_{K'})}(\La,\La)\lra H^1(H,\La^\vee\otimes_{B_e}\La')$$
and, as $B_e$ is of characteristic $0$, we have $H^1(H,\La^\vee\otimes_{B_e}\La')=0$). Therefore, $\La$ is a direct summand of a successive extension of trivialisable $B_e$-representations.
\qed\\

\begin{remk} {\rm The results of this paragraph can also be deduced from the work of Berger (\cite{Be08} and \cite{Be09}) relying $(\ph,\Gamma)$-modules on the Robba ring and $B_e$-pairs.}
\end{remk}

\subsection{The forget functor $\rp_{B_e}(G_K)\to \wh\cC(G_K)$}\label{forgetbe}\

\begin{prop}\label{forget} The forget functor
$$\rp_{B_e}(G_K)\lra \cI\cF(G_K)$$
is fully faithful.
\end{prop} 

\pr Let $\La$ and $\La'$ two $B_e$-representations. We want to prove that any $G_K$-equivariant continuous map 
$$\La\overset{\a}\lra\La'$$ is $B_e$-linear.

Let $K'$ be a finite Galois extension of $K$ contained in $\ov\Q_p$ such that $\La$ and $\La'$ are successive extensions of trivialisable $B_e$-representations of $G_{K'}$. If $H={\rm Gal}(K'/K)$,  we have
$$\Hom_{\rp_{B_e}(G_K)}(\La,\La')=(\Hom_{\rp_{B_e}(G_{K'})}(\La,\La'))^H\ ,\ \Hom_{\cI\cF(G_K)}(\La,\La')=(\Hom_{\cI\cF(G_{K'})}(\La,\La'))^H\ .$$
Therefore, replacing $K$ by $K'$ we may assume again that there is $r\in\N$ and a filtration of $\La$ by sub-$B_e$-representations
$$0=\La_0\subset\La_1\subset\ldots\subset\La_{r-1}\subset\La_r=\La$$
such that each $\La_i/\La_{i-1}$ is trivialisable.

We proceed by induction on $r$, the case $r=0$ being trivial. Assume $r\geq 1$ and that $\La_r/\La_{r-1}=B_e\otimes_{\Q_p}U$ for some $U\in\cC^0(G_K)$. Chose a $B_e$-linear section $s:B_e\otimes U\to \La$ of the projection $\La\to B_e\otimes U$.
We have  a decomposition of $\La$ as a $B_e$-module into a direct sum
$$\La=\La_{r-1}\oplus s(B_e\otimes U)=\La_{r-1}\oplus(B_e\otimes s(U))\ .$$

By induction, the restriction of $\a$ to $\La_{r-1}$ is $B_e$-linear. Hence there is a unique $B_e$-linear map
$$\a_0:\La\lra\La'$$
such that $\a_0(\lambda)=\a(\lambda)$ if $\lambda\in\La_{r-1}$ and $\a_0(s(u))=\a s(u))$ for all $u\in U$. It is easy to check that $\a_0$ is continuous and $G_K$-equivariant. The maps
$$\a,\a_0: \La\lra\La'$$
coincide on $\La_{r-1}\oplus s(U)$ and the map $\a-\a_0$ induces, by going to the quotient, a morphism in $\cI\cF(G_K)$
$$\b:\La/(\La_{r-1}\oplus s(U))\lra \La'\ .$$
Recall (cf. eg \cite{CF}, prop.1.3) that $B_{dR}=B_e+B^+_{dR}$ though $B_e\cap B^+_{dR}=\Q_p$.
Hence, if we set $\wt B_{dR}=B_{dR}/B^+_{dR}$, we can identify $B_e/\Q_p$ to $\wt B_{dR}$. 

Therefore we have
$$\La/(\La_{r-1}\oplus s(U))=(\La_r/\La_{r-1})/U=B_e\otimes U/U=\wt B_{dR}\otimes_{\Q_p}U\ .$$
and $\b\in\Hom_{\cI\cF(G_K)}(\wt B_{dR}\otimes U,\La')$.

We see that $\wt B_{dR}$ is the direct limit of the $B_d(-d)$, for $d\in\N$, hence
$$\wt B_{dR}\otimes U=\limi_{d\in\N}B_d(-d)\otimes_{\Q_p}U\ .$$
Each $B_d(-d)\otimes U$ is an object of $\cC^\infty(G_K)$. Hence, proposition \ref{wlambda}, implies that $$\Hom_{\cI\cF(G_K)}(B_d(-d)\otimes U,\La')=0\ .$$
Therefore $\b=0$ and $\a=\a_0$ is $B_e$-linear. \qed\\

We use this result to identify $\rp_{B_e}(G_K)$ to a full subcategory of $\cI\cF(G_K)$.

\begin{prop} \label{forgete} We have $$\rp_{B_e}(G_K)\subset\wh\cC(G_K)\ .$$
More precisely, for any $B_e$-representation $\La$ of $G_K$, there is a $G_K$-equivariant admissible filtration $(F^n\La)_{n\in\Z}$ with $F^1\La=0$ and $F^n\La\in\cC(G_K)$ for all $n$.
Moreover, we may chose this filtration so that,  if $b\in \fil^{-d}B_e$ and $\lambda\in F^n\La$ (with $d\in\N$, $n\in\Z$), then $b\lambda\in F^{n-d}\La$.
\end{prop}

\pr -- Assume first that $\La$ is a successive extension of trivialisable $B_e$-representations, i.e. that there is $r\in\N$ and a filtration by sub-$B_e$-representations
$$0=\La_0\subset\La_1\subset\ldots\subset\La_{r-1}\subset\La_r=\La$$
such that each $\La_i/\La_{i-1}$ is trivialisable. We proceed by induction on $r$, the case $r=0$ being trivial. Assume $r\geq 1$. Setting $\La_{r-1}=\La'$ and choosing $U\in\cC^0(G_K)$ such that $\La_r/\La_{r-1}\simeq B_e\otimes_{\Q_p}U$, we may assume that we have a short exact sequence of $B_e$-representations
$$0\lra\La'\lra\La\lra B_e\otimes U\lra 0$$
and, using induction hypothesis, that we have an admissible filtration $(F^n\La')_{n\in\Z}$ of $\La'$ satisfying the required properties. Let $s:B_e\otimes U\to\La$ a $B_e$-linear section of the projection $\La\to B_e\otimes U$, so that we have a decomposition of the $B_e$-module $\La$ into a direct sum
$$\La=\La'\oplus s(B_e\otimes U)=\La'\oplus( B_e\otimes s(U))\ .$$
The map
$$\rho :G_K\times U\to \La'\ \ ,\ \ (g,u)\mapsto g(s(u))-s(g(u))$$
is continuous. Therefore, if $T$ is a $G_K$-stable lattice of $U$,  then $\rho(G_K\times T)$ is compact, hence bounded which implies (\cite{Schneider}, prop.5.6) that there exists $m\in\Z$ such that $\rho(G_K\times T)$, hence also $\rho(G_K\times U)$ is contained in $F^m\La'$.

If, for $n\in\Z$, we set
$$F^n\La= \left\lbrace\begin{matrix} F^n\La'\oplus (F^{n-m}B_e\otimes U)&\text{if }n\leq m,\cr
0&\text{if } n>m,\end{matrix}\right.$$
we see that $(F^n\La)_{n\in\N}$ is an admissible filtration satisfying the required properties.

-- In the general case, we chose a finite extension $K'$ of $K$ such that $\La$ is a successive extension of trivialisable $B_e$-representation of $G_{K'}$. Therefore we can find a $G_{K'}$-equivariant decreasing admissible filtration 
$$(F^n_{K'}\La)_{n\in\Z}$$
such that, if $n\in\Z$, then  $F^n_{K'}\La\in\cC(G_{K'})$ and that, if $b\in\fil^{-d}B_e$, for some $d\in\N$ and $\lambda\in F^n_{K'}\La$, then $b\lambda\in F^{n-d}_{K'}\La$. 

For each $n\in\Z$, denote $F^n\La$ the smallest sub-$\Q_p$-vector space of $\La$ containing $F^n_{K'}\La$ and stable under $G_K$. This is also the image of the obvious map $$\Q_p[G_K]\otimes_{\Q_p[G_{K'}]}F^n_{K'}\La\to \La\ .$$
If $h_1,h_2,\ldots,h_m$ is a system of representatives of $G_K/G_{K'}$ in $G_K$, this is also $\sum_{i=1}^mh_i(F^n\La_{K'})\subset \La$ which is still bounded and it is clear that the $(F^n\La)_{n\in\Z}$ satisfy the required properties. \qed

\begin{remk} \label{serrebe} {\rm 
We see immediately that $\rp_{B_e}(G_K)$ is a weak Serre subcategory of $\wh\cC(G_K)$. }\end{remk}

\subsection{Cohomology of coherent $\O_X$-modules}\label{coco}\

We denote $\ov B_{dR}$ the $B_e$-module $B_{dR}/B_e$. It is not of finite type but, as the cokernel of the inclusion $B_e\to B_{dR}$ which is a morphism of $\wh\cC(G_K)$, it can be viewed as an object of this category. The equalities $B_{dR}=B_e+B^+_{dR}$ and $\Q_p=B_e\cap B^+_{dR}$ imply that $\ov B_{dR}$, as an object of $\wh\cC(G_K)$,   can also be identified to $B^+_{dR}/\Q_p$.\\

If $\cF=(\cF_e,\cF^+_{dR},\iota_\cF)\in\oxm$. The map 
$$\cF_e\lra\cF_{dR}=B_{dR}\otimes_{B_e}\cF_e\ \ ,\ \ x\mapsto 1\otimes x$$
is injective, we use it to identify $\cF_e$ to a sub-$B_e$-module of $\cF_{dR}$ and we denote $\ov\cF_{dR}$ the quotient $\cF_{dR}/\cF_e$.

From \cite{FF}, prop. 8.2.3, we know:

\begin{prop} \label{cohcoh} For any $\cF\in\oxm$, we have $H^i(X,\cF)=0$ for $i\not\in\{0,1\}$ and 
$$\begin{matrix}\cF(X)=H^0(X,\cF)\not=0&\iff &\cF^{\geq 0}\not =0,\cr
H^1(X,\cF)\not=0&\iff&\cF^{<0}\not =0.\end{matrix}$$

Moreover, there is a canonical exact sequence of $\Q_p$-vector spaces
$$0\lra H^0(X,\cF)\lra \cF_e\oplus\cF^+_{dR}\overset{d_\cF}\lra \cF_{dR}\lra H^1(X,\cF)\lra 0\leqno{(1)}$$
(where $d_\cF(x,y)=\iota_\cF(y)-x$) which is functorial in $\cF$.
\end{prop}

 We have a commutative diagramme of $\Q_p$-vector spaces
$$\begin{matrix}&0&&0\cr&\da&&\da\cr&\cF_e&=&\cF_e\cr&\da&&\da\cr
0\lra H^0(X,\cF)\lra&\cF_e\oplus\cF^+_{dR}&\lra&\cF_{dR}& \lra H^1(X,\cF)\lra 0 \cr
&\da&&\da\cr& \cF^+_{dR}&\lra&\ov\cF_{dR}\cr&\da&&\da\cr&0&&0
\end{matrix}$$
whose columns and the two first lines are exact. Hence we have also an exact sequence
$$0\lra H^0(X,\cF)\lra \cF^+_{dR}\overset{\ov d_\cF}\lra \ov\cF_{dR}\lra H^1(X,\cF)\lra 0\leqno{(2)}$$
where $\ov d_\cF(y)$ is the image of $\iota_\cF(y)$ in $\ov\cF_{dR}$.\\

\subsection{Cohomology of coherent $\O_X[G_K]$-modules}\label{chch}\

We say that an almost $C_p$-representation is {\it effective} if this object of $\cC(G_K)$ is isomorphic to a sub-object of $\cC^\infty(G_K)$. We denote $\cC^{\geq 0}(G_K)$ the full subcategory of $\cC(G_K)$ whose objects are those which are effective.

\begin{prop} Let $f:W\to V$ a morphism of $\cC(G_K)$ with $W\in\cC^\infty(G_K)$ and $V\in\cC^{\geq 0}(G_K)$. Then the kernel of $f$ belongs to $\cC^\infty(G_K)$.
\end{prop}

\pr By assumption, there exists a monomorphism $g:V\to W'$ in $\cC(G_K)$ with $W'\in\cC^\infty(G_K)$. The kernel of $f$ is the same as the kernel of $gf:W\to W'$. As $W$ and $W'$ are in $\cC^\infty(G_K)$, so is this kernel.\qed

\begin{prop} Let $\cF\in\cM(G_K)$. Then $H^0(X,\cF)\in\cC^{\geq 0}(G_K)$.
\end{prop}

\pr  We see that $\cF_e$, $\cF^+_{dR}$ and $\cF_{dR}$ can be viewed as objects of the abelian category $\wh\cC(G_K)$. The inclusion $\cF_e\hookrightarrow \cF_{dR}$ is a morphism of this category, hence $\ov\cF_{dR}$ can also viewed as an object of $\wh\cC(G_K)$. The map $\ov d_\cF$ of the exact sequence (2) is obviously a morphism of this category, hence
$$H^0(X,\cF)=\ker\ov d_\F\ \text{ and }\ H^1(X,\cF)=\coker\ \ov d_\F$$
are objects of $\wh\cC(G_K)$.\\

For $m\in\N$, big enough, $\cF(-m)_{HN}$ has all its HN-slopes $<0$ and $H^0(X,\cF(-m)_{HN})=0$. But this is the kernel of the map 
$$\cF^+_{dR}(m)\lra \ov \cF_{dR}\ \ ,\ \ b\otimes t^m\mapsto t^mb\ (\mod\ \cF_e)\ .$$
Hence we have a commutative diagramme
$$\begin{matrix}0&\lra& 0&\lra&\cF^+_{dR}(m)&\lra&\ov\cF_{dR}\cr
&&\da&&\da&&\parallel\cr
0&\lra &H^0(X,\cF)&\lra & \cF^+_{dR}&\lra&\ov\cF_{dR}\end{matrix} $$
(the first non zero vertical arrow sends $b\otimes t^m$ to $t^mb$) whose lines are exact.
Therefore, the compositum $H^0(X,\cF)\to \cF^+_{dR}\to \cF^+_{dR}/t^m\cF^+_{dR}$ is injective and $H^0(X,\cF)$, subobject in $\wh\cC(G_K)$ of $\cF^+_{dR}/t^m.\cF^+_{dR}\in\cC^\infty(G_K)$ is in $\cC^{\geq 0}(G_K)$. \qed \\

\section{Hulls and construction of the functor $V\mapsto\cF_V$}\label{fv}

\subsection{Generalities}\

In what follows, $B_?$ is either $B_e$, $B^+_{dR}$ or $B_{dR}$.

\def\rpq{\rp_{B_?}(G_K)}

We know (remarks \ref{serrebe}  and \ref{serrebdr}) that $\rpq$ can be identified to a weak Serre's subcategory of $\wh\cC(G_K)$. We have "inclusions" of weak Serre subcategories
$$\begin{matrix} \rp_{B^+_{dR}}(G_K)\cr
&\searrow\cr&&\rp_{B_{dR}}(G_K)&\lra&\wh\cC(G_K)\cr&\nearrow\cr
\rp_{B_e}(G_K)\end{matrix}  $$

 Let $V$ be an almost $C_p$-representation. We say that {\it $V$ has a $B_?$-hull} if the functor
$$\rpq\lra \Q_p\text{-vector spaces}\ \ ,\ \ W\mapsto \Hom_{\wh\cC(G_K)}(V,W)$$
is representable, i.e. if there is a (necessarily unique up to unique isomorphism) pair $(V_?,\iota^V_?)$, with $V_?$ a $B_?$-representation and $\iota^V_?:V\to V_?$ a $G_K$-equivariant continuous $\Q_p$-linear map, such that, for all $B_?$-representation $W$, the map
$$ \Hom_{\rpq}(V_?,W)\lra \Hom_{\wh\cC(G_K)}(V,W)\ ,$$
induced by $\iota^V_?$, is bijective. 

When it is the case, we call $(V_?,\iota^V_?)$, or abusively $V_?$, the {\it $B_?$-hull of $V$}.

Our purpose is to show that such an hull always exists and to use these hulls to construct a functor
$$\cC(G_K)\lra \cM(G_K)\ \ ,\ \ V\mapsto \cF_V\ .$$

\begin{remk} {\rm Let $V$ be an almost $C_p$-representation and let $I_V$ the class of morphisms
$$\iota: V\to W_\iota$$
 of $\wh\cC(G_K)$ whose source is $V$ and target a $B_?$-representation.  With suitable conventions and abuses, to say that $V$ has a $B_?$-hull means that the directed inverse limit 
 $$V_?=\limp_{\iota\in I_V}W_\iota$$
 exists and that the $B_?$-module underlying this  "pro-$B_?$-representation of $G_K$" is of finite type.}
 \end{remk}\
 
 Restricted to the full subcategory of $\cC(G_K)$ of almost $C_p$-representations admitting a $B_?$-hull, the correspondence $V\mapsto V_?$ is obviously functorial.\\
 
 Let $V\in\cC(G_K)$ such that, with obvious notations,  $(V^+_{dR},\iota^{V,+}_{dR})$ exists, let $M\in\rp_{B_{dR}}(G_K)$ and $f:V\to M$ a morphism in $\wh\cC(G_K)$. We see that the sub $B^+_{dR}$-module $W$ of $M$ generated by $f(V)$ is an object of $\cC^\infty(G_K)$, hence there is a unique morphism (in $\wh\cC(G_K)$ or, in this case, equivalently in $\cC^\infty(G_K)$) $$g: V^+_{dR}\to W\subset M$$ such that $f=g\circ \iota^{V,+}_{dR}$ and we have
 $$\displaylines{\Hom_{\wh\cC(G_K)}(V,M)=\Hom_{\wh\cC(G_K)}(V^+_{dR},M)=\cr\Hom_{\wh\cC(G_K)}(B_{dR}\otimes_{B^+_{dR}}V^+_{dR},M)=\Hom_{\rp_{B_{dR}(G_K)}}(B_{dR}\otimes_{B^+_{dR}}V^+_{dR},M)\ .}$$
 Therefore $V_{dR}$ exists and can be identified to $B_{dR}\otimes_{B^+_{dR}}V^+_{dR}$.
 
The same argument applies to the case where $(V_e,\iota^V_e)$ exists. Hence we have:

\begin{prop} \label{vdr} Let $V\in\cC(G_K)$.

i) If $V^+_{dR}$ exists, then $V_{dR}$ exists and is, canonically and functorially, $B_{dR}\otimes_{B^+_{dR}}V^+_{dR}$,

ii) If $V_e$ exists, then $V_{dR}$ exists and is, canonically and functorially, $B_{dR}\otimes_{B_e}V_e$.
\end{prop}
\
 \begin{prop}\label{quotenv} Let $B_?$ as above and let $V$ be an almost $C_p$-representation of $G_K$ which has a $B_?$-hull $(V_?,\iota^V_?)$. Then:
 
 i)  The image of $\iota^V_?$ generates $V_?$ as a $B_?$-module.
 
 ii) If moreover
 $$0\lra V'\lra V\lra V''\lra 0$$
 is a short exact sequence in $\cC(G_K)$, then $V''$ has a $B_?$-hull which is the quotient of $V_?$ by the sub-$B_?$-module of $V_?$ generated by the image of $V'$.
 
 iii) In this situation, if $V'$ has a $B_?$-hull, then the sequence
 $$V'_?\lra V_?\lra V''_?\lra 0$$
 is exact.
 \end{prop}
 
 \pr i) Let $W_0$ be the sub $B_?$-module of $V_?$ generated by the image of $V$. As $B_?$ is noetherian, this is a $B_?$-module of finite type.  By the universal property of $V_?$, there is a unique morphism $\nu: V_?\to W_0$ such that the map $V\to W_0$ is $\nu\circ\iota^V_?$ and we see that $V_?=W_0\oplus\ker\nu$. The fact that ${\rm id}_{V_?}$ is the unique endomorphism of $V_?$ such that $\nu\circ\iota^V_?=\iota^V_?$ forces $\ker\nu$ to be $0$.
 
 ii) If $W$ is any $B_?$-representation, we have
 $$\displaylines{\Hom_{\wh\cC(G_K)}(V'',W)=\big\{ f\in\Hom_{\wh\cC(G_K)}(V,W)\mid f(V')=0\big\}=\cr \big\{f\in\Hom_{\rp_{B_?}(G_K)}(V_?,W)\mid f(\iota^V_?(V'))=0\big\}=\Hom_{\rp_{B_?}(G_K)}(V_?/B_?\iota^V_?(V'),W)}$$
 
 iii) Let $N$ be the kernel of the projection $V_?\to V''_?$. The image of $V'_?$ in $V_?$ is clearly contained in $N$.  As $N$ is the sub $B_?$-module generated by the image of $V'$, the map $V'_?\to N$ is surjective and 
 $$V'_?\lra V_?\lra V''_?\lra 0$$
 is exact.
  \qed\\

\subsection{Construction of trivialisable almost $C_p$-representations}\label{extens}\

A {\it trivialisation} of an almost $C_p$-representation $V$ is a short exact sequence
$$0\lra U\lra V\lra W\lra 0$$
in $\cC(G_K)$ with $U\in\cC^0(G_K)$ and $W\in\cC^\infty(G_K)$.

An almost $C_p$-representation is {\it trivialisable} if it admits a trivialisation.

If $V\in\cC(G_K)$, if $\wt f:V/U_+\to W/U_-$ is an almost isomorphism with $W\in\cC(G_K)$ and if $\wh V=W\times_{W/U_-}V$, we have, in $\cC(G_K)$, a commutative diagramme
$$\begin{matrix} &&&&0&&0\cr&&&&\da&&\da\cr&&&&U_+&=&U_+\cr&&&&\da&&\da\cr
0&\lra&U_-&\lra&\wh V&\lra&V&\lra&0\cr &&\parallel&&\da&&\da\cr
0&\lra&U_-&\lra&W&\lra&V/U_+&\lra&0\cr &&&&\da&&\da\cr&&&&0&&0\end{matrix}$$
and $V$ is a quotient of $\wh V$ which is trivialisable as it is an extension of $W$ by $U_+\in\cC^0(G_K)$.\\

Given $U\in\cC^0(G_K)$ and $W\in\cC^\infty(G_K)$, it is easy to construct all almost $C_p$-representations which are extensions of $W$ by $U$:

Recall that  $$B_{dR}=B_e+B^+_{dR}\ \ \text{ and }\ \ B_e\cap B^+_{dR}=\Q_p$$
and that we set
$$\wt B_{dR}=B_{dR}/B^+_{dR}=B_e/\Q_p\ .$$

Let $U$ be an object of $\co$  and $W$ an object of $\ci$. Tensoring the exact sequence
$$0\lra \Q_p\lra B_e\lra\wt B_{dR}\lra 0$$
 by $U$ we get a short exact sequence in $\wh\cC(G_K)$
$$0\lra U\lra B_e\otimes_{\Q_p}U\lra\wt B_{dR}\otimes_{\Q_p}U\lra 0$$
inducing a map
$$\begin{matrix}\delta_{U,W}&:& \Hom_{\wh\cC(G_K)}(W,\wt B_{dR}\otimes_{\Q_p}U)&\lra &{\rm Ext}^1_{\wh\cC(G_K)}(W,U)\cr
&&\parallel&&\parallel\cr
&&  \Hom_{\wh\cC^\infty(G_K)}(W,\wt B_{dR}\otimes_{\Q_p}U)&&{\rm Ext}^1_{\cC(G_K)}(W,U)\end{matrix}$$

We have (cf. \cite{Fo03}, prop. 3.7):

\begin{prop} \label{ext}    Let $U\in\cC^0(G_K)$ and $W\in\cC^\infty(G_K)$. The map
$$\delta_{U,W}: \Hom_{\wh\cC^\infty(G_K)}(W,\wt B_{dR}\otimes_{\Q_p}U)\lra {\rm Ext}^1_{\cC(G_K)}(W,U)$$
is an isomorphism.
\end{prop}

Hence if $V$ is a trivialisable almost $C_p$-representation and if 
$$0\lra U\lra V\lra W_0\lra 0\leqno{(T)}$$
is a trivialisation of $V$, there is a unique
$$\rho_T\in\Hom_{\wh\cC^\infty(G_K)}(W_0,\wt B_{dR}\otimes_{\Q_p}U)$$
such that the square
$$\begin{matrix} V&\lra&\mspace{-18mu} W_0\cr
\da&&\da\rho_T\cr B_e\otimes_{\Q_p}U&\lra&\wt B_{dR}\otimes_{\Q_p}U\end{matrix}$$
is cartesian.

\subsection{Construction of the  hulls}\

\begin{prop}\label{prophulls} Any almost $C_p$-representation $V$ has a $B_e$-hull $V_e$,  a $B^+_{dR}$-hull  $V^+_{dR}$ and a $B_{dR}$-hull $V_{dR}$. We have
$$\displaylines{V_{dR}=B_{dR}\otimes_{B_e}V_e=B_{dR}\otimes_{B^+_{dR}}V^+_{dR}\ ,\cr
\rm{rank}_{B^+_{dR}}V^+_{dR}={\rm rank}_{B_e}V_e=\dim_{B_{dR}}V_{dR}\geq h(V)\ }$$
 and equality holds when $V$ is trivialisable.
 
 Moreover

i) if $U\in\cC^0(G_K)$, then $U_e=B_e\otimes_{\Q_p}U$ and $U^+_{dR}=B^+_{dR}\otimes_{\Q_p}U$,

ii) if $W\in\cC^\infty(G_K)$, then $W_e=0$ and $W^+_{dR}=W$,

iii) if 
$$0\lra U\lra V\lra W_0\lra 0\leqno{(T)}$$
is a trivialisation of an almost $C_p$-representation $V$, 

a) the map $U_e=B_e\otimes_{\Q_p}U\to V_e$ is an isomorphism, 

b) we have a short exact sequence
$$0\lra B^+_{dR}\otimes_{\Q_p}U\lra V^+_{dR}\lra W_0\lra 0$$
More precisely, $V^+_{dR}$ is the fiber product $(B_{dR}\otimes_{\Q_p}U)\times_{\wt B_{dR}\otimes_{\Q_p}U} W_0$ (where $W_0\to \wt B_{dR}\otimes_{\Q_p}U$ is the map $\rho_T$).
\end{prop}

\pr From the proposition \ref{vdr}, we see that the existence of $V_e$ and $V^+_{dR}$ implies the existence of $V_{dR}$ and the equalities:
$$\displaylines{V_{dR}=B_{dR}\otimes_{B_e}V_e=B_{dR}\otimes_{B^+_{dR}}V^+_{dR}\ ,\cr
\rm{rank}_{B^+_{dR}}V^+_{dR}={\rm rank}_{B_e}V_e=\dim_{B_{dR}}V_{dR}\ .}$$

i) Let $U\in\cC^0(G_K)$. By adjunction, for any $B_e$-representation $\La$, we have
$$\Hom_{\wh\cC(G_K)}(U,\La)=\Hom_{\rp_{B_e}(G_K)}(B_e\otimes_{\Q_p}U,\La)$$
hence $U_e$ exists and is $B_e\otimes_{\Q_p}U$. Similarly, for any object $W_0\in\cC^\infty(G_K)$, we have
$$\Hom_{\cC(G_K)}(U,W_0)=\Hom_{\wh\cC^\infty(G_K)}(B^+_{dR}\otimes_{\Q_p}U,W_0)$$
hence $U^+_{dR}$ exists and is $B^+_{dR}\otimes_{\Q_p}U$.
 In particular, $\dim_{B_{dR}}U_{dR}=h(U)$.\\ 

ii)  Let $W\in\cC^\infty(G_K)$. For all $B_e$-representation $\La$, we have $\Hom_{\wh\cC(G_K)}(W,\La)=0$ (prop.\ref{wlambda}). Therefore $W_e$ exists and is $=0$. For any  $W_0\in\cC^\infty(G_K)$, we have $\Hom_{\cC(G_K)}(W,W_0)=\Hom_{\cC^\infty(G_K)}(W,W_0)$ (prop.\ref{propforget}) hence $W^+_{dR}$ exists and is $W$. In particular $\dim_{B_{dR}}W_{dR}=0=h(W)$.\\

iii) Let $V$ a trivialisable almost $C_p$-representation and 
$$0\lra U\lra V\lra W_0\lra 0 \leqno{(T)}$$ a trivialisation.\\

a)  Let $\La$ be a $B_e$-representation. The inclusion $U\to V$ induces a map
 $$\a:\Hom_{\wh\cC(G_K)}(V,\La)\lra\Hom_{\wh\cC(G_K)}(U,\La)\overset{\simeq}\lra\Hom_ {\rp_{B_e}(G_K)}(B_e\otimes_{\Q_p}U,\La)=\Hom_{\wh\cC(G_K)}(B_e\otimes_{\Q_p}U,\La)$$
(prop.\ref{forget} and prop.\ref{forgete}).  We have a cartesian square (\S\ref{extens})
 $$\begin{matrix} \mspace{-8mu}V&\lra&\mspace{-18mu} W_0\cr
\da\rho&&\da\rho_T\cr B_e\otimes_{\Q_p}U&\lra&\wt B_{dR}\otimes_{\Q_p}U\end{matrix}\leqno{(S)}$$
 and we may use $\rho$ to get a map
 $$\beta:\Hom_{\wh\cC(G_K)}(B_e\otimes_{\Q_p}U,\La)\lra\Hom_{\wh\cC(G_K)}(V,\La)$$
 Let $f\in\Hom_{\wh\cC(G_K)}(B_e\otimes_{\Q_p}U,\La)$ and $f'=\a(\b(f))$. If  $\sum b_i\otimes u_i\in B_e\otimes_{\Q_p}U$, we have
 $$f'(\sum b_i\otimes u_i)=\sum b_i(\b(f)(u_i))=\sum b_if(u_i)=f(\sum b_i\otimes u_i)$$
 as $f$ is $B_e$-linear, hence $f'=f$.
 
 Let $g\in\Hom_{\wh\cC(G_K)}(V,\La)$ and $g'=\a\b(g))$. If $u\in U$, as $\rho(u)=u$, we have
 $$g'(u)=\b(\a(g))(u)=\a(g)(u)=g(u)$$
 Hence $g'-g$ factors through a morphism in $\wh\cC(G_K)$
 $$W_0\lra \La$$
 which is necessarily $0$ (th.\ref{thalmost}), hence $g'=g$. Therefore we see that $\a$ is an isomorphism. It's implies that $V_e$ exists and is equal to $U_e=B_e\otimes_{\Q_p}U$.\\
 
 b) We want to show  that $V^+_{dR}$ exists and is equal to
  $$W_1=(B_{dR}\otimes_{\Q_p}U)\times_{\wt B_{dR}\otimes_{\Q_p}U} W_0\ .$$
  Using the cartesian square $(S)$ and the inclusion $B_e\otimes_{\Q_p}U\subset B_{dR}\otimes_{\Q_p}U$, we get a morphism  of $\wh\cC(G_K)$
  $$V\to W_1$$
  and we have a commutative diagramme in $\wh\cC(G_K)$
$$\begin{matrix} 0&\lra&U &\lra&V&\lra &W_0&\lra&0\cr
&&\da&&\da&&\parallel\cr
0&\lra&U^+_{dR}&\lra&W_1&\lra&W_0&\lra&0\end{matrix}\leqno{(\ast)}$$
whose lines are exact.

If $W$ is any $B^+_{dR}$-representation, we have a commutative diagramme
$$\begin{matrix}
0&\lra&\Hom(W_0,W)&\lra&\Hom(V,W)&\lra&\Hom(U,W)&\lra&{\rm Ext}^1(W_0,W)\cr
&&\parallel&&\da&&\da\simeq&&\parallel\cr
0&\lra&\Hom(W_0,W)&\lra&\Hom(W_1,W)&\lra&\Hom(U^+_{dR},W)&\lra&{\rm Ext}^1(W_0,W)\end{matrix}$$
(where all the $\Hom$ and ${\rm Ext}^1$ are computed in $\wh\cC(G_K)$)
which implies that $$\Hom_{\cC(G_K)}(V,W)\lra \Hom_{\cC(G_K)}(W_1,W)=\Hom_{\cC^\infty(G_K)}(W_1,W)$$ is an isomorphism. Hence $V^+_{dR}$ exists and is equal to $W_1$.\\

-- Finally, let $V$ be any object of $\cC(G_K)$. We can find an exact sequence
$$0\lra U\lra \wh V\lra V\lra 0$$
with $\wh V$ trivialisable. The existence of $\wh V_e$ and $\wh V^+_{dR}$ implies (prop.\ref{quotenv}) the existence of $V_e$ and $V_{dR}^+$. The exactness of the sequence
$$U_{dR}\lra\wh V_{dR}\lra V_{dR}\lra 0$$
implies that $$\dim_{B_{dR}}V_{dR}\geq \dim_{B_{dR}}\wh V_{dR}-\dim_{B_{dR}}U_{dR}= h(\wh V)-h(U)=h(V)\ .$$ \qed

\subsection{The functor $V\mapsto\cF_V$}\label{functorfv}
\

 For any almost $C_p$-representation $V$, denote
$$\iota_V:V^+_{dR}\lra V_{dR}=B_{dR}\otimes_{B_e}V_e$$
the natural map. It induces an isomorphism $B_{dR}\otimes_{B^+_{dR}}V^+_{dR}\to V_{dR}$. Therefore
$$\cF_V=(V^+_{dR},V_e,\iota_V)$$
is a coherent $\O_X[G_K]$-module. This construction is clearly functorial and we get an additive functor
$$\cC(G_K)\lra \cM(G_K)\ \ ,\ \ V\mapsto \cF_V\ .$$ 
From the universal properties of the functor $V\mapsto V^+_{dR}$ and $V\mapsto V_e$, we deduce the fact that $V\mapsto \cF_V$ is left adjoint to $\cF\mapsto\cF(X)$.

\section{The equivalence $\cM^{\geq 0}(G_K)\lra\cC^{\geq 0}(G_K)$}

\subsection{A characterisation of effective coherent $\O_X[G_K]$-modules}\

\begin{thm} \label{checm} The category $\cM^{\geq 0}(G_K)$ is the smallest strictly full subcategory of $\cM(G_K)$ containing $\cM^0(G_K)$ and $\cM^\infty(G_K)$ and stable under taking extensions and direct summands.
\end{thm}

\begin{lemma} Let $s$ be a positive rational number. There exists $\cG_s\in\cM^s(G_K)$ which is an extension of an object of $\cM^\infty(G_K)$ by an object of $\cM^0(G_K)$.
\end{lemma}

\noindent{\it The lemma implies the theorem:} As a subcategory of $\cM(G_K)$, the category $\cM^{\geq 0}(G_K)$ is obviously stable under taking extensions and direct summands. Hence, it suffices to show that any $\cF\in\cM^{\geq 0}(G_K)$ can be written as a direct summand of successive extensions of direct summands of objects which are extensions of an object of $\cM^\infty(G_K)$ by an object of $\cM^0(G_K)$. Using the Harder-Narasimhan filtration, it is enough to show that, if $\cF$ is semi-stable of slope $s\geq 0$, then $\cF$ is such a direct summand.

If $s=0$, then $\cF\in\cM^0(G_K)$ though, if $s=+\infty$, then $\cF\in\cM^\infty(G_K)$ and we may assume that $s$ is a positive rational number.

Let $\cG_s$ as in the lemma, so that we have a short exact sequence
$$0\lra \cG^0_s\lra\cG_s\lra\cG_s^\infty\lra 0$$
with $\cG^0_s\in\cM^0(G_K)$ ans $\cG^\infty_s\in\cM^\infty(G_K)$.
As $\cG_s$ is a vector bundle (it has no torsion), its dual $\cG_s^\vee$ is well defined and semi-stable of slope $-s$. Therefore $$\cF_0=\cF\otimes\cG_s^\vee$$ is semi-stable of slope $0$. We have a short exact sequence
$$0\lra \cF_0\otimes\cG^0_s\lra\cF_0\otimes\cG_s\lra\cF_0\otimes\cG_s^\infty\lra 0$$
and $\cF_0\otimes\cG_s$ is an extension of $\cF^0\otimes\cG_s^\infty\in\cM^\infty(G_K)$ by $\cF^0\otimes\cG_s^0\in\cM^0(G_K)$.

But, with obvious notations,
$$\cF_0\otimes\cG_s=\cF\otimes\cG_s^\vee\otimes\cG_s=\cF\otimes End(\cG_s)\ .$$
If $End^0(\cG_s)$ is the sub sheaf of elements of trace $0$ in $End(\cG_s)$, we have
$$End(\cG_s)=\O_X\oplus End^0(\cG_s)$$
hence $$\cF_0\otimes\cG_s=\cF\otimes (\O_X\oplus End^0(\cG_s))=\cF\oplus(\cF\otimes End^0(G_s))$$ and $\cF$ is a direct summand of $\cF^0\otimes\cG_s$. \qed\\

\noindent{\it Proof of the lemma:} We may assume $K=\Q_p$.  Recall (\cite{FF}, prop. 10.5.3, see also \cite{CF}, \S 5) that
\begin{itemize}
\item A {\it filtered $\ph$-module over $\Q_p$} is a pair $(D,\fil)$ consisting of

a) a {\it $\ph$-module over $\Q_p$}, i.e. a finite dimensional $\Q_p$-vector space $D$ equipped with an automorphism $\ph :D\to D$,

b) an exhausted and separated decreasing filtration $(\fil^nD)_{n\in\Z}$.

\item There is a fully faithful additive functor $$(D,\fil)\mapsto \cF_{D,\fil}$$ from the category of filtered $\ph$-modules over $\Q_p$ to the category of $G_{\Q_p}$-equivariant vector bundles over $X$ (the essential image consists of those equivariant vector bundles which are {\it crystalline}, i.e. those $\cF$'s such that the natural map
$$B_{cris}\otimes_{\Q_p}\big{(}B_{cris}\otimes_{B_e}\cF_e\big{)}^{G_K}\lra B_{cris}\otimes_{B_e}\cF_e$$
is bijective): we have 
$\cF_{D,\fil}=(\cF_{D,\fil,e},\cF^+_{D,\fil,dR})$ where 

-- $\cF_{D,\fil,e}$ is the $B_e$-module $(B_{cris}\otimes_{\Q_p}D)_{\ph=1}$ which implies that $$\cF_{D,\fil,dR}=B_{dR}\otimes_{B_e}\cF_{D,e}=B_{dR}\otimes_{\Q_p}D\ ,$$

-- $\cF_{D,\fil,dR}^+ =\fil^0(B_{dR}\otimes_{\Q_p}D)=\sum_{n\in\Z} \fil^{-n}B_{dR}\otimes \fil^n D$.

\end{itemize}\

 Set $s=d/h$ with $d,h$ positive integers, prime together.
 
  Consider the  $\ph$-module $D$ over $\Q_p$  whose underlying $\Q_p$-vector space is of dimension $h$, with $(e_r)_{r\in\Z/h\Z}$ as a basis and
 $$\ph(e_r)=\left\lbrace\begin{matrix} e_{r+1}&\text{if} & r+1\not =0,\cr
 p^{-d}e_0&\text{if}& r+1=0\end{matrix}\right.$$
 
 We equip $D$ with two distincts filtrations $\fil$ and $\fil_0$:
 $$\fil^nD = \left\lbrace\begin{matrix} D&\text{if}&n\leq 0,\cr 0&\text{if}&n>0\end{matrix}\right.
 \ \ ,\ \ Fil_0^nD=\left\lbrace\begin{matrix} D &\text{if}& n\leq -d,\cr \oplus_{r\not=0}\Q_pe_r&\text{if}&-d<n\leq 0,\cr 0&\text{if}&n>d\end{matrix}\right.$$
 
 Set $\cG_s=\cF_{D,\fil}$ and $\cG_s^0=\cF_{D,\fil_0}$. Both are coherent $\O_X[G_K]$-module of rank $h$. As the polynomial $X^h-p^{-d}$ is irreducible over $\Q_p$, the $\Q_p[\ph]$-module $D$ is irreducible which implies that $\cG_s$ and $\cG_s^0$ are stable, hence semi-stable. An easy computation shows that $\deg(\cG_s)=d$ though $\deg(\cG_s^0)=0$, hence $\cG_s$ is semi-stable of slope $d/h=s$ though $\cG_s^0$ is semi-stable of slope $0$, hence belongs to $\cM^0(G_K)$. We see that $\cG_{s,e}^0=\cG_{s,e}$ and that $(\cG_s^0)^+_{dR}\subset(\cG_s)_{dR}$. Therefore $\cG_s^0$ is a subobject of $\cG_s$ and the cokernel $\cG^\infty_s$ is torsion, i.e. belongs to $\cM^\infty(G_K)$.\qed
 
\subsection{Some properties of effective almost $C_p$-representations}\

Recall (\S \ref{conv}) that an exact subcategory of an abelian category is a strictly full subcategory containing $0$ and stable under extensions. For instance the previous theorem shows that $\cM^{\geq 0}(G_K)$ is an exact subcategory of $\cM(G_K)$.

\begin{thm} \label{structcp} Let $V\in\cC(G_K)$. The following are equivalents:

i) $V$ is effective (i.e. $V\in\cC^{\geq 0}(G_K)$),

ii) there is a finite extension $K'$ of $K$ contained in $\ov Q_p$ such that $V$, as an object of $\cC(G_{K'})$ is a successive extension of objects belonging either to $\cC^0(G_{K'})$ or to $\cC^\infty(G_{K'})$, 

iii) $V$ belongs to the smallest strictly full subcategory of $\cC(G_K)$ containing $\cC^0(G_K)$ and $\cC^\infty(G_K)$ and stable under taking extensions and direct summands.

Moreover $\cC^{\geq 0}(G_K)$ is an exact subcategory of $\cC(G_K)$.

\end{thm}

Before proving this theorem, let's state an other result. Recall (\S \ref{functorfv}) that, to any $V\in\cC(G_K)$, we associated the coherent $\O_X[G_K]$-module
$$\cF_V=(V^+_{dR},V_e,\iota_V)\ .$$
We have 
$$(\cF_V)^+_{dR}=V^+_{dR}\ \ ,\ \ (\cF_V)_e=V_e\ \text{ and }\ \iota_{\cF_V}=\iota_V\ .$$

Therefore, if we set $\ov V_{dR}=\ov{\cF_V}_{dR}=V_{dR}/V_e$, we have (cf. \S \ref{chch})  an exact sequence
$$0\lra H^0(X,\cF_V)\lra V^+_{dR}\overset{\ov\iota_V}\lra \ov V_{dR}\lra H^1(X,\cF_V)\lra 0\leqno{(C)}$$
(where $\ov\iota_V=\ov\iota_{\cF_V}$ is the compositum of $\iota_V$ with the projection $V_{dR}\to V_{dR}/V_e$) and, as $V\subset V_e$ is injective, the image of $V$ in $V^+_{dR}$ is contained in $\cF_V(X)=H^0(X,\cF_V)$.

\begin{prop} \label{cCplus} Let $V\in\cC^{\geq 0}(G_K)$. Then,

i) we have $h(V)\geq 0$ and $\dim_{B_{dR}}V_{dR}=h(V)$,

ii) we have $V\in\cC^\infty(G_K)\iff h(V)=0$,

iii) the sequence
$$0\lra V\lra V^+_{dR}\overset{\ov\iota_V}\lra \ov V_{dR}\lra 0$$
 is exact, 
 
 iv) the map $V\to H^0(X,\cF_V)$ is bijective and $\cF_V\in\cM^{\geq 0}(G_K)$.

Moreover, the restriction to $\cC^{\geq 0}(G_K)$ of the four functors $\cC(G_K)\to \wh\cC(G_K)$
$$V\mapsto V^+_{dR}\ \ ,\ \ V\mapsto V_e\ \,,\ \ V\mapsto V_{dR} \ \ ,\ \ V\mapsto \ov V_{dR}$$
and of the functor 
$$\cC(G_K)\lra\cM(G_K)\ \ ,\ \ V\mapsto\cF_V$$
are exact.
\end{prop}

\noindent {\it Proof of the theorem and beginning of the proof of the proposition}:  

For any $V\in\cC^{\geq 0}(G_K)$, we denote $d_V$ the infimum of the $d(W)$'s for all $W\in\cC^\infty(G_K)$ such that $V$ is isomorphic to a subobject of $W$ (note that $d(V)\leq d_V$). 

Denote $\cK$ the set of finite extensions $L$ of $K$ contained in $\ov\Q_p$. For any $L\in\cK$, let $\cC^?(G_L)$ the full subcategory of $\cC(G_L)$ whose objects can be written as a successive extension of objects belonging either to $\cC^0(G_L)$ or to $\cC^\infty(G_L)$.

We now show the assertion (i) of the proposition and the implication (i)$\implies$(ii) of the theorem, i.e. that, if $V\in\cC^{\geq 0}(G_K)$, then 
$$\dim_{B_{dR}}V_{dR}=h(V)\ \text{ (which implies that $h(V)\geq 0$) and there exists $K'\in\cK$ such that }V\in \cC^?(G_{K'})\ :$$
 We proceed by induction on $d_V$, the case $d_V=0$ being trivial.

Let $V\subset W$ an embedding of $V$ into an object $W\in\cC^\infty(G_K)$ satisfying $d(W)=d_V>0$. We can find (cf. prop.\ref{drdeviss})  $K_1\in\cK$  and a $G_{K_1}$-stable sub-$B^+_{dR}$-module $W'$ of $W$ of length $1$. Setting $W''=W/W'$, $V'=V\cap W'$ and denoting $V''$ the image of $V$ in $W''$, we get a commutative diagramme in $\cC(G_{K_1})$
$$\begin{matrix} 0&\lra &V'&\lra&V&\lra&V''&\lra&0\cr
&&\da&&\da&&\da\cr
0&\lra &W'&\lra&W&\lra& W''&\lra&0
\end{matrix}$$
whose rows are exact and vertical arrows are injective which implies that $V'$ and $V''$ belong to $\cC^{\geq 0}(G_{K_1})$. We have $d(V')\leq d(W')=1$. From corollary \ref{vw}, we get that either $d(V')=1$ in which case $V'=W'$ or $d(V')=0$ which implies that $V'\in\cC^0(G_{K_1})$.

-- If $V'=W'$, we have $h(V')=0$ and  $(V')^+_{dR}=W'$ hence $V'_{dR}=0$,

-- If $V'\in\cC^0(G_{K_1})$, we have $h(V')=\dim_{\Q_p}V'$ and $V'_{dR}=B_{dR}\otimes_{\Q_p}V'$.

In both cases, we have $\dim_{B_{dR}}V'_{dR}=h(V')$. 
By induction, we have $\dim_{B_{dR}}V''_{dR}=h(V'')$. The exactness of the sequence
$$V'_{dR}\lra V_{dR}\lra V''_{dR}\lra 0$$
implies that $$\dim_{B_{dR}}V_{dR}\leq \dim_{B_{dR}}V'_{dR}+\dim_{B_{dR}}V''_{dR}= h(V')+h(V'')=h(V)\ ,$$
hence, as $\dim_{B_{dR}}V_{dR}\geq h(V)$ (prop.\ref{prophulls}), we get $\dim_{B_{dR}}V_{dR}=h(V)$, i.e the assertion (i) of the proposition.\\

Also by induction,  as $V''$ belongs to $\cC^{\geq 0}(G_{K_1})$, there is $K'\in\cK$ containing $K_1$ such that $V''\in\cC^?(G_{K'})$. Then $V$, as a representation of $G_{K'}$, is an extension of $V''$ by either an object of $\cC^\infty(G_{K'})$ (if $d(V')=1$) or by an object of $\cC^0(G_{K'})$ (if $d(V')=0$). In both cases, 
$V$ belongs to $\cC^?(G_{K'})$.  

Therefore, given $V\in\cC^{\geq 0}(G_K)$, there is $K'\in\cK$ and a filtration of $V$ by subobjects in $\cC^+(G_{K'})$
$$0=V_0\subset V_1\subset\ldots\subset V_{r-1}\subset V_r=V$$
such that, if $i=1,2,\ldots ,r$, then $V_i/V_{i-1}$ belongs either to $\cC^0(G_{K'})$ or to $\cC^\infty(G_{K'})$. This proves the implication (i)$\implies$(ii) of the theorem. \\

In particular, we have $h(V)=\sum_{i=1}^rh(V_i/V_{i-1})$ which is $>0$ unless $h(V_i/V_{i-1})=0$ for all $i$, which means that $V_i/V_{i-1}$ belongs to $\cC^\infty(G_K)$. As $\cC^\infty(G_K)$ is stable under taking extensions, we get the equivalence
$$h(V)=0\iff V\in\cC^{\geq 0}(G_K)$$
which is the assertion (ii) of the proposition.\\

The implication (ii)$\implies$(iii) of the theorem is obvious: If $V$ satisfies (ii),  the induced representation $\Q_p[G_K]\otimes_{\Q_p[G_{K'}]}V$ belongs to 
$\cC^?(G_K)$ and $V$ is a direct summand of this representation.\\

As a full subcategory of $\cC(G_K)$, the category $\cC^{\geq 0}(G_K)$ is obviously stable under taking direct summands. Hence,  we see that the implication (iii)$\implies$(i) of the theorem  and the fact that $\cC^{\geq 0}(G_K)$ is an exact subcategory of $\cC(G_K)$ result from the following:

\begin{lemma} Assume we have a short exact sequence in $\cC(G_K)$
$$0\lra V_0\lra V_1\lra V_2\lra 0\leqno{(1)}$$
with $V_2\in\cC^{\geq 0}(G_K)$ and $V_0$ belonging either to $\cC^0(G_K)$ or to $\cC^\infty(G_K)$. Then $V_1\in\cC^{\geq 0}(G_K)$ and the sequence
$$0\lra V^+_{0,dR}\lra V^+_{1,dR}\lra V^+_{2,dR}\lra 0$$
is exact.
\end{lemma}

\noindent{\it Proof of the lemma}: 

-- Assume first that $V_0$ belongs to $\cC^0(G_K)$: we have a commutative diagramme
$$\begin{matrix} 0&\lra &V_0&\lra&V_1&\lra& V_2&\lra&0\cr
&&\da&&\da&&\da\cr
&&V^+_{0,dR}&\lra&V^+_{1,dR}&\lra& V^+_{2,dR}&\lra&0
\end{matrix}$$
whose rows are exact, the maps $V_0\to V^+_{0,dR}$ and $V_2\to V^+_{2,dR}$ being injective. 
I claim that the map $V^+_{0,dR}\to V^+_{1,dR}$ is injective. As $V^+_{0,dR}=B^+_{dR}\otimes_{\Q_p}V_0$ is a torsion free $B^+_{dR}$-module, it is enough to check that $V_{0,dR}\to V_{1,dR}$ is injective. If it were not true, we would have
$$\dim_{B_{dR}}V_{1,dR}< \dim_{B_{dR}}V_{0,dR}+\dim_{B_{dR}}V_{2,dR}=h(V_0)+h(V_2)=h(V_1)\ .$$
As we have (prop.\ref{prophulls})
$\dim_{B_{dR}}V_{1,dR}\geq h(V_1)$, 
this can't happen. This forces $V_1\to V^+_{1,dR}$ to be also injective, hence $V_1\in\cC^{\geq 0}(G_K)$.\\

-- Assume now $V_0$ belongs to $\cC^\infty(G_K)$. As the sequence (1) almost splits (prop.\ref{almostsplit}), we can find an extension $S$ in $\cC^0(G_K)$ of $V_2$ by a $U\in\cC^0(G_K)$ such that $V_1=V_0\oplus_{U}S$. By what we just saw, $S\in\cC^{\geq 0}(G_K)$ and we have a commutative diagramme
$$\begin{matrix} 0&\lra &U&\lra&S&\lra& V_2&\lra&0\cr
&&\da&&\da&&\da\cr
0&\lra&U^+_{dR}&\lra&S^+_{dR}&\lra& V^+_{2,dR}&\lra&0
\end{matrix}$$
whose line are exacts and vertical arrows are injective. 

We also have a commutative diagramme
$$\begin{matrix} 0&\lra &U&\lra&W\oplus S&\lra& V_1&\lra&0\cr
&&\da&&\da&&\da\cr
&&U^+_{dR}&\lra&W\oplus S^+_{dR}&\lra& V^+_{1,dR}&\lra&0
\end{matrix}$$
(the map $U\to W\oplus S$ send $u$ to $(u,-u)$) whose rows are exact and the two first vertical arrows are injective.

The injectivity of $U^+_{dR}\to S^+_{dR}$ implies the injectivity of $U^+_{dR}\to W\oplus S^+_{dR}$. To finish the proof we only need to show that the map $V_1\to V^+_{1,dR}$ is injective or, with obvious identifications, that inside of $W\oplus S^+_{dR}$, we have
$$U^+_{dR}\cap(W\oplus S)= U\ .$$
Assume $(w,s)\in W\oplus S$ belongs to $U^+_{dR}$. This implies that $s\in S\cap U^+_{dR}$ which is $U$ as the map $V_2\to V^+_{2,dR}$ is injective. We then need $w=-s$ and $(w,s)$ is the image of $-s\in U$. \qed\\

\noindent{\it Proof of the exactness of the functors $V\mapsto V_{dR}$, $V\mapsto V_e$ and $V\mapsto \ov V_{dR}$}:
If $$0\lra V'\lra V\lra V''\lra 0$$ is a short exact sequence in $\cC^{\geq 0}(G_K)$, we know that the sequences
$$\begin{matrix}
V'_e&\lra &V_e&\lra& V''_e&\lra & 0\cr
V'_{dR}&\lra &V_{dR}&\lra& V''_{dR}&\lra & 0\end{matrix}$$
are exact. As
$$\dim_{B_{dR}}V_{dR}=h(V)=h(V')+h(V'')=\dim_{B_{dR}}V'_{dR}=\dim_{B_{dR}}V''_{dR}$$
the map $V'_{dR}\to V_{dR}$ must be injective and the functor $V\mapsto V_{dR}$ is exact.

As the $B_e$-modules $V'_e$, $V_e$ and $V''_e$ are torsion free and as 
$${\rm rank}_{B_e}(V'_e)=\dim_{B_{dR}}V'_{dR}\ ,\ {\rm rank}_{B_e}(V_e)=\dim_{B_{dR}}V_{dR}\ , {\rm rank}_{B_e}(V''_e)=\dim_{B_{dR}}V''_{dR}$$
the same argument shows the exactness of $V\mapsto V_e$.

We then have a commutative diagramme
$$\begin{matrix} &&0&&0&&0\cr &&\da&&\da&&\da\cr
0&\lra&V'_e&\lra&V_e&\lra&V''_e&\lra&0\cr &&\da &&\da&&\da\cr
0&\lra&V'_{dR}&\lra&V_{dR}&\lra&V''_{dR}&\lra&0\cr &&\da &&\da&&\da\cr
0&\lra&\ov V'_{dR}&\lra&\ov V_{dR}&\lra&\ov V''_{dR}&\lra&0\cr &&\da &&\da&&\da\cr
&&0&&0&&0 \end{matrix}$$ 
whose three columns and the two first rows are exact. This implies the exactness of the last row.\\

\begin{lemma} Let 
$$0\lra V'\lra V\lra V''\lra 0$$
a short exact sequence in $\cC^{\geq 0}(G_K)$. 
Assume the sequences $$0\lra V'\lra (V')^+_{dR}\lra \ov V'_{dR}\lra 0$$
and $$0\lra V''\lra (V'')^+_{dR}\lra \ov V''_{dR}\lra 0$$
are exact. Then the sequences 
$$0\lra V\lra V^+_{dR}\lra \ov V_{dR}\lra 0$$
and $$0\lra (V')^+_{dR}\lra V^+_{dR}\lra (V'')^+_{dR}\lra 0$$
are exact.
\end{lemma}

\noindent{\it Proof of the lemma}: We have a commutative diagramme:
$$\begin{matrix} &&0&&0&&0\cr &&\da&&\da&&\da\cr
0&\lra&V'&\lra&V&\lra&V''&\lra&0\cr &&\da &&\da&&\da\cr
0&\lra&(V')^+_{dR}&\lra&V^+_{dR}&\lra&(V'')^+_{dR}&\lra&0\cr &&\da &&\da&&\da\cr
0&\lra&\ov V'_{dR}&\lra&\ov V_{dR}&\lra&\ov V''_{dR}&\lra&0\cr &&\da &&\da&&\da\cr
&&0&&0&&0 \end{matrix}$$ 
whose first and third rows are exact. By assumption, the first and the third columns are also exact. We also now that, except may be in $(V')^+_{dR}$, the second line is exact and, as $V\in\cC^{\geq 0}(G_K)$, that the map $V\to V^+_{dR}$ is injective. By diagramme chasing, we get the fact that the second line and the second column  are also exact. \qed \\

\noindent{\it End of the proof of the proposition}:

We first prove (iii), i.e. for all $V\in\cC^{\geq 0}(G_K)$, the exactness of the sequence
$$0\lra V\lra V^+_{dR}\overset{\ov\iota_V}\lra \ov V_{dR}\lra 0$$

a) If $V\in \cC^\infty(G_K)$, as $V^+_{dR}=V$ and $V_{dR}=\ov V_{dR}=0$, exactness is obvious.

b) If $V\in\cC^0(G_K)$, this sequence can be rewritten 
$$0\lra V\lra B^+_{dR}\otimes_{\Q_p}V\lra \ov B_{dR}\otimes_{\Q_p}V\lra 0$$ 
and exactness is deduced by tensoring with $V$ from the exactness of
$$0\lra \Q_p\lra B^+_{dR}\lra \ov B_{dR}\lra 0$$
(recall that $B_{dR}=B_e+B^+_{dR}$ though $\ov B_{dR}=B_{dR}/B_e$ and that $B_e\cap B^+_{dR}=\Q_p$).

c) In general, we proceed by induction on  the smallest integer $r_V$ such that there is $K'\in\cK$ with the property that $V$ is a successive extension of $r_V$  objects belonging either to $\cC^0(G_{K'})$ or to $\cC^\infty(G_{K'})$. Replacing $K$ by $K'$ if necessary, we may assume $K'=K$. We just proved it's OK if $r_V=1$. Assume $r_V\geq2$, so that we can find a short exact sequence in $\cC^{\geq 0}(G_K)$
$$0\lra V'\lra V\lra V''\lra 0$$
with $r_{V'}$ and $r_{V''}<r_V$. Then, by induction, the sequences
$$\begin{matrix} 
0&\lra &V'&\lra& (V')^+_{dR}&\lra& \ov{V'}_{dR}&\lra& 0\cr
0&\lra &V''&\lra& (V'')^+_{dR}&\lra& (\ov{ V''})_{dR}&\lra& 0\end{matrix}$$
are exact and the result follows from the two assertions of the previous lemma.\\

From the exact sequence $(C)$, we see that $V=H^0(X,\cF_V)$ and that $H^1(X,\cF_V)=0$ hence that $\cF\in\cM^{\geq 0}(G_K)$, which proves (iv).\\

We are left to prove the exactness  of the functors $V\mapsto V^+_{dR}$ and $V\mapsto \cF_V$, i.e. that, if 
$$0\lra V'\lra V\lra V''\lra 0$$
is a short exact sequence in $\cC^{\geq 0}(G_K)$, then the sequences 
$$0\lra (V')^+_{dR}\lra V^+_{dR}\lra (V'')^+_{dR}\lra 0$$
and
$$0\lra\cF_{V'}\lra\cF_V\lra\cF_{V''}\lra 0$$
are exact. As we now know the assertion (iii) of the proposition, the exactness of the first sequence is a consequence of the previous lemma. Finally, we see that exactness of the second is equivalent to the exactness of 
$$0\lra (V')^+_{dR}\lra V^+_{dR}\lra (V'')^+_{dR}\lra 0$$
and of 
$$0\lra V'_e\lra V_e\lra V''_e\lra 0$$
and we are done. \qed\\

\begin{prop} \label{artin} Let $V\in\cC(G_K)$. Any decreasing sequence of subobjects of $V$
$$V_1\supset V_2\supset\ldots\supset V_n\supset V_{n+1}\supset\ldots$$
 is stationary.
 \end{prop}
 
 \pr Chose $\wh V\in\cC^{\geq 0}(G_K)$ such that $V$ is a quotient of $\wh V$. For all $n\in\N$, set $$\wh V_n=\wh V\times_VV_n\ .$$ The $\wh V_n$ form a decreasing sequence of subobject of $\wh V$ and, for all $n\in\N$, we have a canonical isomorphism $\wh V_n/\wh V_{n+1}\simeq V_n/V_{n+1}$. In particular 
 $$V_{n+1}=V_n\iff \wh V_{n+1}=\wh V_n\ .$$
 Replacing $V$ by $\wh V$ and the $V_n$'s by the $\wh V_n$'s if necessary we ay assume that $V$, therefore also the $V_n$'s are in $\cC^{\geq 0}$.

 As $d(V_{n+1})\leq d(V_n)$ and $d(V_n)\geq 0$, there is an integer $m$ such that $d(V_n)=d(V_{n+1})$ for $n\geq m$.
 
 For $n\geq m$, we have $d(V_n/V_{n+1})=0$, hence $V_n/V_{n+1}\in\cC^0(G_K)$ and, if we set $h_n=\dim_{\Q_p}(V_n/V_{n+1})$ $(\in\N)$, we have $h(V_{n+1})=h(V_n)-h_n$. As $V_{n+1}\in\cC^{\geq 0}(G_K)$, we have $h(V_{n+1})\geq 0$. Therefore, there is an integer $m'\geq m$ such that $h_n=0$ if $n\geq m'$. This implies that $V_{n+1}=V_n$.\qed \\
 
 \begin{remk} {\rm On the other hand, there are objects of $\cC(G_K)$ which admit non stationary increasing sequences of subobjects. For instance, it is easy to see that $C_p$ contains infinitely many subobjects belonging to $\cC^0(G_K)$. From that, one can constructs non stationary increasing sequences
 $$V_0\subset V_1\subset\ldots\subset V_n\subset V_{n+1}\subset\ldots$$
 of subobjects of $C_p$ belonging to $\cC^0(G_K)$.}
 \end{remk}

\subsection{The main result}\label{mainresult}\

We may consider the functors
$$\cM^{\geq 0}(G_K)\lra \cC^{\geq 0}(G_K)\ \ ,\ \ \cF\mapsto \cF(X)$$
and
$$\cC^{\geq 0}(G_K)\lra \cM^{\geq 0}(G_K)\ \ ,\ \ V\mapsto \cF_V\ .$$

\begin{thm}\label{mainthm}  The functor 
$$ \cM^{\geq 0}(G_K)\lra \cC^{\geq 0}(G_K)\ \ ,\ \ \cF\mapsto \cF(X)$$
is an equivalence of exact categories and 
$$\cC^{\geq 0}(G_K)\lra \cM^{\geq 0}(G_K)\ \ ,\ \ V\mapsto \cF_V$$
is a quasi-inverse.
 \end{thm}
 
\pr -- As the functor $V\mapsto\cF_V$ is left adjoint to $\cF\mapsto\cF(X)$ (\S \ref{functorfv}), we are reduced to check that \begin{enumerate}
\item if $V\in\cC^{\geq 0}(G_K)$, the map $V\to\cF_V(X)$ coming from adjonction is an isomorphism,
\item if $\cF\in\cM^{\geq 0}(G_K)$, the map $\cF_{\cF_V(X)}\to\cF$ coming from adjonction is an isomorphism,
\item if
$$0\lra V'\lra V\lra V''\lra 0$$
is a short exact sequence of $\cC^{\geq 0}(G_K)$, the sequence
$$0\lra\cF_{V'}\lra\cF_V\lra\cF_{V''}\lra 0$$
is exact,
\item if
$$0\lra\cF'\lra\cF\lra\cF''\lra 0$$
is a short exact sequence of $\cM^{\geq 0}(G_K)$, the sequence
$$0\lra\cF'(X)\lra\cF(X)\lra\cF''(X)\lra 0$$
is exact.
\end{enumerate}\

(1) and (3) have already been proved (prop.\ref{cCplus}) and (4) results from the fact that, if $\cF'\in\cM^{\geq 0}(G_K)$, then $H^1(X,\cF')=0$ (prop.\ref{cohcoh}). 

Let's prove (2): Let $\cM$ the full subcategory of $\cM^{\geq 0}(G_K)$ whose objects are those $\cF$'s for which $\cF_{\cF_V(X)}\to\cF$ is an isomorphism. It is obviously stable under taking direct summands. By exactness of the functors $\cF\to\cF(X)$ and $V\mapsto\cF_V$, it is stable under extensions. It contains $\cM^0(G_K)$ and $\cM^\infty(G_K)$. Then theorem \ref{checm} implies that $\cM=\cM^{\geq 0}(G_K)$. \qed\\

\section{From $\cM(G_K)$ to $\cC(G_K)$ and conversely}

\subsection{Some general nonsense}

Let $\cA$ be an abelian category and $\cB$ be an exact subcategory of $\cA$. Recall  (cf. eg \cite{Lau}, \S 1.1)  that one can define the derived category of bounded complexes of $\cB$ that we denote $D^b_{\cA}(\cB)$:  in the triangulated category $\cK^b(\cB)$ of bounded complexes of $\cB$ up to homotopies, the full subcategory $\cN$ of bounded acyclic complexes (in $\cB$) form a null system and we set
$$\cD^b_{\cA}(\cB)=\cK^b(\cB)/\cN\ .$$
\

Let $\cA$ be an abelian category, $\cB$  an exact subcategory of $\cA$ and $\cD$ a strictly full subcategory of $\cB$ which is a Serre's subcategory of $\cA$ (hence $\cD$ is abelian), 

--  We say that {\it the exact embedding $\cB\hra\cA$ is  left big with respect to $\cD$} if,
\begin{enumerate} \item any quotient in $\cA$ of an object of $\cB$ belongs to $\cB$,
\item for any object $A$ of $\cA$, one can find a short exact sequence
$$0\lra A\lra B\lra D\lra 0$$
of $\cA$ with $B$ an object of $\cB$ and $D$ an object of $\cD$.
\end{enumerate}
 
 -- We say that {\it the exact embedding $\cB\hra\cA$ is  right big with respect to $\cD$} if $\cB^{\rm op}\hra\cA^{\rm op}$ is left big with respect to $\cD^{\rm op}$ which amounts to requiring that,
 \begin{enumerate} \item any subobject in $\cA$ of an object of $\cB$ belongs to $\cB$,
 \item for any object $A$ of $\cA$, one can find a short exact sequence
$$0\lra D\lra B\lra A\lra 0$$
of $\cA$ with $B$ an object of $\cB$ and $D$ an object of $\cD$.
\end{enumerate}

We say that an exact embedding $\cB\hra\cA$ is {\it left big} (resp. {\it right big}) if one can find
a Serre's subcategory  $\cD$ of $\cA$ contained in $\cB$ such that $\cB\hra\cA$ is left big (resp. right big) with respect to $\cD$.

\begin{prop} \label{nonsense} Let $\cB\hra\cA$ an exact embedding which is either left big or right big. Then the natural functor
$$D^b_\cA(\cB)\lra D^b(\cA)$$
is an equivalence of triangulated categories.
\end{prop}

It is a formal consequence of the more precise following statement:

\begin{prop} Let $\cB\hra\cA$ be an exact embedding and $\cD$ a Serre's subcategory of $\cA$ contained in $\cB$  such that  $\cB\hra\cA$ is left big (resp. right big) with respect to $\cD$ and let $A^\bullet$ a bounded complex of $\cA$. Then,

i) there exists a short exact sequence of bounded complexes of $\cA$
$$0\to A^\bullet\to B^\bullet\to D^\bullet\to 0\ \ (\text{resp. }\ 0\to D^\bullet\to B^\bullet\to A^\bullet\to 0\ )$$
with $B^\bullet$ a bounded complex of $\cB$ and $D^\bullet$ an acyclic complex of $\cD$,

ii) if 
$$0\to A^\bullet\to B^{'\bullet}\to D^{'\bullet}\to 0\ (\text{resp. }\ 0\to D^{'\bullet}\to B^{'\bullet}\to A^\bullet\to 0\ )$$
is an other short exact sequence of the same kind, there exists a a third short exact sequence of the same kind
$$0\to A^\bullet\to B^{''\bullet}\to D^{''\bullet}\to 0\ (\text{resp. }\ 0\to D^{''\bullet}\to B^{''\bullet}\to A^\bullet\to 0\ )$$ 
together with morphisms of complexes 
$$B^\bullet\to B^{''\bullet}\text{ and } B^{'\bullet}\to B^{''\bullet}\ (\text{resp. }B^{''\bullet}\to B^\bullet\text{ and } B^{''\bullet}\to B^{'\bullet})$$
such that the diagramme$$\begin{matrix}&&A^\bullet\cr &\swarrow&\da&\searrow\cr
B^{'\bullet}&\to&B^{''\bullet}&\la&B^\bullet\end{matrix}\  \ \big{(}\text{ resp. }\
\begin{matrix}B^{'\bullet}&\to&B^{''\bullet}&\la&B^\bullet\cr &\searrow&\da&\swarrow\cr&&A^\bullet
\end{matrix}\ \big{)}$$
is commutative.
\end{prop}

\noindent{\it Proof of } (i): It is enough to prove it in the case were the strict embedding is right big. Assume it is the case. To prove (i), by induction, we are reduced to prove:

\begin{lemma} Let $r\in\Z$ and let 
$$0\to D_r^\bullet\to B_r^\bullet\to A^\bullet\to 0$$
a short exact sequence of bounded complexes of $\cA$. Assume that $D_r^\bullet$ is an acyclic complex of $\cD$,  that $D_r^n=0$ for $n\geq r$ and that $B_r^n$ is an object of $\cB$ for all $n<r$. Then, there exists a short exact sequence of bounded complexes of $\cA$
$$0\to D_{r+1}^\bullet\to B_{r+1}^\bullet\to A^\bullet\to 0$$
where $D_{r+1}^\bullet$ is an acyclic complex of $\cD$  with $D_{r+1}^n=0$ for $n\geq r+1$ and $B_{r+1}^n$ an object of $\cB$ for all $n<r+1$. 
\end{lemma}

\noindent {\it Proof of the lemma}: We can identify $B_r^n$ to $A^n$ for $n\geq r$. Granted to right bigness of $\cB\hra\cA$, we can find a short exact sequence
$$0 \lra D\lra B\lra A^r\lra 0$$
with $B$ an object of $\cB$ and $D$ an object of $\cD$. Set 
$$B_{r+1}^n=\left\lbrace\begin{matrix} B_r^n&\text{for}& n<r-1,\cr
	B_r^{r-1}\times_{A_r}B&\text{for}& n=r-1,\cr
	B&\text{for}&n=r,\cr
	B_r^n=A^n&\text{for}& n>r.\end{matrix}\right.$$
We have a commutative diagramme 
$$\begin{matrix}&&&0&&0\cr&&&\da&&\da\cr&&&D&=&D\cr&&&\da&&\da\cr
\ldots&\to B_{r+1}^{r-2}&\to &B_{r+1}^{r-1}&\to &B_{r+1}^{r}&\to& B_{r+1}^{r+1}&\to&\ldots\cr &\parallel&&\da&&\da&&\parallel\cr
\ldots&\to B_{r}^{r-2}&\to &B_{r}^{r-1}&\to &B_{r}^{r}&\to &B_{r}^{r+1}&\to&\ldots\cr
&&&\da&&\da\cr&&&0&&0\cr\end{matrix}$$
whose rows are complexes and columns are exact. Hence we have a quasi-isomorphism $B_{r+1}^\bullet\to B_r^\bullet$. Moreover $B_{r+1}^n$ is an object of $\cB$ for all $n<r+1$ (for $n=r-1$, this is due to the fact that $B_{r+1}^{r-1}$ is a subobject of $B_r^{r-1}\oplus B$ which belongs to $\cB$).

The compositum $$B_{r+1}^\bullet\lra B_r^\bullet\lra A^\bullet$$ is a surjective morphism of complexes which is a quasi-isomorphism. Then the kernel $D_{r+1}^\bullet$ is acyclic. As it is the complex
$$\ldots \to D_r^{r-3}\to D_r^{r-2}\to D_{r+1}^{r-1}\to D\to 0\to\ldots \to 0\to\ldots$$
we see that $D_{r+1}^n=0$ for $n\geq r+1$ and that all the $D_{r+1}^n$ belong to $\cD$ (for $n=r-1$, this is due to the fact that we have a short exact sequence
$$0\lra D'' \lra D_{r+1}^{r-1}\lra D\lra 0$$
with $D''=\coker (D_r^{r-3}\to D_r^{r-2})\in\cD$, hence, as $D_{r+1}^{r-1}$ is an extension in $\cA$ of $D\in\cD$ by $D''\in\cD$, it belongs to $\cD$).\qed\\

\noindent{\it Proof of }(ii): We just take, for each $n\in\Z$, the fiber product
$$B^{''n}= B^{'n}\times_{A^n}B^n\ .$$
For each $n$, we have an exact sequence
$$0\lra D^{''n}\lra B^{''n}\lra A^n\lra 0$$
with $D^{''n}=D^{'n}\oplus D^n$ and all the required properties are obviously fulfilled. \qed

\subsection{The equivalence of triangulated categories}\

\begin{thm} \label{thtrig} The equivalence of categories of theorem \ref{mainthm}
extends uniquely to an equivalence of triangulated categories
$$D^b(\cM(G_K))\lra D^b(\cC(G_K))\ .$$
\end{thm}

\pr Uniqueness is obvious. 

Recall (\S \ref{mainresult}) that $\cM^{\geq 0}(G_K)$ is an exact subcategory of $\cM(G_K)$ though $\cC^{\geq 0}(G_K)$ is an exact subcategory of $\cC(G_K)$.\\

 -- The category $\cM^\infty(G_K)$ is a Serre's subcategory of $\cM(G_K)$ contained in $\cM^{\geq 0}(G_K)$ and any quotient $\cF''$ in $\cM(G_K)$ of an object $\cF$ of $\cM^{\geq 0}(G_K)$ is  in $\cM^{\geq 0}(G_K)$ 
$$\text{(as }\cF\in\cM^{\geq 0}(G_K)\iff H^1(X,\cF)=0\implies H^1(X,\cF'')=0\iff \cF''\in\cM^{\geq 0}(G_K)\text{).}$$
 If $\cF\in\cM(G_K)$, for all $n\in\N$, as, for all $n\in\N$, the HN-slopes of $\cF(n)_{HN}$ are the $s+n$ for $s$ describing the HN-slopes of $\cF$ (cf. \S \ref{twists}), for $n\gg0$, we have $\cF(n)_{HN}\in\cM^{\geq 0}(G_K)$.

Tensoring with $\cF$ the short exact sequence (\S \ref{twists})
$$0\lra\O_X\lra \O_X(n)_{HN}\lra (0,B_n(-n))\lra 0$$
we get a short exact sequence
$$0\lra\cF\lra \cF(n)_{HN}\lra (0,\cF^+_{dR}\otimes_{B^+_{dR}}B_n(-n))\lra 0$$
As $\cF(n)_{HN}$ belongs to $\cM^{\geq 0}(G_K)$ and $(0,\cF^+_{dR}\otimes_{B^+_{dR}}B_n(-n))$ belongs to $\cM^\infty(G_K)$, it shows that the exact embedding $\cM^{\geq 0}(G_K)\to \cM(G_K)$ is left big with respect to $\cM^\infty(G_K)$. Therefore (prop.\ref{nonsense})  the natural functor
$$D_{\cM(G_K)}^b(\cM^{\geq 0}(G_K))\lra D^b(\cM(G_K))$$
is an equivalence of triangulated categories.\\

-- Similarly, the category  $\cC^0(G_K)$ is a Serre's subcategory of $\cC(G_K)$ contained in $\cC^{\geq 0}(G_K)$ and any subobject in $\cC(G_K)$ of an object of $\cC^{\geq 0}(G_K)$ belongs to $\cC^{\geq 0}(G_K)$.

 Let $V\in\cC(G_K)$ and chose an almost  isomorphism $V/U_+\simeq W/U_-$ with 
 $W\in\cC^\infty(G_K)$ (cf. \S\ref{secalmostcp}). Set 
 $$\hat V=V\times_{W/U_-}W$$
 (where the map $V\to W/U_-$ is the compositum of the projection $V\to V/U_+$ with the isomorphism $V/U_+\to W/U_-$). 
 
 We have a diagramme
$$\begin{matrix} &&&&0\cr&&&&\da\cr&&&& U_+\cr &&&&\da\cr
0&\lra& U_-&\lra&\wh V&\lra&V&\lra&0\cr
&&&&\da\cr&&&& W\cr &&&&\da\cr&&&& 0
\end{matrix}$$
whose line and column are exacts. The column shows that $\wh V\in\cC^{\geq 0}(G_K)$ and, therefore, the line shows that $V$ is a quotient of an object of $\cC^{\geq 0}(G_K)$ by an object of $\cC^0(G_K)$. In other words, the exact embedding $\cC^{\geq 0}(G_K)\to \cC(G_K)$ is right big with respect to $\cC^0(G_K)$.  Hence (prop.\ref{nonsense}) the natural functor
$$D_{\cC(G_K)}^b(\cC^{\geq 0}(G_K))\lra D^b(\cC(G_K))$$
is an equivalence of triangulated categories.\\

As the equivalence $\cM^{\geq 0}(G_K)\overset{\simeq}\lra \cC^{\geq 0}(G_K)$ is an equivalence of exact categories, it extends uniquely to an equivalence of triangulated categories
$$D^b_{\cM(G_K)}(\cM^{\geq 0}(G_K))\lra D^b_{\cC(G_K)}(\cC^{\geq 0}(G_K))\ .$$
\

-- It is now clear that there is a unique equivalence of triangulated categories
$$D^b(\cM(G_K))\lra D^b(\cC(G_K))$$
such that the square
$$\begin{matrix} D^b_{\cM(G_K)}(\cM^{\geq 0}(G_K))&\lra &D^b_{\cC(G_K)}(\cC^{\geq 0}(G_K))\cr
\da&&\da\cr D^b(\cM(G_K))&\lra &D^b(\cC(G_K))\end{matrix}
$$
is commutative and that this equivalence extends the equivalence of the theorem \ref{mainthm}. \qed

\subsection{The equivalence $\cM^{<0}(G_K)\to \cC^{<0}(G_K)$}\label{cmocgo}

We say that a coherent $\O_X[G_K]$-module is {\it co-effective} if all its HN slopes are $<0$. We saw (prop.\ref{cohcoh}) that $\cF\in\cM(G_K)$ is co-effective if and only if $H^0(X,\cF)=0$. The full subcategory of $\cM(G_K)$ whose objects are co-effective is $\cM^{<0}(G_K)$ and is stable under taking subobjects and extensions.

Any $\cF\in\cM(G_K)$ as a biggest quotient $\cF^{<0}$ belonging to $\cM^{<0}(G_K)$and the sequence 
$$0\lra \cF^{\geq 0}\lra\cF\lra\cF^{<0}\lra 0$$
is exact.\\

We say that an almost $C_p$-representation $V$ is {\it co-effective} if, for all $W\in\cC^\infty(G_K)$, we have $\hom_{\cC(G_K)}(V,W)=0$. We denote $\cC^{<0}(G_K)$ the full subcategory of $\cC(G_K)$ whose objects are co-effective. It is obviously stable undertaking quotients and extensions.

\begin{prop} Let $V$ be an almost $C_p$-representation. The following are equivalent

i)  $V$ is co-effective,

ii) $V^+_{dR}=0$,

iii) $\cF_V=0$. 

If this is the case, then
$$V_e=V_{dR}=\ov V_{dR}=0\ .$$
\end{prop}

\pr The equivalence (i)$\iff$(ii) results from the universal property of $V^+_{dR}$ and (ii)$\iff$(iii) is trivial. If $V^+_{dR}=0$, we have $V_{dR}=B_{dR}\otimes_{B^+_{dR}}V^+_{dR}=0$, hence also $V_e=0$ as the map $V_e\to V_{dR}$ is injective and therefore $\ov V_{dR}=V_{dR}/V_e=0$.\qed\\

\begin{prop} \label{vcoef} Let $V\in\cC(G_K)$. The set of subobjects of $V$ in $\cC(G_K)$ belonging to $\cC^{<0}(G_K)$ has a biggest element $V^{<0}$ and the set of quotients of $V$ in $\cC(G_K)$ belonging to $\cC^{\geq 0}(G_K)$ as a biggest element $V^{\geq 0}$. Moreover $V^{<0}$ (resp. $V^{\geq 0}$) is the kernel (resp. the image) of the natural map $V\to V^+_{dR}$.
The sequence
$$0\lra V^{<0}\lra V\lra V^{\geq 0}\lra 0$$
is exact.
\end{prop}

\pr  If $V'$ and $V''$ are subobjects of $V$ belonging to $\cC^{<0}(G_K)$, we see that $V'+V''$ also. Hence to show the existence of $V^{<0}$ it is enough to show that any increasing sequence
$$V_0\subset V_1\subset\ldots\subset V_n\subset V_{n+1}\subset\ldots$$
of subobjects of $V$ belonging to $\cC^{<0}(G_K)$ is stationary. As the sequence of the integers $d(V_n)$ is increasing and bounded by $d(V)$, there exists $m\in\N$ such that $d(V_n)=d(V_m)$ for all $n\geq m$. For such an $n$, we have $d(V_{n+1}/V_n)=0$, hence $V_{n+1}/V_n\in\cC^0(G_K)$. This implies $V_{n+1}=V_n$ as, otherwise, the compositum of the projection of $V_{n+1}$ onto $V_n$ with the injective map 
$$V_{n+1}/V_n\lra B^+_{dR}\otimes_{\Q_p}(V_{n+1}/V_n)\ \ ,\ \ v\mapsto 1\otimes v$$
would be a non zero morphism from $V_{n+1}$ to an object of $\cC^\infty(G_K)$.\\

If $\ov V'$ and $\ov V''$ are quotients of $V$ belonging to $\cC^{\geq 0}(G_K)$, then the image of $V\to \ov V'\oplus\ov V''$ also (as it is a subobject of $\ov V'\oplus\ov V''\in\cC^{\geq 0}(G_K)$). Hence to show the existence of $V^{\geq 0}$ it suffices to show that any sequence
$$\ldots\lra \u V_{n+1}\lra \u V_n\lra\ldots\lra \u V_1\subset \u V_0$$
of quotients of $V$ (belonging to $\cC^{<0}(G_K)$) is stationary. If $\wt V_n$ is the kernel of the projection $V\to\ov V_n$, the sequence $(\wt V_n)_{n\in\N}$ is a decreasing sequence of objects of $\cC(G_K)$, hence is stationary (prop.\ref{artin}), therefore the sequence of the $\ov V_n$'s also.
\\

Set $V_0=\ker(V\to V^+_{dR})$. We obviously have $V^{<0}\subset V_0$ and to show the equality it is enough to show that $V_0\in\cC^{<0}(G_K)$. Otherwise, we could find a non zero morphism $f:V^0\to W$ with $W\in\cC^\infty(G_K)$. Let $V_1=\ker f$ and consider the short exact sequence
$$0\lra V_0/V_1\lra V/V_1\lra V/V_0\lra 0$$
As $V_0/V_1$ injects into $W$, it belongs to $\cC^{\geq 0}(G_K)$. As $V/V_0$ injects into $V^+_{dR}$, it also belongs to $\cC^{\geq 0}(G_K)$. Therefore, as $\cC^{\geq 0}(G_K)$ is stable under extensions, $V/V_1\in\cC^{\geq 0}(G_K)$. Hence the sequence
$$0\lra (V_0/V_1)^+_{dR}\lra (V/V_1)^+_{dR}\lra (V/V_0)^+_{dR}\lra 0$$
is exact. As obviously $(V/V_1)^+_{dR}=(V/V_0)^+_{dR}=V^+_{dR}$, it contradicts the fact that, as $V_0/V_1$ is a non zero object of $\cC^{\geq 0}(G_K)$, we have $(V_0/V_1)^+_{dR}\not = 0$.\\

Let $V_2=\im(V\to V^+_{dR})$. As the map $V_2\to V^+_{dR}$ is injective, $V_2$ belongs to $\cC^{\geq 0}(G_K)$ and is, therefore a quotient of $V^{\geq 0}$. The kernel $V_3$ of the projection $V^{\geq 0}\to V_2$ belongs also to $\cC^{\geq 0}(G_K)$ (as this category is stable under taking subobjects) and we have an exact sequence in $\cC^{\geq 0}(G_K)$
$$0\lra V_3\lra V^{\geq 0}\lra V_2\lra 0$$
Therefore the sequence 
$$0\lra V^+_{3,dR}\lra V^{\geq 0,+}_{dR}\lra V^+_{2,dR}\lra 0$$
is also exact. 

 As $V^{\geq 0}$ is a quotient of $V$, we see that $V^{\geq 0,+}_{dR}$ is a quotient of $V^+_{dR}$. But clearly $V^+_{2,dR}=V^+_{dR}$. Therefore $V^{\geq 0,+}_{dR}=V^+_{dR}$ and $V^+_{3,dR}=0$. As $V_3\in\cC^{\geq 0}(G_K)$, this implies $V_3=0$, hence $V^{\geq 0}=V_2$.\\
 
 The exactness of
 $$0\lra V^{<0}\lra V\lra V^{\geq 0}\lra 0$$
 is now clear. \qed

\begin{remks} {\rm i) To any $V\in\cC(G_K)$, we just associated the canonical short exact sequence
$$0\lra V^{<0}\lra V\lra V^{\geq 0}\lra 0$$
 It's worth to compare with the canonical short exact sequence 
$$0\lra \cF^{\geq 0}\lra \cF\lra \cF^{<0}\lra 0$$
associated to any $\cF\in\cM(G_K)$.

ii) We know that, for any $\cF\in\cM(G_K)$, the natural map $\cF^{\geq 0}(X)\to \cF(X)$ is an isomorphism.The two previous propositions together imply that, for any $V\in\cC(G_K)$, the natural map $\cF_V\mapsto \cF_{V^{\geq 0}}$ is an isomorphism. In particular, $\cF_V$ always belongs to $\cM^{\geq 0}(G_K)$.}
\end{remks}

It is clear that $\cM^{<0}(G_K)$ is an exact subcategory of $\cM(G_K)$, though $\cC^{<0}(G_K)$ is an exact subcategory of $\cC(G_K)$.

\begin{prop}\label{equivcoeff} If $\cF\in\cM(G_K)$, then $H^1(X,\cF)\in\cC^{<0}(G_K)$ and the map $H^1(X,\cF)\mapsto H^1(X,\cF^{<0})$ is an isomorphism.

Moreover, the  functor
$$\cM^{<0}(G_K)\lra \cC^{<0}(G_K)\ \ ,\ \ \cF\mapsto H^1(X,\cF)$$
is an equivalence of exact categories.
\end{prop}

\pr

If $\cF\in\cM(G_K)$, we may find a short exact sequence in $\cM(G_K)$
$$0\lra\cF\lra \cF^0\lra\cF^1\lra 0$$
with $\cF^0\in\cM^{\geq 0}(G_K)$.
As $H^1(X,\cF^0)=0$, we see that $H^1(X,\cF)$ is the cokernel of $H^0(X,\cF^0)\to H^0(X,\cF^1)$, hence belongs to $\cC(G_K)$. 

We know that $H^1(X,\cF)$ is a quotient of $\ov\cF_{dR}$ therefore also of $\cF_{dR}$. If $f:H^1(X,\cF)\to W$ were a non zero morphism of $\cC(G_K)$ with $W\in\cC^\infty(G_K)$, the compositum
$\cF_{dR}\to H^1(X,\cF)\to W$ would be a non zero morphism in $\wh\cC^\infty(G_K)$ and, therefore, would be $B^+_{dR}$-linear. As multiplication by $t$ is invertible in $\cF_{dR}$ and nilpotent in $W$, the map must be $0$ which implies that $H^1(X,\cF)\in\cC^{<0}(G_K)$.\\

If $A$ is an object of an abelian category and $d\in\Z$, we denote $A[d]$ the bounded complex in $\cA$ which is $A$ in degree $-d$ and $0$ elsewhere.

Denote $$\Gamma:D^b(\cM(G_K))\lra D^b(\cC(G_K))\ \text{ (resp. } \Delta:D^b(\cC(G_K))\lra D^b(\cM(G_K))\text{ )}$$
the functor extending $\cF\mapsto \cF(X)$ (resp. $V\mapsto\cF_V$). If $\cF\in\cM^{<0}(G_K)$ and if 
$$0\lra\cF\lra \cF^0\lra\cF^1\lra 0$$
is as above (observe that $\cF^0\in\cM^{\geq 0}(G_K)\implies\cF^1\in\cM^{\geq 0}(G_K)$), we see that (with obvious conventions)
 $$\Gamma(\cF[0])= \Gamma(\cF^0\to\cF^1)=(H^0(X,\cF^0)\to H^0(X,\cF^1))=H^1(X,\cF)[-1]$$
 (as $\cF\in\cM^{<0}(G_K)$ and $\cF^0\in\cM^{\geq 0}(G_K)$, the sequence
 $$0\lra H^0(X,\cF^0)\lra H^0(X,\cF^1)\lra H^1(X,\cF)\lra 0$$
 is exact).\\
 
 Let $V\in\cC^{<0}(G_K)$. We can find a short exact sequence in $\cC(G_K)$
 $$0\lra V^0\lra V^1\lra V\lra 0$$
 with $V^1\in\cC^{\geq 0}(G_K)$ which implies $V^0\in\cC^{\geq 0}(G_K)$. With obvious conventions, we have
 $$\Delta(V[-1])=\Delta(V^0\to V^1)=(\cF_{V^0}\to\cF_{V^1})=\cF[0]$$
with $\cF$ the kernel of $\cF_{V^0}\to\cF_{V^1}$ (as $V\in\cC^{<0}(G_K)$, we have $V^+_{dR}=V_e=0$ which implies that
$$\cF_{V^0}=(V^{0+}_{dR},V^0_e,\iota_{V^0})\lra \cF_{V^1}=(V^{1+}_{dR},V^1_e,\iota_{V^1})$$ is an epimorphism).

We have a commutative diagramme
$$\begin{matrix} &&0&&0&&0\cr&&\da&&\da&&\da\cr
0&\lra&H^0(X,\cF)&\lra&V^0&\lra&V^1\cr&&\da&&\da&&\da\cr
0&\lra&\cF^+_{dR}&\lra&V^{0+}_{dR}&\lra&V^{1+}_{dR}\cr
&&\da&&\da&&\da\cr 0&\lra&\cF_e&\lra&V^0_e&\lra&V^1_e
\end{matrix}$$
whose rows and columns are exact. The injectivity of $V^0\to V^1$ implies that $H^0(X,\cF)=0$, i.e. that $\cF\in\cM^{<0}(G_K)$.\\

Finally, we see that, if we view 

-- $\cM^{<0}(G_K)$ as the full subcategory of $D^b(\cM(G_K))$ whose objects are of the form $\cF[0]$ with $\cF\in\cM^{<0}(G_K)$,

-- $\cC^{<0}(G_K)$ as the full subcategory of $D^b(\cC(G_K))$ whose objects are of the form $V[-1]$ with $V\in\cC^{<0}(G_K)$,

then $\Gamma$ induces the required equivalence of categories. \qed

\subsection{$t$-Structures and hearts}\label{heart}\

The functors
$$\Gamma: D^b(\cM(G_K))\to D^b(\cC(G_K))\ \text{ and }\ \Delta : D^b(\cC(G_K))\to D^b(\cM(G_K))$$
are as in the proof of the previous proposition.

Let $(D_\cM^{\leq 0},D_\cM^{\geq 0})$ be the canonical $t$-structure on $D^b(\cM(G_K))$: we see that $D_\cM^{\leq 0}$ (resp. $D_\cM^{\geq 0}$) is the full subcategory of $D^b(\cM(G_K))$ whose objects are those $\cF^\bullet$ such that $H^i(\cF^\bullet)=0$ for $i>0$ (resp. $i<0$).
Therefore if we denote $\Gamma(D_\cM^{\leq 0})$  (resp.  $\Gamma(D_\cM^{\geq 0})$) the essential image under $\Gamma$ of  $D_\cM^{\leq 0}$ (resp.  $D_\cM^{\geq 0}$), we see that $(\Gamma(D_\cM^{\leq 0}),\Gamma(D_\cM^{\geq 0}))$ is a $t$-structure on $D^b(\cC(G_K))$ whose heart $\Gamma(D_\cM^{\leq 0})\cap \Gamma(D_\cM^{\geq 0})$ is an abelian category equivalent via $\Delta$ to $\cM(G_K)$.

Similarly, let $(D_\cC^{\leq 0},D_\cC^{\geq 0})$ the canonical $t$-structure on $D^b(\cC(G_K))$: hence $D_\cC^{\leq 0}$ (resp . $D_\cC^{\geq 0}$) is the full subcategory of $D^b(\cC{G_K)}$ whose objects are those $V^\bullet$ such that $H^i(V^\bullet)=0$ for $i>0$ (resp. $i<0$).
Therefore if we denote $\Delta(D_\cC^{\leq 0})$  (resp.  $\Delta(D_\cC^{\geq 0})$) the essential image under $\Delta$ of  $D_\cC^{\leq 0}$ (resp.  $D_\cC^{\geq 0}$), we see that $(\Delta(D_\cC^{\leq 0}),\Delta(D_\cC^{\geq 0}))$ is a $t$-structure on $D^b(\cM(G_K))$ whose heart $\Delta(D_\cC^{\leq 0})\cap \Delta(D_\cC^{\geq 0})$ is an abelian category equivalent via $\Gamma$ to $\cC(G_K)$.

\begin{prop} \label{tstruc} i) $\Gamma(D_\cM^{\geq 0})$ (resp. $\Gamma(D_\cM^{\leq 0})$) is the full subcategory of $D^b(\cC(G_K))$ whose objects are those $V^\bullet$'s such that $H^r(V^\bullet)=0$ for $r<0$ and $H^0(V^\bullet)\in \cC^{\geq 0}(G_K)$ (resp. $H^r(V^\bullet)=0$ for $r>1$ and $H^1(V^\bullet)\in\cC^{<0}(G_K)$).

ii) $\Delta(D_\cC^{\geq 0})$ (resp. $\Delta(D_\cC^{\leq 0})$) is the full subcategory of $D^b(\cM(G_K))$ whose objects are those $\cF^\bullet$'s such that $H^r(\cF^\bullet)=0$ for $r<-1$ and $H^{-1}(\cF^\bullet)\in \cM^{< 0}(G_K)$ (resp. $H^r(\cF^\bullet)=0$ for $r> 0$ and $H^0(\cF^\bullet)\in\cM^{\geq 0}(G_K)$).
\end{prop}

\pr Let's prove that the description of $\Gamma(D_\cM^{\geq 0})$ is correct (the proof of the three other statements are similar):

Any object $\u\cF$ of $D_\cM^{\geq 0}$ can be represented by a bounded complex $\cF^\bullet$ such that $\cF^ i=0$ for $i<0$. From the fact that, for any $\cF\in\cM(G_K)$, one can find a short exact sequence
$$0\lra \cF\lra\cF_0\lra\cF_1\lra 0$$
with $\cF_0,\cF_1\in\cM^{\geq 0}(G_K)$ and the fact that any quotient, in $\cM(G_K)$, of an object of $\cM^{\geq 0}(G_K)$ still belongs to $\cM^{\geq 0}(G_K)$, one easily deduces that the complex $\cF^\bullet$ is quasi-isomorphic to a bounded complex $\cF_0^\bullet$ with $\cF_0^r=0$ for $r<0$ and $\cF_0^r\in \cM^{\geq 0}(G_K)$ for all $r\in\N$. Therefore $\Gamma(\u\cF)$ is represented by the bounded complex 
$$\ldots\lra 0\lra\ldots\lra 0\lra \cF_0^0(X)\lra \cF_0^1(X)\lra\ldots\lra\cF_0^r(X)\lra\cF_0^{r+1}(X)\lra\ldots$$
whose all term belong to $\cC^{\geq 0}(G_K)$. In particular, as $\cC^{\geq 0}(G_K)$ is stable under taking subobjects in $\cC(G_K)$, we see that $\Gamma(\u\cF)$ belongs to the full subcategory $D_{\cC,\cM}^{\geq 0}$ of $D^b(\cC(G_K))$ whose objects are those $\u V$'s such that $H^r(\u V)=0$ for $r<0$ and $H^0(\u V)\in\cC^{\geq 0}(G_K)$.

Conversely, any object $\u V$ of $D_{\cC,\cM}^{\geq 0}(G_K)$ can be represented by a complex $V_0^\bullet$ such that $V_0^r=0$ for $r<0$ and  that the kernel of $V^0\to V^1$ belongs to $\cC^{\geq 0}(G_K)$. Using the fact that, for any $V\in \cC(G_K)$ one can find a short exact sequence in $\cC(G_K)$
$$0\lra V_1\lra V_0\lra V\lra 0$$
with $V_1,V_0\in\cC^{\geq 0}(G_K)$,
one easily deduces that the complex $V_0^\bullet$ is quasi-isomorphic to a bounded complex $V^\bullet$ with $V^r=0$ for $r<0$ and $V^r\in\cC^{\geq 0}(G_K)$ for  $r>0$. 

We have a short exact sequence
$$0\lra (V^0_{d=0})\lra V^0\lra dV^0\lra 0$$
The inclusion $dV^0\subset V^1$ implies that $dV^0\in\cC^{\geq 0}(G_K)$.
 As $V^0_{d=0}=H^0(V^\bullet)$, we have $(V^0)_{d=0}\in \cC^{\geq 0}(G_K)$. We know that $\cC^{\geq 0}(G_K)$, as a full subcategory of $\cC(G_K)$, is stable under extension. Therefore $V^0\in\cC^{\geq 0}(G_K)$.

As all the $V^r$'s belong to $\cC^{\geq 0}(G_K)$, we see that $\Delta(\u V)$ is represented by the bounded complex 
$$\ldots \lra 0\lra\ldots\lra 0\lra \cF_{V^0}\lra\cF_{V^1}\lra \ldots\ \lra \cF_{V^r}\lra\cF_{V^{r+1}}\lra\ldots$$
hence belong to $D^{\geq 0}_\cM$. \qed

\subsection{Torsion pairs in $\cM(G_K)$ and in $\cC(G_K)$}\

The langage of torsion pairs  (cf. \cite{HRS}, chap.1)
 is very convenient to give an explicite description of the way to go from $\cM(G_K)$ to $\cC(G_K)$ and conversely. The results of this subsection and of the next one are independent of those of the previous one and give an other proof of the description of the heart of the $t$-structures we considered (prop.\ref{tstruc}).\\
 
 Recall ({\it loc. cit.}) that a {\it torsion pair }in an abelian category $\cA$ is a pair $t=(\cA^+,\cA^-)$ of full subcategories of $\cA$ containing $0$ such that:
\begin{enumerate}
\item If $B$ is an object of $\cA^+$ and $C$ is an object of $\cA^-$, then $\Hom_\cA(B,C)=0$,

\item for any object $A$ of $\cA$, there is a short exact sequence in $\cA$
$$0\lra A^+\lra A\lra A^-\lra 0$$
with $A^+\in\Ob(\cA^+)$ and $A^-\in\Ob(\cA^-)$.
\end{enumerate}

Condition (1) implies that the exact sequence of (2) is unique up to a unique isomorphism and that the correspondences $A\mapsto A^+$ and $A\mapsto A^-$ are functorial.\\

\def\hea{{\rm He}(\cA^+,\cA^-)}
\def\heoa{{\rm He}_0(\cA^+,\cA^-)}
\def\heb{{\rm He}(\cB^+,\cB^-)}
\def\heob{{\rm He}_0(\cB^+,\cB^-)}

We define the {\it heart $\cA^t$} of $t$ as the full subcateogry of the derived category $D^b(\cA)$ whose objects are those $A^\bullet$ such that
$$H^{-1}(A^\bullet)\in\Ob(\cA^-)\ , \ H^0(A^\bullet)\in\Ob(\cA^+)\ , \ H^n(A^\bullet)=0 \text{ if } n\not\in\{-1,0\}\ .$$
\

\begin{prop} Let $t=(\cA^+,\cA^-)$ be a torsion pair in an abelian category $\cA$.
Consider the full subcategories $D^{\leq 0}=D^{\leq 0}_t(\cA)$ and $D^{\geq 0}=D^{\geq 0}_t(\cA)$ of $D=D^b(\cA)$ defined by
\begin{enumerate}
\item ${\rm Ob}(D^{\leq 0})=\{A^\bullet\in\Ob(D^b(\cA))\mid H^1(A^\bullet)\in\Ob(\cA^+)\text{ and } H^n(A^\bullet)=0 \text{ for all } n>1\}$,
\item ${\rm Ob}(D^{\geq 0})=\{A^\bullet\in\Ob(D^b(\cA))\mid H^0(A^\bullet)\in\Ob(\cA^-)\text{ and } H^n(A^\bullet)=0 \text{ for all } n<0\}$.
\end{enumerate} 

Then $(D^{\leq 0},D^{\geq 0})$ is a $t$-structure on $D$ whose heart is $\cA^t$.
\end{prop}

\pr To show that $(D^{\geq 0},D^{\leq 0})$ is a $t$-structure, we have to check (cf. \cite{KS}, def. 10.1.1) that (with standard notations)
\begin{enumerate}
\item $D^{\leq -1}\subset D^{\leq 0}$ and $D^{\geq 1}\subset D^{\geq 0}$,
\item $\Hom_D(X,Y)=0$ for $X\in\Ob(D^{\leq 0})$ and $Y\in\Ob(D^{\geq 1})$,
\item For any $X\in\Ob(D)$, there exists a distinguished triangle $X_0\lra X\lra X_1\underset{+1}\lra$ in $D$ with $X_0\in\Ob(D^{\geq 0})$ and $X_1\in\Ob(D^{\geq 1})$.
\end{enumerate}
 (1) is obvious.  (2) is clear as, if $f:X\to Y$ with $X\in\Ob(D^{\leq 0})$ and $Y\in\Ob(D^{\geq 1})$, we have $H^n(f)=0$ for $n\leq 0$ (as $H^n(Y)=0$), for $n>1$ (as $H^n(X)=0$) and for $n=1$ (as $H^1(X)\in\Ob(\cA^+)$ and $H^1(Y)\in\Ob(\cA^-)$).
Let's check (3): we have $H^1(X)= X^1_{d=0}/dX^0$. Let $U= (\widehat{H^1(X)})^+$ where $\widehat{H^1(X)}$ is the inverse image of $H^1(X)$ in $X^1_{d=0}$. We have a short exact sequence of complexes
$$ 0 \lra X_0\lra X\lra X_1\lra 0$$
where $$X_0^n=\left\lbrace\begin{matrix} X^n& \text{if }  n<1,\cr
U &\text{if } n=1,\cr 0&\text{if } n >1\end{matrix}\right.\ \ \text{ and }\ \ 
X_1^n=\left\lbrace\begin{matrix} 0& \text{if }  n<1,\cr
X^1/U &\text{if } n=1,\cr X^n&\text{if } n >1\end{matrix}\right.$$

which gives the desired distinguished triangle.

We have $\cA^t=D^{\leq 0}\cap D^{\geq 0}$ and the last assertion is obvious. \qed\\

In particular, $\cA^t$ is an abelian category (\cite{KS}, prop. 10.1.11).\\

Denote $\cA^t_0$ the full subcategory of $\cA^t$ whose objects are those $A^\bullet$ such that $A^n=0$ for $n\not\in\{0,1\}$. To give an object $A^\bullet$ of $\cA^t_0$ amounts to give a morphism $$d_A=d^0_{A^\bullet}:A^0\to A^1$$ of $\cA$ such that $\ker(d_A)$ is an object of $\cA^-$ and $\coker(d_A)$ an object of $\cA^+$. 

The inclusion functor $\cA^t_0\to\cA^t$ is obviously an equivalence of categories: there is even a canonical quasi-inverse
$$\cA^t\to \cA^t_0\ :$$
it sends $A^\bullet$ to $A^{-1}/dA^{-2}\to (A^0)_{d=0}$.\\

We have an obvious functor
$$\iota^+_t : \cA^+\to \cA_0^t\ \ ,\ \ A\mapsto (0\lra A)\ .$$
It is easy to check that this functor is fully faithful and we denote $\cA^{t,-}_0$ its essential image.

Similarly, it is easy to check that the functor
$$\iota^-_t: \cA^-\to \cA^t_0\ : \ A\mapsto (A\lra 0)$$
is fully faithful and we denote $\cA_0^{t,+}$ its essential image.

It is also easy to check that $\wt t= (\cA_0^{t,+},\cA_0^{t,-})$ is a torsion pair in $\cA_0^t$. \\

\begin{prop} i) $t=(\cM^{\geq 0}(G_K),\cM^{<0}(G_K))$ is a torsion pair in $\cM(G_K)$.

ii) $t'=(\cC^{<0}(G_K),\cC^{\geq 0}(G_K))$ is a torsion pair in $\cC(G_K)$.

\end{prop}

\pr i) We already know (\S \ref{cmocgo}) that, for any object $\cF$ of $\cM(G_K)$, we have a canonical exact sequence
$$0\lra\cF^{\geq 0}\lra \cF\lra \cF^{<0}\lra 0$$
with $\cF^{\geq 0}\in\cM^{\geq 0}(G_K)$ and $\cF^{<0}\in\cM^{<0}(G_K)$.

If $f:\cF\to\cG$ is a morphism of $\cM(G_K)$, it sends $\cF^{\geq 0}$ to $\cG^{\geq 0}$. Therefore if $\cF\in\cM^{\geq 0}(G_K)$ ($\iff \cF^{\geq 0}=\cF$) and if $\cG\in\cM^{<0}(G_K)$ ($\iff \cG^{\geq 0}=0$), we have $f=0$.

ii) We already know (prop.\ref{vcoef}) that, for any object $V$ of $\cC(G_K)$, we have a canonical exact sequence
$$0\lra V^{<0}\lra V\lra V^{\geq 0}\lra 0$$
with $V^{< 0}\in\cC^{< 0}(G_K)$ and $V^{\geq 0}\in\cC^{\geq 0}(G_K)$.
Let $f:V_1\to V_2$ be a morphism of $\cC(G_K)$ with $V_1\in\cC^{<0}(G_K)$ and $V_2\in\cC^{\geq 0}(G_K)$. We can find a monomorphism $V_2\to W$ with $W\in\cC^\infty(G_K)$. As any morphism from $V_1$ to $W$ is $0$, the compositum $V_1\to V_2\to W$ is $0$, hence $f=0$.\qed\\

Denote ${\rm Ar}^t(\cM(G_K))$ the full subcategory of the categories of arrows  of $\cM^{\geq 0}(G_K)$ whose objects are those $d_\cF:\cF^0\to\cF^1$ such that $\ker d_\cF\in\cM^{<0}(G_K)$. Denote $(\cM(G_K))_{00}^t$ the full subcategory of $(\cM(G_K))_0^t$ whose objects are of the form
$$d_\cF:\cF^0\to\cF^1$$
with $\cF^0$ and $\cF^1$ objects of  $\cM^{\geq 0}(G_K)$.

As $\cM^{\geq 0}(G_K)$ is stable by taking quotients,  $(\cM(G_K))_{00}^t$  and ${\rm Ar}^t(\cM(G_K))$ have the same objects. With obvious conventions, $(\cM(G_K))_{00}^t$  is the category deduced from ${\rm Ar}^t(\cM(G_K))$ by working up to homotopies and inverting quasi-isomorphisms.

\begin{prop}
The inclusion functor
$$(\cM(G_K))_{00}^t\lra (\cM(G_K))_{0}^t$$
is an equivalence of categories.
\end{prop}

\pr  It means that any object $d_\cF:\cF^0\to\cF^1$ of $(\cM(G_K))_{0}^t$ is quasi-isomorphic to an object of $(\cM(G_K))_{00}^t$. Indeed, we may find a monomorphism $\cF^0\to\cG^0$ of $\cM(G_K)$ with $\cG^0\in\cM^{\geq 0}(G_K)$. Set 
$$\cG^1=\cG^0\oplus_{\cF^0}\cF^1$$
We have a short exact sequence
$$0\lra \ov\cG^0\lra\cG^1\lra\coker d_\cF\lra 0$$
where $\ov\cG^0$ is a quotient of $\cG^0$. Then $\coker d_\cF\in\cM^{\geq 0}(G_K)$ by assumption and $\ov\cG^0$ also because $\cM^{\geq 0}(G_K)$ is stable under taking quotients. As it is also stable under extensions, $\cG^1$ also belongs to $\cM^{\geq 0}(G_K)$. Hence, $\cG^0\to\cG^1$ is an object of $(\cM(G_K))_{00}^t$ which is quasi-isomorphic to $\cF^0\to\cF^1$. \qed\\

Similarly, denote ${\rm Ar}^{t'}(\cC(G_K))$ the full subcategory of the categories of arrows  of $\cC^{\geq 0}(G_K)$ whose objects are those $d_V:V^0\to V^1$ such that $\coker d_V\in\cC^{<0}(G_K)$. Denote $(\cC(G_K))_{00}^{t'}$ the full subcategory of $(\cC(G_K))_0^{t'}$ whose objects are of the form
$$d_V:V^0\to V^1$$
with $V^0$ and $V^1$ objects of  $\cC^{\geq 0}(G_K)$.

As   $\cC^{\geq 0}(G_K)$ is stable by taking subobjects,  $(\cC(G_K))_{00}^{t'}$  and ${\rm Ar}^{t'}(\cC(G_K))$ have the same objects. With obvious conventions, $(\cC(G_K))_{00}^{t'}$  is the category deduced from ${\rm Ar}^{t'}(\cC(G_K))$ by working up to homotopies and inverting quasi-isomorphisms.

\begin{prop}
The inclusion functor 
$$(\cC(G_K))_{00}^{t'}\lra (\cC(G_K))_{0}^{t'}$$
is an equivalence of categories.
\end{prop}

\pr The proof is entirely similar to the proof of the previous proposition: It means that any object $d_V:V^0\to V^1$ of $(\cC(G_K)_{0}^{t'}$ is quasi-isomorphic to an object of $(\cC(G_K))_{00}^{t'}$. Indeed, we may find an epimorphism $W^1\to V^1$ of $\cC(G^K)$ with $V^1\in\cC^{\geq 0}(G_K)$. Set 
$$W^0=V_0\times_{V^1}W^1$$
We have a short exact sequence
$$0\lra\ker d_V \lra W^0\lra W'\lra 0$$
where $W'$ is a subobject  of $\cG^0$. Then $\ker d_V\in\cC^{\geq 0}(G_K)$ by assumption and $W'$ also because $\cC^{\geq 0}(G_K)$ is stable under taking subobjects. As it is also stable under extensions, $W^0$ also belongs to $\cC^{\geq 0}(G_K)$. Hence, $V^0\to V^1$ is an object of $(\cC(G_K))_{00}^{t'}$ which is quasi-isomorphic to $V^0\to V^1$. \qed\\

\begin{thm} i) The functor
$$\wh\Gamma: {\rm Ar}^t(\cM(G_K))\lra \cC(G_K)\ \ ,\ \ (d_\cF:\cF^0\lra\cF^1)\mapsto \coker(\cF^0(X)\lra\cF^1(X))$$
factors uniquely through a functor
$$\Gamma : \cM(G_K)_{00}^t\lra \cC(G_K)$$
and $\Gamma$ is an equivalence of categories.

ii) The functor
$$\wh\Delta:{\rm Ar}^{t'}(\cC(G_K))\lra \cM(G_K)\ \ ,\ \ (d_V:V^0\lra V^1)\mapsto \ker(\cF_{V^0}\lra\cF_{V^1})$$
factors uniquely through a functor
$$\Delta : \cC(G_K)_{00}^{t'}\lra \cM(G_K)$$
and $\Delta$ is an equivalence of categories.
\end{thm}

\pr Let's prove  (i). Set $\wh\cM={\rm Ar}^t(\cM(G_K))$ and $\cM=\cM(G_K)^t_{00}$. If $d_\cF=\cF^0\to\cF^1$ is an object of one of these categories we denote it also $d_\cF$ or $\cF^0\to\cF^1$.  

We see that $\cM$ has an obvious structure of an exact category and that the natural functor $\wh\cM\to\cM$ is exact.

-- Let $\wh\cM^+$ (resp. $\cM^+$) the full subcategory of $\wh\cM$ (resp. $\cM$) whose objects are those $d_\cF$'s  such that $\coker d_\cF=0$. For such an object, as $\ker d_\cF\in\cM^{<0}(G_K)$ though $\cF^0$ and $\cF^1$ belong to $\cM^{\geq 0}(G_K)$, the long exact sequence of coherent cohomology associated to the exact sequence of $\cM(G_K)$
$$0\lra \ker d_\cF\lra \cF^0\lra\cF^1\lra 0$$
is reduced to 
$$0\lra \cF^0(X)\lra\cF^1(X)\lra H^1(X,\ker d_\cF)\lra 0$$
Granted to the proposition \ref{equivcoeff}, this shows that the restriction of $\wh\Gamma$ to $\wh\cM^+$ factors through a functor $$\Gamma^+:\cM^+\to\cC^{<0}(G_K)$$
which is an equivalence of categories.

-- Let $\wh\cM^-$ (resp. $\cM^-$) the full subcategory of $\wh\cM$ (resp. $\cM$) whose objects are those $d_\cF$ such that $\cF^0=0$. The natural functor $\wh\cM^-\lra\cM^-$ is an equivalence of categories and, granted to the theorem \ref{mainthm}, the restriction of $\wh\Gamma$ to $\wh\cM^+$ factors through  an equivalence of categories
$$\Gamma^-: \cM^-\lra\cC^{\geq 0}(G_K)$$

-- For any $d_\cF\in\wh M$, we have a canonical short exact sequence
$$0\lra d_{\cF_+}\lra d_\cF\lra d_{\cF_-}\lra 0$$
with $d_{\cF_+}=(\cF^0\to\im\ d_\F)\in\wh\cM^+$ and $d_{\cF_-}=(0\to\cF^1)\in\wh\cM^-$ and this construction is functorial. Moreover, we see that the sequence
$$0\lra \wh\Gamma(d_{\cF_+})\lra\wh\Gamma(d_\cF)\lra\wh\Gamma(d_{\cF_-})\lra 0$$
is exact. 

From these facts, we see that $\wh\Gamma$ factors through a functor $\Gamma:\cM\to\cC(G_K)$ and that this functor is faithful. It is also straightforward to check that it is exact. 

We are left to check the essential surjectivity: Let $V\in\cC(G_K)$. We can find a short exact sequence in 
$\cC(G_K)$
$$0\lra U\lra \wh V\lra V\lra 0$$
with $U\in\cC^0(G_K)$ and $\wh V\in\cC^{\geq 0}(G_K)$. Let $\cF^-$ be the kernel of the morphism $\cF_U\to\cF_{\wh V}$ of $\cM(G_K)$. As the functor global section is left exact, we have an exact sequence
$$0\lra \cF^-(X)\lra\cF_U(X)\lra \cF_{\wh V}(X)$$
But $\cF_U(X)=U$, $\cF_{\wh V}(X)=\wh V$ and the map $U\to\wh V$ is the given map which is injective. Therefore $\cF^-(X)=0$ which means that $\cF^-\in\cM^{<0}(G_K)$ and
$$d_\cF= (\cF_U\lra\cF_{\wh V})$$
is an object of $\cM$. Clearly $\Gamma(d_\cF)=V$, i.e. $\Gamma$ is essentially surjective.

The proof of (ii) is entirely similar and we leave it to the reader. \qed

\begin{remk} {\rm The category $\cM(G_K)^t_{00}$ is a full subcategory of $D^b(\cM(G_K))$ though $\cC(G_K)$ is a full subcategory of $D^b(\cC(G_K))$. The functor $\Gamma$ of the previous theorem is the restriction to $\cM(G_K)^t_{00}$ of the functor $\Gamma: D^b(\cM(G_K))\to D^b(\cC(G_K))$ considered in \S \ref{heart}. Similarly, the functor $\Delta$ of the previous theorem is the restriction to $\cC(G_K)^{t'}_{00}$ of the functor $\Delta:D^b(\cC(G_K))\to D^b(\cM(G_K))$ considered in \S \ref{heart}.}
\end{remk}


\

\noindent{\small Jean-Marc Fontaine, Universit\'e de Paris-Sud, Math\'ematique, Bâtiment 307, 91405 Orsay Cedex, France.}
\end{document}